\documentclass[graybox]{svmult}

\usepackage{type1cm}        
\usepackage{makeidx}         
\usepackage{graphicx}        
\usepackage{multicol}        
\usepackage[bottom]{footmisc}

\usepackage{newtxtext}       %
\usepackage[varvw]{newtxmath}       

\makeindex             

\usepackage{amsmath,amsfonts,eucal,mathtools} 
\usepackage{subfig,tikz} 
\usetikzlibrary{matrix, patterns}
\usepackage{lineno} 

\usepackage[norelsize,ruled,vlined,titlenumbered,linesnumbered]{algorithm2e}
\SetKwFor{For}{for}{}{end for}
\SetAlgoCaptionSeparator{}
\makeatletter
\newcommand{\nosemic}{\renewcommand{\@endalgocfline}{\relax}}
\newcommand{\dosemic}{\renewcommand{\@endalgocfline}{\algocf@endline}}
\let\oldnl\nl
\newcommand{\nonl}{\renewcommand{\nl}{\let\nl\oldnl}}
\makeatother

\newcommand{\defeq}{\coloneq}

\newcommand{\NN}{\mathbb{N}}
\newcommand{\PP}{\mathbb{P}}
\newcommand{\RR}{\mathbb{R}}
\let\SS\relax
\newcommand{\SS}{\mathbb{S}}
\newcommand{\cE}{\mathcal{E}}
\newcommand{\cM}{\mathcal{M}}
\newcommand{\cV}{\mathcal{V}}
\newcommand{\cW}{\mathcal{W}}
\newcommand{\cF}{\mathcal{F}}
\newcommand{\cB}{\mathcal{B}}
\newcommand{\cL}{\mathcal{L}}
\newcommand{\cR}{\mathcal{R}}

\DeclareMathOperator{\supp}{supp}
\DeclareMathOperator{\spn}{span}
\newcommand{\m}[1]{\mathbf{#1}}

\begin{document}

\title*{A Hitchhiker's Guide to LR B-Splines}
\titlerunning{A Hitchhiker's Guide to LR B-Splines}
\author{Carla Manni \and Francesco Patrizi \and Hendrik Speleers}
\authorrunning{Carla Manni \and Francesco Patrizi \and Hendrik Speleers}
\institute{Carla Manni \and Francesco Patrizi \and Hendrik Speleers\at Department of Mathematics, University of Rome Tor Vergata, Italy \\
\email{\{manni, patrizi, speleers\}@mat.uniroma2.it}
}
%
%
\maketitle

\abstract{
We review the theoretical foundations of locally refined (LR) B-splines, a natural extension of univariate B-splines to multivariate, locally refined, and generally unstructured box partitions.
LR B-splines are defined using standard tools from B-spline theory, such as knot insertion, and retain key properties including locality, non-negativity, and a piecewise polynomial structure.
Our presentation starts from the definition of spline spaces over box partitions and discusses the instability of their dimension formula due to the presence of a homological term. Conditions ensuring the stabilization of the dimension naturally motivate the construction of LR meshes and LR B-splines.
We subsequently analyze the spanning properties of LR B-splines, with particular emphasis on completeness, linear dependence, and local linear independence. We also address known characterizations and refinement strategies ensuring that the resulting LR B-splines preserve the desirable properties of tensor-product B-splines.
}

\section{Introduction}\label{sec:introduction}
Splines are piecewise functions where the pieces are smoothly joined together. Most commonly, the pieces are polynomials. To efficiently work with splines, one needs a suitable basis for their representation.
B-splines are the most prominent spline basis functions because they possess several properties that are important from both theoretical and computational perspectives. In particular, the local support, non-negativity, and partition of unity of B-splines allow for an easy control and modification of the geometries they describe. This motivates their undisputed success as modeling tools for objects with complex shapes in engineering; see, e.g., \cite{cohen,cagdbook}.
Thanks to their ability to locally reproduce polynomials, B-splines also play a central role in the construction of approximation schemes; see, e.g., \cite{deboor,schumaker}.
Furthermore, the advent of isogeometric analysis has integrated spline technologies into finite element analysis; see, e.g., \cite{igabook}. Local linear independence and high global smoothness are additional B-spline properties that are well appreciated in this context.

Tensor-product B-splines are the most common multivariate spline basis functions. Their popularity roots in their
simple, elegant, and efficient construction; they are nothing but tensor products of univariate B-splines.
The tensor structure of the underlying mesh, however, is the major weakness of tensor-product B-splines as it hinders adequate
local refinement, forcing the use of unnecessarily large discrete spaces and leading to a significant loss in efficiency. This is regarded as a major limitation in numerical simulation.

The request for adaptivity triggered the interest in more efficient spline technologies, and related  B-spline-like bases, still based on local tensor structures.
The underling meshes, the so-called \emph{T-meshes}, remain axis-aligned and preserve locally the simplicity of the tensor approach, but support local refinement. The most popular examples are (analysis-suitable) T-splines \cite{ast,sederberg}, (truncated) hierarchical B-splines \cite{hb,thb}, PHT-splines \cite{pht}, and LR B-splines \cite{lr}.

The definition of LR B-splines, short for locally refined B-splines, is inspired by the knot insertion refinement process of univariate (and tensor-product) B-splines;
it is consistent with the tensor-product B-spline definition when the underlying mesh is a tensor mesh and reduces to standard B-splines in the univariate case.
This makes them one of the most elegant extensions of univariate B-splines on local tensor structures.
Since their introduction in \cite{lr}, LR B-splines have found interesting applications in several contexts, ranging from data approximation \cite{scattereddata,n2s2,skytt} to numerical simulation \cite{hollweck,johannessen2014,zimmermann}.

In this paper, we review the theoretical foundations of LR B-splines. By collecting the most recent results, we offer a coherent and organic presentation of the state of the art of this spline technology, with a special focus on the similarities and differences between LR and tensor-product B-splines.
The presentation is accompanied and enriched by a large number of figures, illustrating the main concepts and used terminology. We mainly focus on bivariate LR B-splines, since the multivariate case is a direct extension.

The remainder of the paper is divided into four sections.
\begin{itemize}
	\item Section~\ref{sec:spline-spaces} presents the definition of spline spaces on multivariate box partitions and discusses their dimension formula, which may exhibit instability. Conditions ensuring the stabilization of the dimension formula are provided. The construction of LR meshes and LR B-splines then naturally emerges from these conditions.
	\item Section~\ref{sec:lr-definition} introduces LR meshes and LR B-splines. It reviews their definition and basic properties, as well as the core algorithms for generating them, properly scaled to yield a convex partition of unity.
	\item Section~\ref{sec:spanning-linind} discusses important aspects related to the span of LR B-splines. It analyzes when their span coincides with the full spline space on the underlying mesh and when they exhibit linear or local linear dependence. Furthermore, it outlines refinement strategies ensuring the local linear independence of the resulting LR B-splines.
	\item Section~\ref{sec:conclusion} collects some final remarks and describes some possible extensions of LR B-splines and their applications.
\end{itemize}

Although the presentation is more or less self-contained, being familiar with the definition and main properties of univariate B-splines, in particular with the knot insertion procedure, can be helpful. Comprehensive introductions can be found in \cite{foundation} and \cite{standard}, as well as in the books \cite{deboor} and \cite{schumaker}.

Throughout the paper, we use capital Latin calligraphic letters to denote sets and the cardinality of the set $\cal{A}$ is denoted by $|\cal{A}|$.

\section{Spline Spaces}\label{sec:spline-spaces}
In this section, we introduce the geometric entities underlying the definition of LR B-splines as well as the associated spline space. We also recall its dimension formula \cite{splinedim} and the corresponding dimension increasing formula \cite{lr}. The dimension formula can be decomposed into two parts: a combinatorial component, depending only on the polynomial degrees, continuity constraints, and topological quantities, and a second component given by a homological term.
The latter may lead to an unstable dimension of the spline space in the multivariate setting. We present sufficient conditions that ensure the vanishing of this homological term, thereby stabilizing the dimension formula.

We first consider the univariate case, then the bivariate case, and finally the general multivariate setting.

\subsection{Univariate Spline Space}
In the univariate case, a spline space is identified by a knot vector, defined as follows.

\begin{definition}\label{def:knotvector}
Given two integers $p \geq 0$ and $m \geq 1$, let $\pmb{\tau} \defeq (\tau_1, \ldots, \tau_m)$ be a strictly increasing sequence of real numbers and let  $\mu\colon\pmb{\tau} \to \mathbb{N}$ be a function such that $1 \leq \mu(\tau_i) \leq p+1$ for all $i$. The corresponding triple $\pmb{\tau}_p^{\mu} \defeq (\pmb{\tau}, \mu, p)$ is called a \textbf{global knot vector}, or simply a \textbf{knot vector}.
Moreover, we use the following terminology:
\begin{itemize}
\item the elements of $\pmb{\tau}$ are the \textbf{breakpoints}; 
\item $\mu(\tau_i)$
is the \textbf{multiplicity} of $\tau_i$; a breakpoint $\tau_i$ has \textbf{full multiplicity} if $\mu(\tau_i) = p+1$, while it is \textbf{simple} if $\mu(\tau_i) = 1$;
\item  $\pmb{\tau}_p^{\mu}$ is \textbf{open} (or \textbf{clamped}) if both $\tau_1$ and $\tau_m$ have full multiplicity;
\item $[\tau_1, \tau_m]$ is the \textbf{knot span}, and the intervals $[\tau_i, \tau_{i+1}]$, $i = 1, \ldots, m-1$, are the \textbf{knot intervals}.
\end{itemize}
\end{definition}

Given a non-negative integer $p$, let $\PP_p$ be the vector space of polynomials of degree at most $p$.
The univariate spline space is then defined as follows.

\begin{definition}\label{def:splinespace1d}
Given a knot vector $\pmb{\tau}_p^{\mu} \defeq (\pmb{\tau}, \mu, p)$, with $\pmb{\tau} \defeq (\tau_1, \ldots, \tau_m)$, the \textbf{univariate spline space on} $\pmb{\tau}_p^\mu$, denoted by $\SS[\pmb{\tau}_p^\mu]$,
is the set of functions $f\colon\mathbb{R} \to \mathbb{R}$ such that
\begin{itemize}
    \item $f$ is zero outside $[\tau_1, \tau_m]$;
    \item the restriction of $f$ to  $[\tau_i, \tau_{i+1})$,  $i = 1, \ldots, m-2$, and to $[\tau_{m-1}, \tau_m]$ belongs to $\PP_p$;
    \item $f$ is $C^{p-\mu(\tau_i)}$-continuous at $\tau_i$, $i=1,\dots,m$.
\end{itemize}
\end{definition}

It is not difficult to prove that
\begin{equation}\label{eq:splinedim1d}
	\dim \SS[\pmb{\tau}_p^\mu] = \max\left\{  - (p+1)+ \sum_{i=1}^m \mu(\tau_i), \, 0 \right\};
\end{equation}
see, e.g., \cite[Theorem~4.4]{schumaker}.
This formula shows that sufficiently many breakpoints, with sufficiently large multiplicity, are required to obtain a non-trivial spline space. In particular, if $\pmb{\tau}_p^{\mu}$ is  an open knot vector then $\dim \SS[\pmb{\tau}_p^\mu] \geq p+1$.

\begin{remark}\label{rmk:homology1d}
With a simple manipulation, the dimension formula in \eqref{eq:splinedim1d} can be rewritten as
\begin{equation}\label{eq:splinedim1d-homological}
\dim\SS[\pmb{\tau}_p^\mu] =(m-1) (p+1) -\sum_{i=1}^m(p- \mu(\tau_i)+1) + \text{homological term},
\end{equation}
where the homological term is
\begin{equation*}
\begin{cases}
	\ 0, & \text{if }\displaystyle \sum_{i=1}^m \mu(\tau_i) > p + 1,\\
	\displaystyle\ p + 1 - \sum_{i=1}^m \mu(\tau_i), &  \text{otherwise};
\end{cases}
\end{equation*}
see \cite[Lemma~15]{splinedim}.
The expression in \eqref{eq:splinedim1d-homological} is the analog of the dimension formulas for spline spaces in higher dimensions that will be addressed in the subsequent sections.
\end{remark}

\subsection{Bivariate Spline Space} \label{sec:space-2d}
The key geometric objects for the definition of spline spaces in the bivariate setting are box partitions and the associated meshes, defined as follows.

\begin{definition}
Given an axis-aligned rectangular domain $\Omega \subset \RR^2$, a \textbf{box partition} of $\Omega$ is a finite collection $\mathcal{E}$ of axis-aligned rectangles in $\Omega$ such that
\begin{itemize}
\item $\mathring{\beta}_1 \cap \mathring{\beta}_2 = \varnothing$ for any $\beta_1, \beta_2 \in \mathcal{E}$ and $\beta_1 \neq \beta_2$, i.e., they have disjoint interiors;
\item $\bigcup_{\beta \in \mathcal{E}} \beta = \Omega$, i.e., they cover the whole domain.
\end{itemize}
The elements of $\mathcal{E}$ are called \textbf{boxes}, \textbf{cells}, or \textbf{elements}.
\end{definition}

\begin{figure}[t!]
\centering
\subfloat[]{
\begin{tikzpicture}[scale=2.25]
\fill[red!60] (0,.66) -- (.5,.66) -- (.5, 1) -- (0,1) -- cycle;
\fill[cyan!50] (0,0) -- (.5,0) -- (.5, .66) -- (0,.66) -- cycle;
\fill[orange!60] (.5, 0) -- (1,0) -- (1,.33) -- (.5, .33) --cycle;
\fill[yellow!60] (.5, .33) -- (1,.33) -- (1,.66) -- (.5, .66) --cycle;
\fill[green!60!black] (.5, .66) -- (1.5,.66) -- (1.5,1) -- (.5, 1) --cycle;
\fill[purple!60] (1, 0) -- (1.5,0) -- (1.5,.66) -- (1, .66) --cycle;
\draw (0,0) -- (1.5,0) -- (1.5, 1) -- (0,1) --cycle;
\draw (0,.66) --(1.5, .66);
\draw (.5, .33) -- (1, .33);
\draw (.5, 0) -- (.5, 1);
\draw (1, 0) -- (1, .66);
\draw[white] (0,-.08);
\end{tikzpicture}}\quad
\subfloat[]{
\begin{tikzpicture}[scale=2.25]
\tikzset{%
  block/.style    = {draw=black, fill=white, rectangle, minimum height = .05cm,
    minimum width = .05cm}
    }
\draw (0,0) -- (1.5,0) -- (1.5, 1) -- (0,1) --cycle;
\draw (0,.66) --(1.5, .66);
\draw (.5, .33) -- (1, .33);
\draw (.5, 0) -- (.5, 1);
\draw (1, 0) -- (1, .66);
\draw (0,.33) node[block]{\footnotesize{$1$}};
\draw (0,.845) node[block]{\footnotesize{$1$}};
\draw (.25,0) node[block]{\footnotesize{$1$}};
\draw (.25,.66) node[block]{\footnotesize{$1$}};
\draw (.25,1) node[block]{\footnotesize{$1$}};
\draw (.5,.165) node[block]{\footnotesize{$1$}};
\draw (.5, .495) node[block]{\footnotesize{$3$}};
\draw (.5,.845) node[block]{\footnotesize{$1$}};
\draw (.75,.66) node[block]{\footnotesize{$2$}};
\draw (.75,.33) node[block]{\footnotesize{$4$}};
\draw (.75,0) node[block]{\footnotesize{$2$}};
\draw (1.25,.66) node[block]{\footnotesize{$2$}};
\draw (1,1) node[block]{\footnotesize{$1$}};
\draw (1,.165) node[block]{\footnotesize{$1$}};
\draw (1,.495) node[block]{\footnotesize{$1$}};
\draw (1.5,.33) node[block]{\footnotesize{$1$}};
\draw (1.5,.845) node[block]{\footnotesize{$1$}};
\draw (1.25,0) node[block]{\footnotesize{$2$}};
\end{tikzpicture}
}\quad
\subfloat[]{
\begin{tikzpicture}[scale = 2.25]
\draw[ultra thick, red!60] (0,0) -- (1.5,0);
\draw[ultra thick, blue!60] (1.5, 0) -- (1.5, 1);
\draw[ultra thick, purple!60] (1.5, 1) -- (0,1);
\draw[ultra thick, magenta!80] (0, 1) -- (0, 0);
\draw[ultra thick, cyan!50] (0,.66) --(1.5, .66);
\draw[ultra thick, orange!60] (.5, .33) -- (1, .33);
\draw[ultra thick, yellow!60] (.5, 0) -- (.5, 1);
\draw[ultra thick, green!50!black] (1, 0) -- (1, .66);
\filldraw[draw = black, fill = red!80] (0.5, 0.33) circle (0.0333);
\filldraw[draw = black, fill = red!80] (1, 0.33) circle (0.0333);
\filldraw[draw = black, fill = red!80] (1, 0.66) circle (0.0333);
\draw[white] (0,-.08);
\end{tikzpicture}
}
\caption{A box partition (a) and an associated $\mu$-extended mesh (b). The minimal meshlines are indicated by squares reporting their assigned multiplicities. The distinct meshlines (or splits) of the mesh are identified in (c) by different colors and the T-vertices by dot markers. The mesh does not have constant meshlines.}\label{fig:boxpartition}
\end{figure}
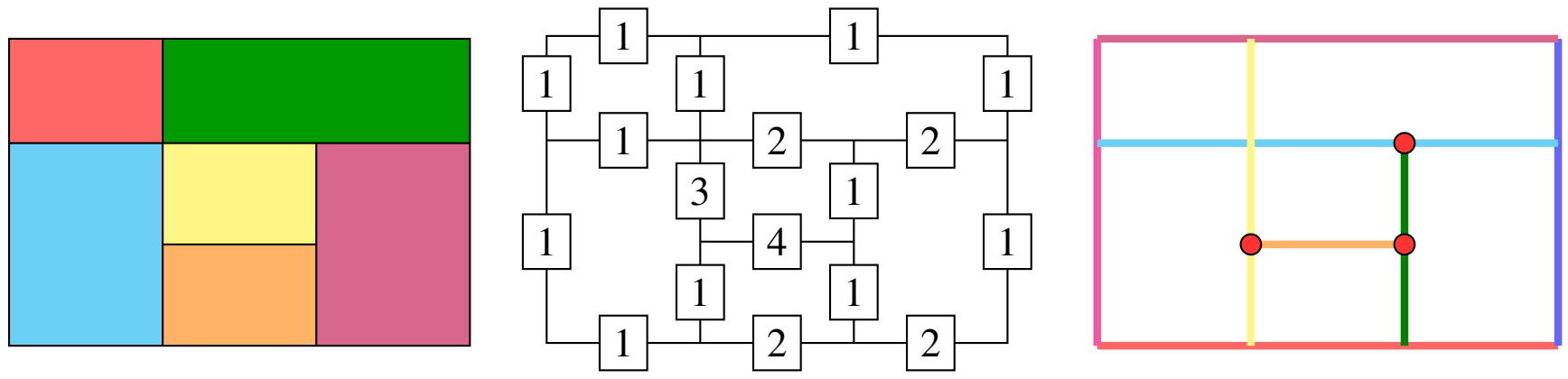

Figure~\ref{fig:boxpartition}(a) depicts an example of a box partition $\cE$.

\begin{definition}
Given a box partition $\mathcal{E}$, we define the \textbf{vertices} of $\mathcal{E}$ as the vertices of its elements. A vertex of $\mathcal{E}$ is called a \textbf{T-vertex} if it is the intersection of edges of three elements. The set of vertices of $\mathcal{E}$ is denoted by $\mathcal{V}$.
\end{definition}

\begin{definition}\label{def:boxmesh2D}
Given a box partition $\mathcal{E}$ of a rectangle $\Omega \subset \RR^2$, a \textbf{minimal meshline} of $\cE$ is a segment contained in an edge of an element of $\mathcal{E}$, connecting two, and only two, vertices in $\mathcal{V}$ at its endpoints. The collection of all minimal meshlines of the box partition is the \textbf{box mesh}, or simply the \textbf{mesh}, denoted by $\cM$, associated with $\mathcal{E}$. Moreover, we use the following terminology:
\begin{itemize}
\item a \textbf{meshline}, or \textbf{split}, of $\cM$ is a maximal union of contiguous and collinear minimal meshlines of $\cE$;
\item a \textbf{multiplicity function} $\mu\colon \cM \to \NN$ is a function that assigns a positive integer to every minimal meshline in $\cM$;
\item a \textbf{$\pmb{\mu}$-extended box mesh}, or simply a \textbf{$\pmb{\mu}$-extended mesh}, denoted by $\cM^\mu$, is a mesh $\cM$ equipped with a multiplicity function $\mu$;
\item a $\mu$-extended mesh $\cM^\mu$ has \textbf{constant} (or \textbf{uniform}) \textbf{meshlines}, or \textbf{constant} (or \textbf{uniform}) \textbf{splits}, if each meshline $\gamma$ in $\cM$ consists of minimal meshlines all having the same multiplicity. In this case, the multiplicity of a meshline $\gamma$ in $\cM$ is defined as the common multiplicity of its constituent minimal meshlines.
\end{itemize}
\end{definition}

A meshline can be expressed as the Cartesian product of a point in $\RR$ and a finite interval, namely $\{c\} \times [a,b]$ or $[a,b] \times \{c\}$, corresponding to vertical and horizontal meshlines, respectively, with $a,b,c \in \RR$. Let $k \in \{1,2\}$ denote the position of point $\{c\}$ in the Cartesian product. We sometimes refer to it as a \textbf{$k$-meshline} to stress the meshline's direction, namely vertical if $k=1$ and horizontal if $k = 2$. Occasionally, we write \textbf{$(k, c)$-meshline} to specify the coordinate value in the $k$th direction as well.
\begin{definition}
The set of minimal $k$-meshlines in a box mesh $\cM$ is denoted by $\cM_{[k]}$, $k=1,2.$
\end{definition}

Figure~\ref{fig:boxpartition}(b) depicts an example of a $\mu$-extended mesh associated with the box partition in Figure~\ref{fig:boxpartition}(a).
A tensor mesh is a very special instance of a $\mu$-extended mesh.
\begin{definition}\label{def:tensor-mesh}
A $\mu$-extended mesh $\cM^\mu$ is a \textbf{tensor mesh} if the set of vertices $\cV$ contains no T-vertices and $\cM^\mu$ has constant meshlines.
\end{definition}

In order to define a spline space on a given mesh, we need to relate the mesh to a (bi)degree.
\begin{definition}\label{def:spline-mesh}
A \textbf{spline mesh} in $\RR^2$ is a pair $\cM_{\m{p}}^{\mu} \defeq (\cM^\mu, \m{p})$, where $\m{p} \defeq (p_1,p_2)$ is a pair of integers, with $p_1, p_2 \geq 0$, and $\cM^\mu$ is a $\mu$-extended mesh arising from a box partition $\cE$, with multiplicities satisfying $1 \leq \mu(\gamma) \leq p_k+1$ for every minimal $k$-meshline $\gamma \in \cM$.
Moreover, we use the following terminology:
\begin{itemize}
\item  a minimal $k$-meshline $\gamma$ has \textbf{full multiplicity} if $\mu(\gamma) = p_k+1$, while it is \textbf{simple} if $\mu(\gamma) = 1$;  
\item  $\cM_\m{p}^\mu$ is \textbf{open} if every boundary minimal meshline has full multiplicity;
\item $\cM_\m{p}^\mu$ is a \textbf{tensor spline mesh} if $\cM^\mu$ is a tensor mesh.
\end{itemize}
\end{definition}

Given a bidegree $\m{p} \defeq (p_1, p_2)$, let $\PP_\m{p}$ be the bivariate polynomial space of degree at most $p_k$ for $k=1,2$ componentwise, i.e., the space spanned by the monomials $x_1^{i_1} x_2^{i_2}$ such that $0 \leq i_k \leq p_k$ for $k = 1, 2$. We are now ready to define the spline space in the bivariate setting.

\begin{definition}
Given a spline mesh $\cM_{\m{p}}^\mu$ corresponding to a box partition $\cE$ of $\Omega \defeq [a_1,b_1]\times[a_2,b_2]$, for any element $\beta \in \cE$, with
$\beta = I_1 \times I_2$ and $I_k = [a_{\beta,k}, b_{\beta,k}]$, we set
\begin{equation*}
\tilde{\beta} \defeq \tilde{I}_1 \times \tilde{I}_2 \quad \text{with} \quad
\tilde{I}_k \defeq
\begin{cases}
[a_{\beta,k}, b_{\beta,k}), & \text{if } b_{\beta,k} < b_k, \\
[a_{\beta,k}, b_{\beta,k}], & \text{if } b_{\beta,k} = b_k.
\end{cases}
\end{equation*}
The \textbf{bivariate spline space on} $\cM_\m{p}^\mu$, denoted by $\SS[\cM_{\m{p}}^\mu]$, is the set of functions $f\colon\RR^2\to \RR$ such that
\begin{itemize}
\item $f$ is zero outside $\Omega$;
\item for each element $\beta \in \cE$, the restriction $f|_{\tilde{\beta}}$ belongs to $\PP_\m{p}$;
\item for each minimal $k$-meshline $\gamma \in \cM$, $f$ is $C^{p_k-\mu(\gamma)}$-continuous across $\gamma$, i.e., with respect to the $k$th direction.
\end{itemize}
\end{definition}

In order to state a dimension formula for $\SS[\cM_{\m{p}}^\mu]$, we need to assign multiplicities to the vertices of the partition as well.
\begin{definition}\label{def:mu1mu2}
Let $\cM^\mu$ be a $\mu$-extended mesh  corresponding to a box partition $\cE$. For any vertex $\m{v} \in \cV$, we define for $k = 1, 2$,
\begin{equation*}
	\mu_k(\m{v}) \defeq \max\left\{\,\mu(\gamma): \m{v} \in \gamma \text{ and } \gamma \text{ is a minimal } k\text{-meshline in } \cM\,\right\}.
\end{equation*}
\end{definition}

\begin{figure}[t!]
	\centering
	\begin{tikzpicture}[scale=3]
		\tikzset{%
			block/.style    = {draw=black, fill=white, rectangle, minimum height = .05cm,
				minimum width = .05cm},
		}
		\draw (0,0) -- (1.5, 0) -- (1.5, 1) -- (0,1) -- cycle;
		\draw (0,.5) -- (1.5,.5);
		\draw (.5, 0) -- (.5, 1);
		\draw (1, .5) -- (1, 1);
		\draw (0,.25) node[block]{\footnotesize{$1$}};
		\draw (0,.75) node[block]{\footnotesize{$1$}};
		\draw (.25,0) node[block]{\footnotesize{$1$}};
		\draw (.25,.5) node[block]{\footnotesize{$1$}};
		\draw (.25,1) node[block]{\footnotesize{$1$}};
		\draw (.5,.25) node[block]{\footnotesize{$1$}};
		\draw (.5,.75) node[block]{\footnotesize{$1$}};
		\draw (.75,.5) node[block]{\footnotesize{$2$}};
		\draw (1.25,.5) node[block]{\footnotesize{$2$}};
		\draw (.75,1) node[block]{\footnotesize{$1$}};
		\draw (1.25,1) node[block]{\footnotesize{$1$}};
		\draw (1,0) node[block]{\footnotesize{$1$}};
		\draw (1,.75) node[block]{\footnotesize{$1$}};
		\draw (1.5,.25) node[block]{\footnotesize{$1$}};
		\draw (1.5,.75) node[block]{\footnotesize{$1$}};
		\fill (.5, .5) node[below]{\footnotesize{$\quad~\m{v}_1$}} circle (.02);
		\fill (1, .5) node[below]{\footnotesize{$\m{v}_2$}} circle (.02);
	\end{tikzpicture}
	\caption{A $\mu$-extended mesh with two internal vertices $\m{v}_1$ and $\m{v}_2$. According to Definition~\ref{def:mu1mu2}, the multiplicities of these vertices are $\mu_1(\m{v}_i) = 1$, $\mu_2(\m{v}_i) = 2$, for both $i \in \{1, 2\}$.}\label{fig:exmu1mu2}
\end{figure}
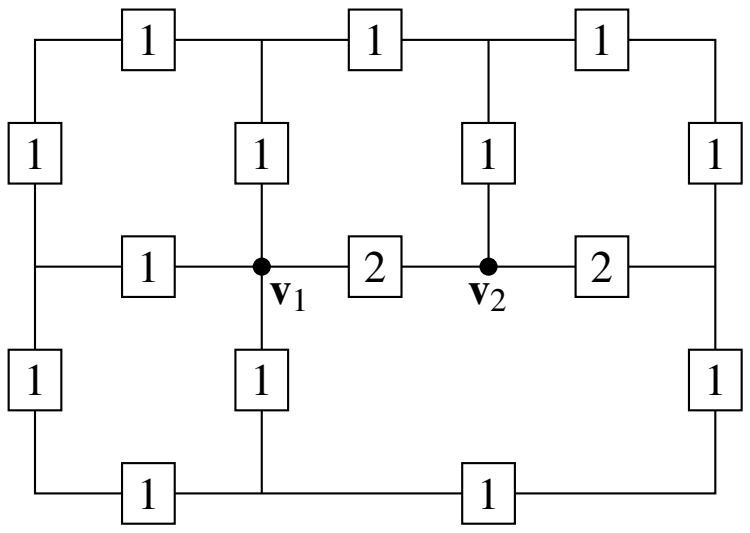

An example of the multiplicities of two vertices in a given $\mu$-extended mesh can be found in Figure~\ref{fig:exmu1mu2}.
The following result has been shown in {\cite[Theorem~2]{splinedim}}.
\begin{theorem}
Let $\cM_\m{p}^\mu$ be a spline mesh corresponding to a box partition $\mathcal{E}$.
Then,
\begin{equation}\label{eq:splinedim2d}
	\begin{aligned}
		&\dim\,\SS[\cM_\m{p}^\mu] \\
		&\quad = |\cE|(p_1+1)(p_2+1) + \sum_{\m{v}\in\mathcal{V}} (p_1-\mu_1(\m{v}) +1)(p_2-\mu_2(\m{v})+1)\\
		&\quad\quad -(p_2+1)\sum_{\gamma \in \cM_{[1]}} (p_1-\mu(\gamma)+1)  -(p_1+1)\sum_{\gamma \in \cM_{[2]}}(p_2-\mu(\gamma)+1) \\
		&\quad\quad + \text{homological term}.
	\end{aligned}
\end{equation}
\end{theorem}

The dimension formula in \eqref{eq:splinedim2d} consists of two parts: a combinatorial component and a homological term.
The latter may lead to an unstable dimension, as it may depend on the sizes of the boxes in $\cE$, rather than solely on their topological arrangement. As a consequence, spline spaces defined on meshes with the same topological structure, the same minimal meshline multiplicities, and the same degrees may have different dimensions, as shown in the following example.

\begin{example}{Example: Unstable Dimension}%
We consider two open spline meshes with $\m{p} = (2,2)$ and meshes as in Figure~\ref{fig:dimunstable}(a)--(b), with internal meshlines of multiplicity $1$. It is easy to check that the combinatorial part of the dimension formula is equal to $36$ in both cases. In contrast, it can be shown that the homological part is equal to $1$ in (a) and $0$ in (b); see \cite[Example~1]{chendim}. Therefore, the corresponding spline spaces have dimension $37$ in (a) and $36$ in (b).

\begin{figure}[t!]
\centering
\subfloat[]{
\begin{tikzpicture}[scale=3]
\draw node[below]{0}(0,0);
\draw node[left]{0}(0,0);
\draw (.2,0)node[below]{\color{red}{1}};
\draw (.4,0)node[below]{2};
\draw (.6,0)node[below]{3};
\draw(.8,0) node[below]{4};
\draw (1,0)node[below]{5};
\draw (1.2,0)node[below]{6};
\draw (0,.2)node[left]{1};
\draw (0,.4)node[left]{2};
\draw (0,.6)node[left]{3};
\draw (0,.8)node[left]{4};
\draw (0,1)node[left]{5};
\draw (0,1.2)node[left]{6};
\draw (0,0) -- (1.2,0) -- (1.2,1.2) --(0,1.2) -- cycle;
\draw (0,.8) -- (1.2,.8);
\draw (0,1) -- (1.2,1);
\draw (.4,.4) -- (1,.4);
\draw (.4,0) -- (.4,1.2);
\draw (.6,.4) -- (.6,1);
\draw (.8,.2) -- (.8,.8);
\draw (1,0) -- (1,1.2);
\draw (0,.2)--(1.2,.2);
\draw[red] (.2,0) -- (.2,1.2);
\draw (.2,.6) -- (.8,.6);
\end{tikzpicture}
}\qquad
\subfloat[]{
\begin{tikzpicture}[scale=3]
\draw node[below]{0}(0,0);
\draw node[left]{0}(0,0);
\draw (.28,0)node[below]{\color{red}{1.5}};
\draw (.4,0)node[below]{2};
\draw (.6,0)node[below]{3};
\draw(.8,0) node[below]{4};
\draw (1,0)node[below]{5};
\draw (1.2,0)node[below]{6};
\draw (0,.2)node[left]{1};
\draw (0,.4)node[left]{2};
\draw (0,.6)node[left]{3};
\draw (0,.8)node[left]{4};
\draw (0,1)node[left]{5};
\draw (0,1.2)node[left]{6};
\draw (0,0) -- (1.2,0) -- (1.2,1.2) --(0,1.2) -- cycle;
\draw (0,.8) -- (1.2,.8);
\draw (0,1) -- (1.2,1);
\draw (.4,.4) -- (1,.4);
\draw (.4,0) -- (.4,1.2);
\draw (.6,.4) -- (.6,1);
\draw (.8,.2) -- (.8,.8);
\draw (1,0) -- (1,1.2);
\draw (0,.2)--(1.2,.2);
\draw[red] (.3,0) -- (.3,1.2);
\draw (.3,.6) -- (.8,.6);
\end{tikzpicture}
}
\caption{Two meshes with the same topology but different parametrization.}\label{fig:dimunstable}
\end{figure}
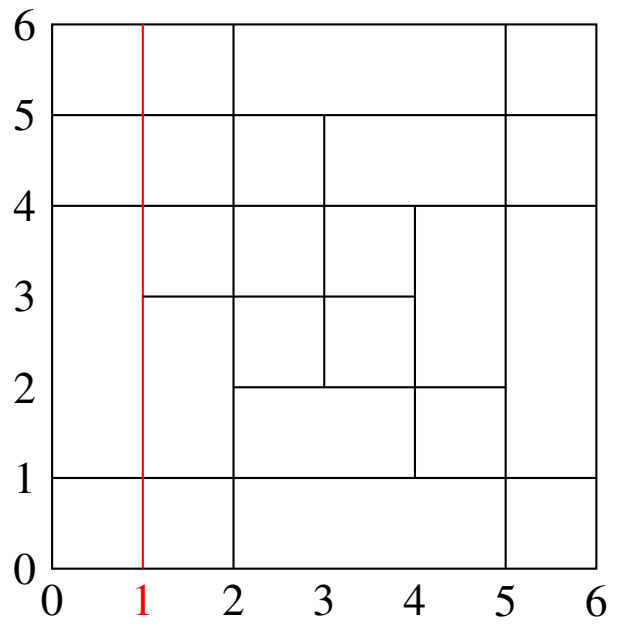
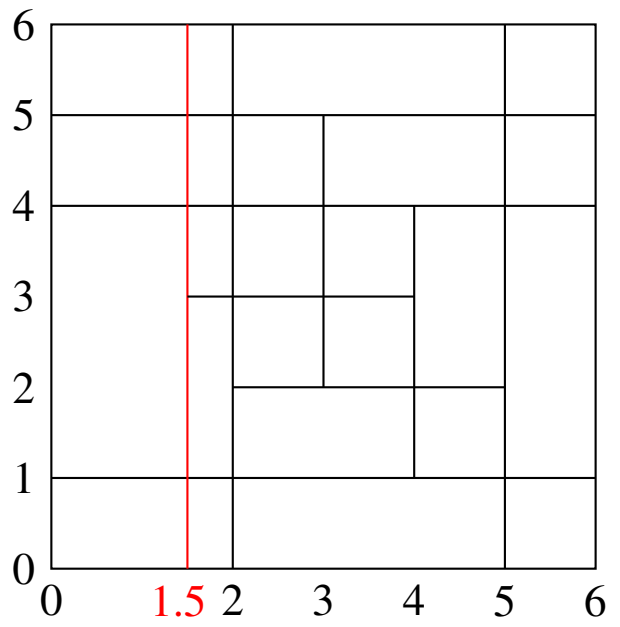
\end{example}
\vspace*{-0.3cm} 
\begin{example}{Example: Tensor Meshes}%
Let $\cM_\m{p}^\mu$  be a tensor spline mesh. According to Definitions~\ref{def:tensor-mesh} and \ref{def:spline-mesh},  
there exist two knot vectors $\pmb{\tau}_{k,p_k}^{\mu_k}$
consisting of $m_k$  breakpoints,  $k=1,2$, such that
\begin{equation*}
	\begin{split}
&\cE=\left\{[\tau_{1,i},\tau_{1,i+1}]\times[\tau_{2,j},\tau_{2,j+1}],\ i=1,\dots, m_1-1,\
j=1,\dots,m_2-1\right\},
\\
&\mathcal{V}=\left\{(\tau_{1,i},\tau_{2,j}),\ i=1,\dots, m_1,\
j=1,\dots,m_2\right\},
\\
&\mu_1(\m{v})=\mu_1(\tau_{1,i}), \quad \mu_2(\m{v})=\mu_2(\tau_{2,j}),\quad \m{v} \in\cV;
\end{split}
\end{equation*}
and the collections of minimal $k$-meshlines are given by
\begin{equation*}
\begin{split}
\cM_{[1]}&= \left\{\tau_{1,i}\times[\tau_{2,j},\tau_{2,j+1}],\ j=1,\dots,m_2-1, \ i=1,\dots,m_1\right\}, \\
\cM_{[2]}&= \left\{[\tau_{1,i},\tau_{1,i+1}]\times \tau_{2,j},\ i=1,\dots,m_1-1, \ j=1,\dots,m_2\right\}.
\end{split}
\end{equation*}
Thus, the sum of the combinatorial terms in the right-hand side of \eqref{eq:splinedim2d} reads
\begin{equation*}
	\begin{split}
&(m_1-1)(m_2-1)(p_1+1)(p_2+1) \\
&\quad +\sum_{i=1}^{m_1}\sum_{j=2}^{m_2} (p_1-\mu_1(\tau_{1,i}) +1)(p_2-\mu_2(\tau_{2,j})+1)
\\
&\quad - (p_2+1)(m_2-1)\sum_{i=1}^{m_1} (p_1-\mu_1(\tau_{1,i})+1) \\
&\quad -(p_1+1)(m_1-1)\sum_{j=1}^{m_2}(p_2-\mu_2(\tau_{2,j})+1).
	\end{split}
\end{equation*}
From \eqref{eq:splinedim1d-homological}, we see that this agrees with
$\dim\,\SS[\pmb{\tau}_{1,p_1}^{\mu_1}]\dim\,\SS[\pmb{\tau}_{2,p_2}^{\mu_2}]$ if the univariate spline spaces $\SS[\pmb{\tau}_{k,p_k}^{\mu_k}]$, $k=1,2$ are non-trivial.
Since it is well known that
\begin{equation}\label{eq:dimtensorproduct}
\dim\,\SS[\cM_\m{p}^\mu]=\dim\,\SS[\pmb{\tau}_{1,p_1}^{\mu_1}]\dim\,\SS[\pmb{\tau}_{2,p_2}^{\mu_2}],
\end{equation}
we deduce that the homological term in the right-hand side of \eqref{eq:splinedim2d} vanishes for a tensor spline mesh identifying a non-trivial spline space.
\end{example}

Since the focus of this survey is LR spline theory, we do not explore the underlying homological framework in full generality. Instead, we recall sufficient conditions on the spline mesh that guarantee the vanishing of this term and thereby restore stability. These conditions restrict our attention to a subclass of spline meshes.
Roughly speaking, the homological term vanishes if the spline mesh $\cM_\m{p}^\mu$ is constructed as an \emph{iterative refinement} of a tensor spline mesh, with a non-trivial associated spline space, by inserting \emph{sufficiently long} meshlines.

We start by describing how a spline mesh can be refined.

\begin{definition}
Given a box partition $\mathcal{E}$ and an axis-aligned segment $\gamma$, we say that $\gamma$ \textbf{traverses} $\beta \in \mathcal{E}$ if the interior of $\beta$ is split into two parts by $\gamma$, i.e., $\beta\setminus \gamma$ is not connected. 
\end{definition}
\begin{definition}
Given a box partition $\mathcal{E}$ and an axis-aligned segment $\gamma$, we say that $\gamma$ can be \textbf{inserted} into $\mathcal{E}$ if its endpoints lie on the mesh $\cM$ associated with $\mathcal{E}$.
\end{definition}

\begin{figure}[t!]
	\centering
	\subfloat[]{
		\begin{tikzpicture}[scale = 2.25]
			\tikzset{%
				block/.style    = {draw=black, fill=white, rectangle, minimum height = .05cm,
					minimum width = .05cm},
			}
			\draw[step = 0.5] (0, 0) grid (1.5, 1.5);
			\draw (0.25, 0) -- (0.25, 1);
			\draw (0.25, 0.25) node[block]{$3$};
			\draw (0.25, 0.75) node[block]{$1$};
			\draw (0, 1.25) --node[block]{$1$} (0.5, 1.25);
			\draw (1, 1.25) --node[block]{$2$} (1.5, 1.25);
			\draw[white] (0, -0.175);
		\end{tikzpicture}
	}
	\subfloat[]{
		\begin{tikzpicture}[scale = 2.25]
			\draw[step = 0.5] (0, 0) grid (1.5, 1.5);
			\draw (0.25, 0) -- (0.25, 1);
			\draw[dashed, thick, xshift = 0.5] (0.25, 0) -- (0.25, 1);
			\draw[dashed, thick, xshift = - 0.5] (0.25, 0) -- (0.25, 0.5);
			\draw[xshift = 3.75] (0.25, 0) node[below]{$\gamma_2$};
			\draw[xshift = - 3] (0.25, 0) node[below]{$\gamma_4$};
			\draw (0, 1.25) -- (0.5, 1.25);
			\draw (1, 1.25) -- (1.5, 1.25);
			\draw[dashed, thick, yshift = 0.5] (0, 1.25)node[left]{$\gamma_1$} -- (1.5, 1.25);
			\draw[dashed, thick] (0.5, 0.75) -- (1.5, 0.75)node[right]{$\gamma_3$};
			\draw (0.75, 1.5)node[above]{\scalebox{0.8}{$\mu_{\gamma_1} = 1, \mu_{\gamma_2} = 2, \mu_{\gamma_3} = 3, \mu_{\gamma_4} = 1$}};  
		\end{tikzpicture}
	}
	\subfloat[]{
		\begin{tikzpicture}[scale = 2.25]
			\tikzset{%
				block/.style    = {draw=black, fill=white, rectangle, minimum height = .05cm,
					minimum width = .05cm},
			}
			\draw[step = 0.5] (0, 0) grid (1.5, 1.5);
			\draw (0.25, 0) -- (0.25, 1);
			\draw (0.25, 0.25) node[block]{$6$};
			\draw (0.25, 0.75) node[block]{$3$};
			\draw (0, 1.25) --node[block]{$2$} (0.5, 1.25);
			\draw (1, 1.25) --node[block]{$3$} (1.5, 1.25);
			\draw (0.5, 1.25) --node[block]{$1$} (1, 1.25);
			\draw (0.5, 0.75) --node[block]{$3$} (1, 0.75);
			\draw (1, 0.75) --node[block]{$3$} (1.5, 0.75);
			\draw[white] (0, -0.175);
		\end{tikzpicture}
	}
	\caption{Starting from the initial $\mu$-extended mesh in (a), with some multiplicities of its meshlines reported, we insert the four segments, indicated with dashed lines, with their associated multiplicities as in (b). This results in the $\mu$-extended mesh in (c).}\label{fig:meshlineinsertion}
\end{figure}
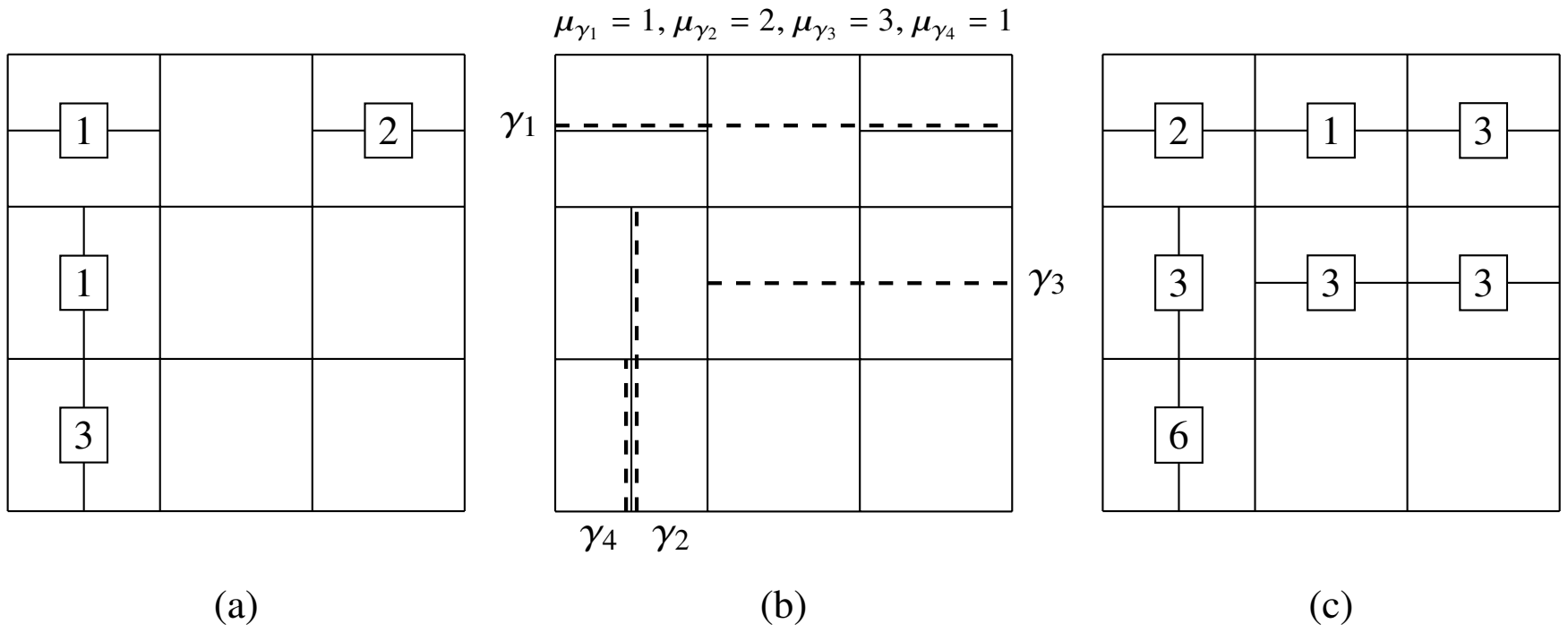
When an axis-aligned segment $\gamma$ can be inserted into the box partition $\mathcal{E}$, it either overlaps minimal meshlines of the associated mesh $\cM$, or it traverses some boxes of $\cE$, or both. Then, after insertion, each traversed $\beta \in \mathcal{E}$ is replaced by the two subboxes $\beta_1, \beta_2$, defined as the closures of the connected components of $\beta \setminus \gamma$. The resulting box partition is denoted by $\mathcal{E}+\gamma$ and its associated mesh by $\cM+\gamma$.

Let $\mu_\gamma$ be a positive integer associated with $\gamma$, i.e., its multiplicity. Then, the multiplicities of the minimal meshlines in $\cM\cap (\cM+\gamma)$ that are not contained in $\gamma$ remain unchanged, while the multiplicities of those contained in $\gamma$ are increased by $\mu_\gamma$. The new minimal meshlines in $(\cM+\gamma)\setminus \cM$ are assigned multiplicity $\mu_\gamma$. If $\mu$ denotes the multiplicity function associated with $\cM$, then the multiplicity function on the refined mesh $\cM + \gamma$ is denoted by $\mu+\mu_\gamma$. Figure~\ref{fig:meshlineinsertion} illustrates the changes in the box partition, mesh, and multiplicities after the insertion of four segments.

An interesting class of spline meshes is identified by the following definition.
\begin{definition}\label{def:iterativesplinemesh}
A spline mesh $\cM_\m{p}^\mu$ is an \textbf{iterative spline mesh} if there exists an integer $N \geq 0$ such that $\cM_\m{p}^\mu$ is the $N$th element of a sequence of spline meshes $\cM_{i,\m{p}}^{\mu_i} \defeq (\cM_i^{\mu_i}, \m{p})$ satisfying
\begin{equation*}
\cM_i \defeq \cM_{i-1} + \gamma_i, \quad
\mu_i \defeq \mu_{i-1} + \mu_{\gamma_i},
\end{equation*}
where $\{\gamma_i\}_{i=1}^N$ is a sequence of axis-aligned segments with associated multiplicities $\{\mu_{\gamma_i}\}_{i=1}^N$, inserted into an initial tensor spline mesh $\cM_{0,\m{p}}^{\mu_0}$.
\end{definition}

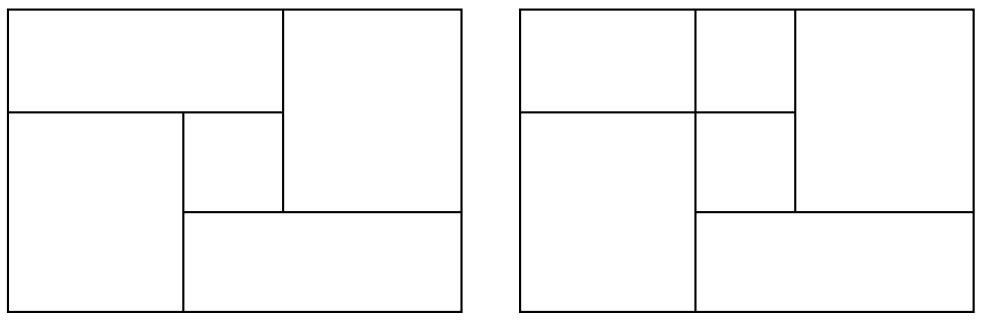
\begin{figure}[t!]
	\centering
	\subfloat[]{
		\begin{tikzpicture}[scale=2]
			\draw (-.25,0) -- (1.25,0) -- (1.25,1) -- (-.25,1) -- cycle;
			\draw (.66,1) -- (.66,.33);
			\draw (-.25,.66) -- (.66,.66);
			\draw (.33,.66) -- (.33,0);
			\draw (.33,.33) -- (1.25,.33);
	\end{tikzpicture}}\quad
	\subfloat[]{
		\begin{tikzpicture}[scale=2]
			\draw (-.25,0) -- (1.25,0) -- (1.25,1) -- (-.25,1) -- cycle;
			\draw (.66,1) -- (.66,.33);
			\draw (-.25,.66) -- (.66,.66);
			\draw (.33,1) -- (.33,0);
			\draw (.33,.33) -- (1.25,.33);
	\end{tikzpicture}}
	\caption{A spline mesh that is not iterative (a) and a modification of it that can be constructed iteratively (b).}\label{fig:noiterativemesh}
\end{figure}

Figure~\ref{fig:noiterativemesh}(a) depicts a spline mesh that is not iterative. Indeed, there is no sequence of segment insertions that can generate the mesh starting from the tensor mesh equal to the domain boundary. Regardless of which segment is considered for insertion first, it will always end inside the initial box $\cE_0$ (the whole domain) and will not traverse it. In contrast, Figure~\ref{fig:noiterativemesh}(b) shows that a slightly finer initial tensor spline mesh allows for the construction of an iterative spline mesh resembling the one in~(a).

Since $\cM_\m{p}^\mu$ is a spline mesh, the multiplicity function $\mu = \mu_N$ satisfies 
$1 \leq \mu_N(\gamma) \leq p_k + 1$ for any minimal $k$-meshline of $\cM = \cM_N$. 
Consequently, any update $\mu_i = \mu_{i-1} + \mu_{\gamma_i}$ must preserve the upper bound.

The initial tensor spline mesh $\cM_{0, \m{p}}^{\mu_0}$ has constant meshlines by definition. If we want to maintain the constant meshline property, for $i = 1, \ldots, N$,
the $i$th segment insertion in the iterative construction of $\cM_\m{p}^\mu$ must satisfy one of the following conditions:
\begin{itemize}
\item $\gamma_i$ coincides entirely with an existing meshline in $\cM_{i-1}^{\mu_{i-1}}$;
\item $\gamma_i$ constitutes a completely new meshline in $\cM_i^{\mu_i}$, i.e., it traverses boxes of $\cE_{i-1}$ but is not attached to any collinear meshline of $\cM_{i-1}^{\mu_{i-1}}$;
\item $\gamma_i$ extends an existing meshline in $\cM_{i-1}^{\mu_{i-1}}$, traversing boxes of $\cE_{i-1}$, with associated multiplicity $\mu_{\gamma_i}$ equal to that of the meshline it extends.
\end{itemize} 
In particular, each $\gamma_i$ cannot simultaneously coincide with an existing meshline and traverse boxes of $\cE_{i-1}$.

\begin{remark}
In the literature, an iterative spline mesh with constant meshlines is sometimes called an \emph{LR mesh}. We prefer to reserve this term to identify a proper subclass; see Section~\ref{sec:bivariatelr}.
\end{remark}

It is useful for later to define knot vectors on the meshlines.
\begin{definition}\label{def:knotvecongamma}
Let $\gamma$ be a $(k, c)$-meshline in a spline mesh $\cM_\m{p}^\mu$. All vertices where $\gamma$ intersects meshlines of $\cM$ have their $k$th coordinates equal to $c$, while their $(3-k)$th coordinates are distinct. The \textbf{global knot vector on} $\gamma$, denoted by $\pmb{\tau}_{p_{3-k}}^{\mu_{3-k}}$, is the knot vector with breakpoints given by these $(3-k)$th coordinates, and  multiplicity function given by the $\mu_{3-k}$ multiplicities of the corresponding vertices, as in Definition \ref{def:mu1mu2}.
\end{definition}

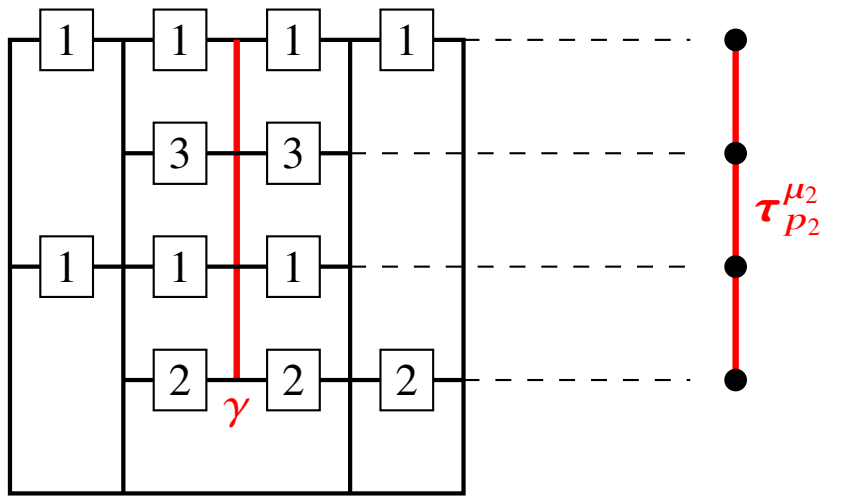
\begin{figure}[t!]
	\centering
	\begin{tikzpicture}[scale = 3]
		\tikzset{%
			block/.style    = {draw=black, fill=white, rectangle, minimum height = .05cm,
				minimum width = .05cm},
		}
		\draw[thick] (0, 0) rectangle (1, 1);
		\draw[thick] (0.25, 0) -- (0.25, 1);
		\draw[red, very thick] (0.5, 0.25)node[below]{$\gamma$} -- (0.5, 1);
		\draw[thick] (0.75, 0) -- (0.75, 1);
		\draw[thick] (0.25, 0.25) -- (1, 0.25);
		\draw[thick] (0, 0.5) -- (0.75, 0.5);
		\draw[thick] (0.25, 0.75) -- (0.75, 0.75);
		\draw[dashed] (1, 1) -- (1.5, 1);
		\draw[dashed] (0.75, 0.75) -- (1.5, 0.75);
		\draw[dashed] (0.75, 0.5) -- (1.5, 0.5);
		\draw[dashed] (1, 0.25) -- (1.5, 0.25);
		\draw[red, very thick] (1.6, 0.25) --node[right]{$\pmb{\tau}_{p_2}^{\mu_2}$} (1.6, 1);
		\fill (1.6, 0.25) circle (0.025);
		\fill (1.6, 0.5) circle (0.025);
		\fill (1.6, 0.75) circle (0.025);
		\fill (1.6, 1) circle (0.025);
		\draw (0.375, 0.25) node[block]{$2$};
		\draw (0.625, 0.25) node[block]{$2$};
		\draw (0.875, 0.25) node[block]{$2$};
		\draw (0.125, 0.5) node[block]{$1$};
		\draw (0.375, 0.5) node[block]{$1$};
		\draw (0.625, 0.5) node[block]{$1$};
		\draw (0.375, 0.75) node[block]{$3$};
		\draw (0.625, 0.75) node[block]{$3$};
		\draw (0.125, 1) node[block]{$1$};
		\draw (0.375, 1) node[block]{$1$};
		\draw (0.625, 1) node[block]{$1$};
		\draw (0.875, 1) node[block]{$1$};
	\end{tikzpicture}
	\caption{The global knot vector $\pmb{\tau}_{p_2}^{\mu_2}$ on a $1$-meshline $\gamma$ of an iterative spline mesh $\cM_{\m{p}}^\mu$ of constant splits. The multiplicities needed for the definition of $\pmb{\tau}_{p_2}^{\mu_2}$ are reported in square boxes on the minimal meshlines orthogonal to and intersecting $\gamma$.}\label{fig:knotvecongamma}
\end{figure}

Figure~\ref{fig:knotvecongamma} depicts an example of a global knot vector on a $1$-meshline $\gamma$ of an iterative spline mesh $\cM_{\m{p}}^\mu$ with constant splits. In a tensor spline mesh, all the vertical (resp., horizontal) meshlines share the same global knot vector $\pmb{\tau}_{2,p_2}^{\mu_2}$ (resp., $\pmb{\tau}_{1,p_1}^{\mu_1}$).

We now consider additional constraints on iterative spline meshes in order to ensure that the corresponding spline spaces have stable dimension, i.e., no homological term appears in their dimension formula. In view of the later construction of LR B-splines, these constraints are called \textbf{LR rules}.

\begin{itemize}
\item \textbf{LR rule 1:} The spline mesh is an iterative spline mesh with constant meshlines.
\item \textbf{LR rule 2:} The spline space $\SS[\cM_{0, \m{p}}^{\mu_0}]$ over the initial tensor spline mesh $\cM_{0, \m{p}}^{\mu_0}$ is non-trivial, i.e., $\SS[\cM_{0, \m{p}}^{\mu_0}] \neq \{0\}$.
\item \textbf{LR rule 3:} For $i = 1, \ldots, N$, and $k \in \{1,2\}$, the spline space $\SS[\pmb{\tau}_{p_{3-k}}^{\mu_{3-k}}]$ defined over the global knot vector $\pmb{\tau}_{p_{3-k}}^{\mu_{3-k}}$ associated with any $k$-meshline of $\cM_{i, \m{p}}^{\mu_i}$ is non-trivial, i.e., $\SS[\pmb{\tau}_{p_{3-k}}^{\mu_{3-k}}] \neq \{0\}$.
\end{itemize} 

\begin{example}{Example: LR Rules}%
The open spline meshes with $\m{p} = (2,2)$ and simple internal meshlines in Figure~\ref{fig:dimunstable}(a)--(b) satisfy LR rules 1--2 but not LR rule 3. Indeed, regardless of the insertion order, the first of the four shorter meshlines has a global knot vector of length $2$, and hence defines a trivial univariate spline space.
\end{example}

The next result follows from \cite[Corollary~2, Lemma~17, and Theorem~6]{splinedim}.
\begin{theorem}
Let $\cM_\m{p}^\mu$ be a spline mesh  obeying  LR rules 1--3. Then, the bivariate dimension formula \eqref{eq:splinedim2d} has vanishing homological term.
 \label{thm:pettersen}
\end{theorem}

When $\cM_\m{p}^{\mu}$ is a spline mesh obeying the LR rules, another way to compute the dimension of $\SS[\cM_\m{p}^{\mu}]$ is by describing the increase in the spline space dimension after insertion of a new segment.
In this perspective, the next theorem relates the dimension of the refined spline space to that of the current spline space and to the dimension of the univariate spline space defined on the global knot vector associated with the new meshline created by the insertion; see \cite[Theorem~5.5]{lr}.
\begin{theorem}\label{thm:splinedimincrease2d}
Let $\cM_\m{p}^\mu$ be an iterative spline mesh with constant meshlines such that $\cM = \cM' + \gamma$ and $\mu = \mu' + \mu_{\gamma}$, where $\gamma$ is a $k$-meshline not contained in $\cM'$.  Let $\pmb{\tau}_{p_{3-k}}^{\tilde{\mu}_{3-k}}$ be the global knot vector on $\gamma$, with multiplicity function $\tilde{\mu}_{3-k}$ equal to $\mu_{3-k}$, as in Definition~\ref{def:knotvecongamma}, except when $\gamma$ extends a meshline of $\cM'$. In that case, the breakpoints in $\pmb{\tau}$ corresponding to the junction vertices with $\cM'$ are assigned full multiplicity.
Then,
\begin{equation*}
    \dim \SS[\cM_\m{p}^{\mu}] = \dim \SS[(\cM')_{\m{p}}^{\mu'}] + \dim \SS[\pmb{\tau}_{p_{3-k}}^{\tilde{\mu}_{3-k}}].
\end{equation*}
\end{theorem}

\begin{example}{Example: Computation of Spline Dimension}%
Let $p \geq 0$ and consider the bidegree $\m{p} = (p,p)$. Let $(\cM')_{\m{p}}^{\mu'}$ be the open iterative spline mesh shown in Figure~\ref{fig:exsplinedim2d}(a)--(b), with simple internal meshlines.
According to the tables in Figure~\ref{fig:exsplinedim2d}(c)--(d), we have
\begin{itemize}
\item $\displaystyle \sum_{\m{v} \in \cV'} (p - \mu_1(\m{v}) + 1)(p - \mu_2(\m{v}) + 1) =
2p^2$,
\item $\displaystyle -(p+1)\sum_{\gamma \in \cM'} (p - \mu(\gamma) +1) =
-(p + 1)6p = -6p^2 -6p$,
\item $|\cE'|(p + 1)^2 = 5(p + 1)^2 = 5p^2 + 10p + 5$.
\end{itemize}
Moreover, the LR rules are easily verified for $(\cM')_{\m{p}}^{\mu'}$. Indeed,
\begin{itemize}
\item LR rule 1: The mesh $\cM'$ can be obtained by inserting two horizontal segments into the initial tensor spline mesh shown in Figure~\ref{fig:exsplinedim2d1}(a). 
\item LR rule 2: The dimension of the spline space on the initial tensor mesh in Figure~\ref{fig:exsplinedim2d1}(a) can be computed by means of \eqref{eq:dimtensorproduct} and \eqref{eq:splinedim1d}, which gives $(p + 3)(p + 1) = p^2 + 4p + 3 > 0$ for all $p \geq 0$. 
\item  LR rule 3: Both spline spaces associated with the global knot vectors on the meshlines of $\cM'$ created by the insertion of the two segments have dimension $1$. 
\end{itemize}
Therefore, from Theorem~\ref{thm:pettersen}, the homological term in \eqref{eq:splinedim2d} vanishes, and we conclude
\begin{equation*}
\dim \SS[(\cM')_{\m{p}}^{\mu'}] = 2p^2 - 6p^2 - 6p + 5p^2 + 10p + 5 = p^2 + 4p + 5.
\end{equation*}

\begin{figure}[t!]
\centering
\subfloat[]{
\begin{tikzpicture}[scale = 3]
\tikzset{%
  block/.style    = {draw=black, fill=white, rectangle, minimum height = .05cm,
    minimum width = .05cm},
    }
\draw (0, 1) --node[block]{$\gamma_1$} (0.5, 1);
\draw (0.5, 1) --node[block]{$\gamma_2$} (1, 1);
\draw (1, 1) --node[block]{$\gamma_3$} (1.5, 1);

\draw (0, 0) --node[block]{$\gamma_4$} (0.5, 0);
\draw (0.5, 0) --node[block]{$\gamma_5$} (1, 0);
\draw (1, 0) --node[block]{$\gamma_6$} (1.5, 0);

\draw (0, 0) --node[block]{$\gamma_7$} (0, 0.5);
\draw (0, 0.5) --node[block]{$\gamma_8$} (0, 1);
\draw (0, 0.5) --node[block]{$\gamma_9$} (0.5, 0.5);
\draw (0.5, 0) --node[block]{$\gamma_{10}$} (0.5, 0.5);
\draw (0.5, 0.5) --node[block]{$\gamma_{11}$} (0.5, 1);

\draw (1, 0) --node[block]{$\gamma_{12}$} (1, 0.5);
\draw (1, 0.5) --node[block]{$\gamma_{13}$} (1, 1);
\draw (1, 0.5) --node[block]{$\gamma_{14}$} (1.5, 0.5);
\draw (1.5, 0) --node[block]{$\gamma_{15}$} (1.5, 0.5);
\draw (1.5, 0.5) --node[block]{$\gamma_{16}$} (1.5, 1);

\draw[white] (0, -0.125);
\end{tikzpicture}
}\quad
\subfloat[]{
\begin{tikzpicture}[scale = 3]
\tikzset{%
  block/.style    = {draw=black, fill=white, rectangle, minimum height = .05cm,
    minimum width = .05cm},
    }
\draw (0, 1) --node[block]{$p + 1$} (0.5, 1);
\draw (0.5, 1) --node[block]{$p + 1$} (1, 1);
\draw (1, 1) --node[block]{$p + 1$} (1.5, 1);

\draw (0, 0) --node[block]{$p + 1$} (0.5, 0);
\draw (0.5, 0) --node[block]{$p + 1$} (1, 0);
\draw (1, 0) --node[block]{$p + 1$} (1.5, 0);

\draw (0, 0) --node[block]{$p + 1$} (0, 0.5);
\draw (0, 0.5) --node[block]{$p + 1$} (0, 1);
\draw (0, 0.5) --node[block]{$1$} (0.5, 0.5);
\draw (0.5, 0) --node[block]{$1$} (0.5, 0.5);
\draw (0.5, 0.5) --node[block]{$1$} (0.5, 1);

\draw (1, 0) --node[block]{$1$} (1, 0.5);
\draw (1, 0.5) --node[block]{$1$} (1, 1);
\draw (1, 0.5) --node[block]{$1$} (1.5, 0.5);
\draw (1.5, 0) --node[block]{$p + 1$} (1.5, 0.5);
\draw (1.5, 0.5) --node[block]{$p + 1$} (1.5, 1);

\fill (0, 1)node[above]{$\m{v}_1$} circle (0.015);
\fill (0.5, 1)node[above]{$\m{v}_2$} circle (0.015);
\fill (1, 1)node[above]{$\m{v}_3$} circle (0.015);
\fill (1.5, 1)node[above]{$\m{v}_4$} circle (0.015);

\fill (0, 0)node[below]{$\m{v}_5$} circle (0.015);
\fill (0.5, 0)node[below]{$\m{v}_6$} circle (0.015);
\fill (1, 0)node[below]{$\m{v}_7$} circle (0.015);
\fill (1.5, 0)node[below]{$\m{v}_8$} circle (0.015);

\fill (0, 0.5)node[left]{$\m{v}_9$} circle (0.015);
\fill (0.5, 0.5)node[right]{$\m{v}_{10}$} circle (0.015);
\fill (1, 0.5)node[left]{$\m{v}_{11}$} circle (0.015);
\fill (1.5, 0.5)node[right]{$\m{v}_{12}$} circle (0.015);
\end{tikzpicture}
}\\
\subfloat[]{
\scalebox{0.85}{
\renewcommand{\arraystretch}{1.5}
\begin{tabular}{|c|c|c|c|c|c|c|c|c|c|c|c|c|}
\hline
  \,~~~$i$~~~\, &  1 & 2 & 3 & 4 & 5 & 6 & 7 & 8 & 9 & 10 & 11 & 12 \\ \hline
$\mu_1(\m{v}_i)$ & $p + 1$ & 1 & 1 & $p + 1$ & $p + 1$ & 1 & 1 & $p + 1$ & $p + 1$ & \,~~~1~~~\, & \,~~~1~~~\, &  $p + 1$\\
$\mu_2(\m{v}_i)$ & $p + 1$ & $p + 1$ & $p + 1$ & $p + 1$ & $p + 1$ & $p + 1$ & $p + 1$ & $p + 1$ & 1 & 1 & 1 &  1 \\
$V_i$ & 0 & 0 & 0 & 0 & 0 & 0 & 0 & 0 & 0 & $p^2$ & $p^2$ &  0\\ \hline
\end{tabular}}
}\\
\subfloat[]{
\scalebox{0.825}{
\renewcommand{\arraystretch}{1.5}
\begin{tabular}{|c|c|c|c|c|c|c|c|c|c|c|c|c|c|c|c|c|}
\hline
  \,~~~$j$~~~\, &  1 & 2 & 3 & 4 & 5 & 6 & 7 & 8 & 9 & 10 & 11 & 12 & 13 & 14 & 15& 16\\ \hline
$\mu(\gamma_j)$ & $p + 1$ & $p + 1$ & $p + 1$ & $p + 1$ & $p + 1$ & $p + 1$ & $p + 1$ & $p + 1$ & \,~~~1~~~\, & \,~~~1~~~\, & \,~~~1~~~\, & \,~~~1~~~\, & \,~~~1~~~\, & \,~~~1~~~\, & $p + 1$ & $p + 1$ \\
$G_j$ & 0 & 0 & 0 & 0 & 0 & 0 & 0 & 0 & $p$ & $p$ & $p$ & $p$ & $p$ & $p$ & 0 & 0 \\ \hline
\end{tabular}}
}
\caption{An open iterative spline mesh of bidegree $\m{p} = (p, p)$, where the minimal meshlines are indexed in (a) and the vertices are indexed in (b). Also the minimal meshline multiplicities are reported in (b). The horizontal and vertical multiplicities of the vertices are reported in (c), together with $V_i \defeq (p - \mu_1(\m{v}_i) + 1)(p - \mu_2(\m{v}_i) + 1)$,
while the meshline multiplicities are reported in (d), together with $G_j \defeq p - \mu(\gamma_j) + 1$.	}\label{fig:exsplinedim2d}
\end{figure}

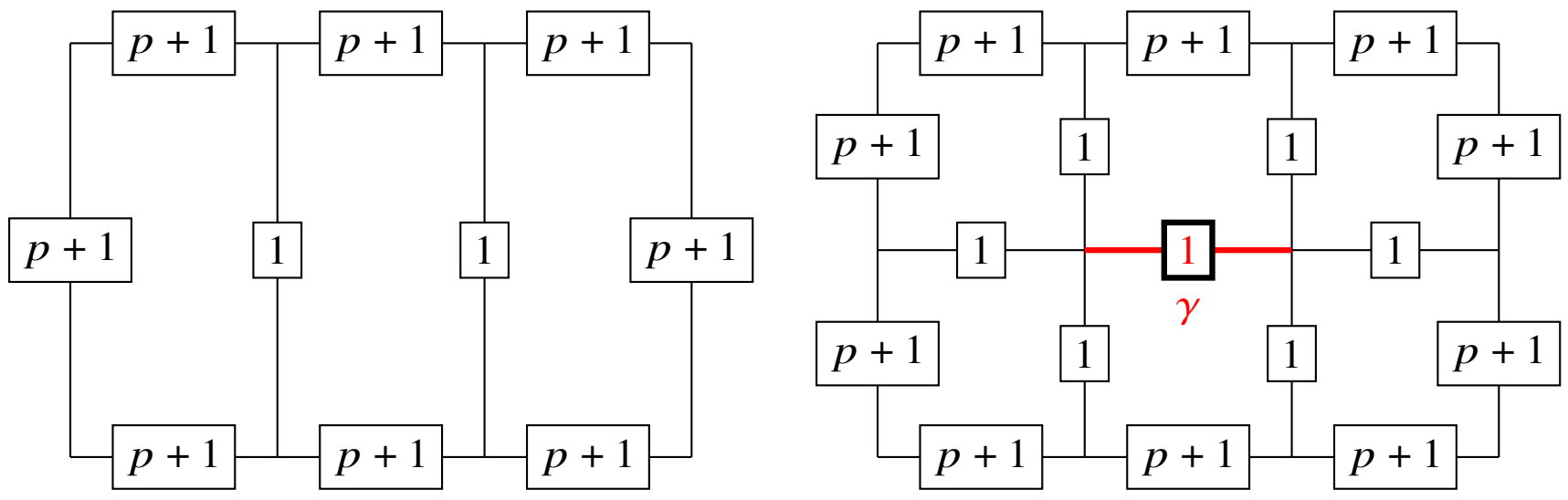
\begin{figure}[t!]
\centering
\subfloat[]{
\begin{tikzpicture}[scale = 3]
\tikzset{%
  block/.style    = {draw=black, fill=white, rectangle, minimum height = .05cm,
    minimum width = .05cm},
    }
\draw (0, 1) --node[block]{$p + 1$} (0.5, 1);
\draw (0.5, 1) --node[block]{$p + 1$} (1, 1);
\draw (1, 1) --node[block]{$p + 1$} (1.5, 1);

\draw (0, 0) --node[block]{$p + 1$} (0.5, 0);
\draw (0.5, 0) --node[block]{$p + 1$} (1, 0);
\draw (1, 0) --node[block]{$p + 1$} (1.5, 0);

\draw (0, 0) --node[block]{$p + 1$} (0, 1);
\draw (0.5, 0) --node[block]{$1$} (0.5, 1);

\draw (1, 0) --node[block]{$1$} (1, 1);
\draw (1.5, 0) --node[block]{$p + 1$} (1.5, 1);
\end{tikzpicture}
}\quad
\subfloat[]{
\begin{tikzpicture}[scale = 3]
\tikzset{%
  block/.style    = {draw=black, fill=white, rectangle, minimum height = .05cm,
    minimum width = .05cm},
    }
\draw (0, 1) --node[block]{$p + 1$} (0.5, 1);
\draw (0.5, 1) --node[block]{$p + 1$} (1, 1);
\draw (1, 1) --node[block]{$p + 1$} (1.5, 1);

\draw (0, 0) --node[block]{$p + 1$} (0.5, 0);
\draw (0.5, 0) --node[block]{$p + 1$} (1, 0);
\draw (1, 0) --node[block]{$p + 1$} (1.5, 0);

\draw (0, 0) --node[block]{$p + 1$} (0, 0.5);
\draw (0, 0.5) --node[block]{$p + 1$} (0, 1);
\draw (0, 0.5) --node[block]{$1$} (0.5, 0.5);
\draw (0.5, 0) --node[block]{$1$} (0.5, 0.5);
\draw (0.5, 0.5) --node[block]{$1$} (0.5, 1);

\draw (1, 0) --node[block]{$1$} (1, 0.5);
\draw (1, 0.5) --node[block]{$1$} (1, 1);
\draw (1, 0.5) --node[block]{$1$} (1.5, 0.5);
\draw (1.5, 0) --node[block]{$p + 1$} (1.5, 0.5);
\draw (1.5, 0.5) --node[block]{$p + 1$} (1.5, 1);

\draw[red, very thick] (0.5, 0.5) --node[block]{$1$} (1, 0.5);
\node at (0.75, 0.35) {\textcolor{red}{$\gamma$}};
\end{tikzpicture}
}
\caption{The iterative spline mesh in Figure~\ref{fig:exsplinedim2d} is generated from the initial tensor spline mesh in (a). The tensor spline mesh in (b) is a refinement of it obtained by the insertion of $\gamma$, which connects two local meshlines.}\label{fig:exsplinedim2d1}
\end{figure}

Assume now that a new segment $\gamma$ is inserted as in Figure~\ref{fig:exsplinedim2d1}(b), yielding a spline mesh $\cM_{\m{p}}^\mu$ with $\cM = \cM' + \gamma$ and $\mu = \mu' + \mu_\gamma$, where $\mu_\gamma = 1$ in order to preserve the constant meshline property. We compute the dimension of $\SS[\cM_{\m{p}}^\mu]$ using the dimension increasing formula of Theorem~\ref{thm:splinedimincrease2d}.
The segment $\gamma$ connects two meshlines of $\cM'$. The associated global knot vector $\pmb{\tau}_p^{\tilde{\mu}_1}$ has two breakpoints $ (\tau_1, \tau_2)$ corresponding to the vertices where $\gamma$ is attached to $\cM'$. By definition of $\tilde{\mu}_1$, we have $\tilde{\mu}_1(\tau_k) = p + 1$ for $k = 1,2$.
By \eqref{eq:splinedim1d}, the spline space $\SS[\pmb{\tau}_p^{\tilde{\mu}_1}]$ has dimension $p + 1$. Therefore, by Theorem~\ref{thm:splinedimincrease2d},
\begin{equation*}
\begin{split}
\dim \SS[\cM_{\m{p}}^\mu] &= \dim \SS[(\cM')_{\m{p}}^{\mu'}] + \dim \SS[\pmb{\tau}_p^{\tilde{\mu}_1}] \\
&= (p^2 + 4p + 5) + (p + 1) = p^2 + 5p + 6.
\end{split}
\end{equation*}
Since $\cM$ is a tensor mesh, $\dim \SS[\cM_{\m{p}}^\mu]$ can also be computed as the product of the dimensions of the univariate spline spaces associated with the corresponding global knot vectors; see \eqref{eq:dimtensorproduct}. Using \eqref{eq:splinedim1d}, these dimensions are $p + 3$ and $p + 2$ along the two directions; their product agrees with the computed value of $\dim \SS[\cM_{\m{p}}^\mu]$.
\end{example}

\subsection{Multivariate Spline Space}
In this section, we briefly describe how the notions of box partition, mesh, and spline space extend to $\RR^d$. The corresponding dimension formula involves contributions from all $\ell$-dimensional boxes, with $\ell \in \{0, \ldots, d\}$, and generalizes the univariate and bivariate cases discussed above. We also provide sufficient conditions ensuring the vanishing of the homological term in the dimension formula.

We begin with some additional definitions. An interval $I \defeq [a, b]$, with $a, b \in \RR$, is called \textbf{trivial} if $a = b$, and \textbf{non-trivial} if $a<b$.

\begin{definition}\label{def:boxdd}
Given an integer $d \geq 1$ and $\ell \in \{0, \ldots, d\}$, an \textbf{$\ell$-box} in $\RR^d$ is a Cartesian product of finite intervals in $\RR$,
\begin{equation}\label{eq:box_def}
\beta \defeq I_1 \times \cdots \times I_d \subset \RR^d,
\end{equation}
in which exactly $\ell$ intervals are non-trivial and $d-\ell$ intervals are trivial.
The integer $\ell$ is called the \textbf{dimension} of $\beta$, and we write $\dim \beta = \ell$.
\end{definition} 

In particular, $0$-boxes correspond to vertices, $1$-boxes to edges, etc. 
In analogy with the bivariate case, we state the following definition.
\begin{definition}
Let $\Omega \subset \RR^d$ be a $d$-box. A finite collection $\cE$ of $d$-boxes in $\RR^d$ is said to be a \textbf{box partition} of $\Omega$ if
\begin{itemize}
\item $\mathring{\beta}_1 \cap \mathring{\beta}_2 = \varnothing$ for any $\beta_1, \beta_2 \in \cE$ and $\beta_1 \neq \beta_2$;
\item $ \bigcup_{\beta \in \cE} \beta = \Omega$.
\end{itemize}
The elements of $\cE$ are called \textbf{boxes}, or \textbf{cells}, or \textbf{elements}.
\end{definition}

\begin{definition}
Given a box partition $\cE$ of a $d$-box $\Omega \subset \RR^d$, for any $\beta \in \cE$, with Cartesian representation
$\beta = I_1 \times \cdots \times I_d$, $I_k = [a_{\beta,k}, b_{\beta,k}]$, we set $I_k^- \defeq \{a_{\beta,k}\}$, $I_k^+ \defeq \{b_{\beta,k}\}$, for each $k=1,\dots,d$,
and we define
\begin{equation*}
\begin{split}
\cF_\ell(\beta) \defeq \Biggl\{ \beta' &= \bigtimes_{k=1}^d J_k \,:\, J_k \in \left\{I_k, I_k^-, I_k^+\right\},\ k = 1, \ldots, d, \text{ and}\\
&\quad \text{ exactly $\ell$ of these intervals are non-trivial} \Biggr\}, \quad \ell=0,\dots,d.
\end{split}
\end{equation*}
Moreover, we define
\begin{equation*}
\begin{aligned}
\cF_0(\cE) &\defeq \bigcup_{\beta \in \cE} \cF_0(\beta),
\\
\cF_\ell(\cE) &\defeq \Biggl\{ \beta' \in \bigcup_{\beta \in \cE} \cF_\ell(\beta) \,:\,
\mathring{\beta}' \cap \gamma = \varnothing \ \ \forall \gamma \in \bigcup_{\ell'=0}^{\ell-1} \cF_{\ell'}(\cE) \Biggr\},
\quad \ell=1,\dots,d.
\end{aligned}
\end{equation*}
\end{definition}

For $\ell = 1, \ldots, d$, the set $\cF_\ell(\cE)$ consists of $\ell$-boxes whose interiors do not intersect any lower-dimensional entities from the sets $\cF_{\ell'}(\cE)$, $\ell'<\ell$.

\begin{definition}\label{def:boxmesh3D}
Given a box partition $\cE$ of a $d$-box $\Omega \subset \RR^d$, the \textbf{box mesh}, or simply the \textbf{mesh}, associated with the box partition is the collection $\cM \defeq \cF_{d-1}(\cE)$. The elements of $\cM$ are called \textbf{minimal meshrectangles}.
 Moreover, we use the following terminology:
 \begin{itemize}
\item maximal collections of minimal meshrectangles in $\cM$ whose union forms a $(d-1)$-box are called \textbf{meshrectangles}, or \textbf{meshblocks}, of $\cM$;
\item a \textbf{multiplicity function} $\mu \colon \cM \to \NN$ is a function that assigns a positive integer to each minimal meshrectangle;
\item a mesh $\cM$ equipped with a multiplicity function is called a \textbf{$\pmb{\mu}$-extended box mesh}, or simply \textbf{$\pmb{\mu}$-extended mesh}, and is denoted by $\cM^\mu$;
\item a $\mu$-extended mesh $\cM^\mu$ has \textbf{constant} (or \textbf{uniform}) \textbf{meshrectangles}, or \textbf{constant} (or \textbf{uniform}) \textbf{meshblocks}, if each intersecting pair of minimal $k$-meshrectangles has the same multiplicity. In this case, the multiplicity of a meshrectangle in $\cM$ is defined as the common multiplicity of its constituent minimal meshrectangles.
\end{itemize}
\end{definition}

We sometimes refer to a meshrectangle as a $k$\textbf{-meshrectangle}, where $k \in \{1, \ldots, d\}$ is the index of the trivial interval $I_k = \{c\}$, for some $c \in \RR$, in the Cartesian representation \eqref{eq:box_def} of such a meshrectangle. We may also write $(k,c)$\textbf{-meshrectangle} to specify the parameter value in the $k$th direction.

\begin{definition}
The set of minimal $k$-meshrectangles in a box mesh $\cM$ is denoted by $\cM_{[k]}$, $k=1,\dots, d$.
\end{definition}

\begin{figure}[t!]
	\centering
	\raisebox{-0.075\textwidth}{\subfloat[]{
			\includegraphics[width = 0.15\textwidth, page = 1]{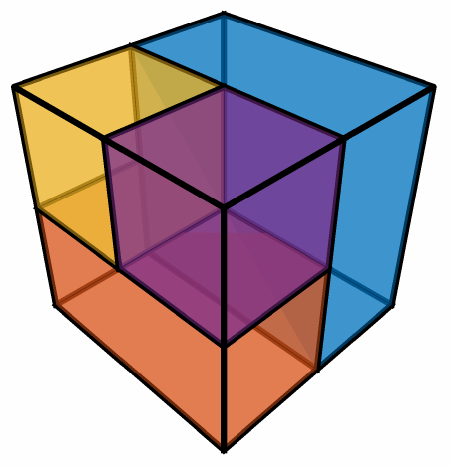}
	}}\hspace{-0.3cm}
	\subfloat[]{
		\raisebox{-0.15\textwidth}{
			\begin{tikzpicture}
				\draw[thick,-stealth] (0,0) -- (0.6,0.325) node[above right=-3]{$x_2$};
				\draw[thick,-stealth] (0,0) -- (0.5,-0.45) node[right=-1]{$x_1$};
				\draw[thick,-stealth] (0,0) -- (-0.05,0.75) node[above = -2]{$x_3$};
		\end{tikzpicture}}\hspace{-1.1cm}
		$
		\begin{array}{c}
			\includegraphics[width = 0.15\textwidth, page = 2]{boxes}~
			\includegraphics[width = 0.15\textwidth, page = 3]{boxes}\vspace{-0.15cm}\\
			\includegraphics[width = 0.15\textwidth, page = 4]{boxes}
		\end{array}
		$
	}
	\subfloat[]{
		$\begin{array}{c}
			\includegraphics[width = 0.15\textwidth, page = 5]{boxes}~
			\includegraphics[width = 0.15\textwidth, page = 6]{boxes}\vspace{-0.15cm}\\
			\includegraphics[width = 0.15\textwidth, page = 7]{boxes}
		\end{array}
		$
	}
	\raisebox{-0.075\textwidth}{\subfloat[]{
			\includegraphics[width = 0.15\textwidth, page = 8]{boxes}
	}}
	\caption{All sets $\cF_\ell(\cE)$ for a box partition $\cE$ of a $3$-box in $\RR^3$. The box partition $\cE = \cF_3(\cE)$ is shown in (a) and the box mesh $\cM = \cF_2(\cE)$ in (b), which contains the minimal meshrectangles, i.e., the facets in the three-dimensional setting: the minimal $1$-meshrectangles (top left), the minimal $2$-meshrectangles (top right), and the minimal $3$-meshrectangles (bottom). The set $\cF_1(\cE)$ of minimal edges is shown in (c) and the set $\cF_0(\cE)$ of vertices in (d).}\label{fig:exboxpartition3d}
\end{figure}

When $d = 1$, the boxes of a box partition correspond to knot intervals, and the (minimal) meshrectangles correspond to breakpoints. When $d = 2$, the boxes are rectangles and the minimal meshrectangles correspond to minimal meshlines.
Figure~\ref{fig:exboxpartition3d} depicts an example of a three-dimensional box partition $\cE$ together with the associated sets $\cF_\ell(\cE)$ for all $\ell = 0, \ldots, 3$.

\begin{remark}
The definition of a $\mu$-extended mesh with constant meshrectangles given in Definition~\ref{def:boxmesh3D} is consistent with the one in Definition~\ref{def:boxmesh2D} for the bivariate case. Indeed, in two dimensions, any two intersecting minimal $k$-meshlines necessarily belong to the same $k$-meshline. In the multivariate case, however, two intersecting minimal $k$-meshrectangles do not need to belong to the same $k$-meshrectangle.
Therefore, our definition of a box mesh with constant meshrectangles imposes a stronger multiplicity constraint than the uniform-multiplicity requirement for meshrectangles in \cite{splinedim}, where more general multiplicity distributions are allowed.
\end{remark}

\begin{definition}
A $\mu$-extended mesh $\cM^\mu$ is a \textbf{tensor mesh} if it has constant meshrectangles and, for each $k = 1, \ldots, d$, there exists a sequence $(a_{k,1}, \ldots, a_{k,m_k})$ in $\RR$, with $a_{k,1} < \cdots < a_{k,m_k}$, such that the associated box partition $\cE$ has the form
\begin{equation*}
\cE = \Biggl\{ \beta = \bigtimes_{k=1}^d [a_{k,i_k}, a_{k,i_k+1}] \,:\, 1 \leq i_k \leq m_k - 1,\ k = 1, \ldots, d \Biggr\}.
\end{equation*}
\end{definition}

With this geometric structure, the definition of a spline mesh in the bivariate setting generalizes directly to the multivariate case as follows.
\begin{definition}
A \textbf{spline mesh} in $\RR^d$ is a pair $\cM_{\m{p}}^\mu \defeq (\cM^\mu, \m{p})$, where $\m{p} \defeq (p_1, \ldots, p_d)$ is a sequence of positive integers and $\cM^\mu$ is a $\mu$-extended mesh associated with a box partition $\cE$, with multiplicities satisfying $1 \leq \mu(\gamma) \leq p_k + 1$ for every minimal $k$-meshrectangle $\gamma \in \cM$.
Moreover, we use the following terminology:
\begin{itemize}
\item $\cM_{\m{p}}^\mu$ is \textbf{open} if every boundary minimal $k$-meshrectangle has 
multiplicity $p_k+1$; 
\item $\cM_{\m{p}}^\mu$ is a \textbf{tensor spline mesh} if $\cM^\mu$ is a tensor mesh.
\end{itemize}
\end{definition}

Given a multidegree $\m{p} \defeq (p_1, \ldots, p_d)$, let $\PP_\m{p}$ be the multivariate polynomial space of degree at most $p_k$, $k=1,\dots,d$, componentwise, i.e., the space spanned by the monomials
\begin{equation*}
x_1^{i_1} \cdots x_d^{i_d}
\quad \text{such that} \quad
0 \leq i_k \leq p_k, \quad k = 1, \ldots, d.
\end{equation*}
We are now ready to define spline spaces in the multivariate setting.

\begin{definition}
Given a spline mesh $\cM_{\m{p}}^\mu$ corresponding to a box partition $\cE$ of $\Omega = [a_1,b_1] \times \cdots \times [a_d, b_d]$, for any element $\beta \in \cE$, with Cartesian representation
$\beta = I_1 \times \cdots \times I_d$, $I_k = [a_{\beta,k}, b_{\beta,k}]$, we define
\begin{equation*}
\tilde{\beta} \defeq \tilde{I}_1 \times \cdots \times \tilde{I}_d, \quad\text{where}\quad
\tilde{I}_k \defeq
\begin{cases}
[a_{\beta,k}, b_{\beta,k}), & \text{if } b_{\beta,k} < b_k, \\
[a_{\beta,k}, b_{\beta,k}], & \text{if } b_{\beta,k} = b_k.
\end{cases}
\end{equation*}
The \textbf{spline space on} $\cM_{\m{p}}^\mu$, denoted by $\mathbb{S}[\cM_{\m{p}}^\mu]$, is the set of functions $f \colon \RR^d \to \RR$ such that
\begin{itemize}
\item $f$ is zero outside $\Omega$;
\item for each element $\beta \in \cE$, the restriction $f|_{\tilde{\beta}}$ belongs to $\PP_\m{p}$;
\item for each minimal $k$-meshrectangle $\gamma \in \cM$, $f$ is $C^{p_k - \mu(\gamma)}$-continuous across $\gamma$, i.e., with respect to the $k$th direction.
\end{itemize}
\end{definition}

The following theorem provides the dimension of multivariate spline spaces on arbitrary spline meshes; see \cite[Theorem~2]{splinedim}.

\begin{theorem}
Let $\cM_{\m{p}}^\mu$ be a spline mesh corresponding to a box partition $\cE$ in $\RR^d$. Then,
\begin{equation}\label{eq:splinedimdd}
\begin{split}
    &\dim\,\mathbb{S}[\cM_\m{p}^\mu] \\
    &\quad=|\cE|\prod_{k = 1}^d (p_k+1)+ \sum_{\ell = 0}^{d-1} (-1)^{d-\ell} \left(\sum_{\beta \in \cF_\ell(\cE)} \prod_{k=1}^d(p_k - \mu_k(\beta) + 1)\right)\\
    &\quad\quad +   \text{homological term},
 \end{split}
\end{equation}
where 
for $\beta \in \cF_\ell(\cE)$ and $k = 1, \ldots, d$, the quantity $\mu_k(\beta)$ is defined as
\begin{equation*}
\mu_k(\beta) \defeq 
\begin{cases}
\mu(\gamma), & \text{if } \beta \subseteq \gamma \text{ for some } \gamma \in \cM_{[k]}, \\
0, & \text{otherwise}.
\end{cases}
\end{equation*}
\end{theorem}

\begin{example}{Example: Dimension of Bivariate Spline Spaces}%
Let $d=2$. Since $\cF_0(\cE)$ is the set of vertices of $\cM$ and $\cF_1(\cE)$ is the set of minimal meshlines of $\cM$, we have
\begin{equation*}
	\begin{split}
		&\sum_{\ell = 0}^{d-1} (-1)^{d-\ell} \left(\sum_{\beta \in \cF_\ell(\cE)} \prod_{k=1}^d(p_k - \mu_k(\beta) + 1)\right) \\
		&\quad =
		  \sum_{\m{v}\in\mathcal{V}} (p_1-\mu_1(\m{v}) +1)(p_2-\mu_2(\m{v})+1)\\
		&\quad\quad -\sum_{\gamma \in \cM_{[1]}} (p_1-\mu_1(\gamma)+1) (p_2-\mu_2(\gamma)+1) \\
		&\quad\quad -\sum_{\gamma \in \cM_{[2]}}(p_1-\mu_1(\gamma)+1)(p_2-\mu_2(\gamma)+1).
	\end{split}
\end{equation*}
Thus, the dimension formula \eqref{eq:splinedimdd} reduces to \eqref{eq:splinedim2d} because $\mu_{3-k}(\gamma)=0$ for $\gamma\in \cM_{[k]}$.
\end{example}

In the multivariate setting, the homological term in \eqref{eq:splinedimdd} can be ensured to vanish under suitable sufficient conditions, which generalize those imposed for $d=2$ (see Section~\ref{sec:space-2d}).
We require the mesh to be constructed iteratively from an initial tensor mesh by consecutive insertions of $(d-1)$-boxes. The refinement procedure after the insertion of a $(d-1)$-box is essentially the same as in the bivariate case, with meshrectangles replacing meshlines. Accordingly, the definition of \textbf{iterative spline mesh} (see Definition \ref{def:iterativesplinemesh}) extends naturally to the multivariate setting.
Therefore, LR rules 1--2 extend straightforwardly.

To find a proper formulation for LR rule 3, we consider any meshrectangle $\gamma$ of the spline mesh $\cM_{\m{p}}^\mu$. On  $\gamma$, we have an induced box partition $\cE_\gamma$ consisting of $(d-1)$-boxes, and the minimal meshrectangles of the induced mesh $\cM_\gamma$ are $(d-2)$-boxes. Their multiplicities are inherited from those of the meshrectangles in $\cM$ intersecting $\gamma$, analogous to Definition~\ref{def:mu1mu2} for $d=2$.
This construction reduces to the global knot vector on $\gamma$ in the bivariate case (see Definition~\ref{def:knotvecongamma}). Then, in view of Remark~\ref{rmk:homology1d}, a multivariate extension of LR rule 3 is to require that the dimension formula of the spline space over the induced spline mesh must have vanishing homological term.

Moreover, for $d \geq 3$, an additional condition is needed to guarantee the vanishing of the homological term in \eqref{eq:splinedimdd}. This condition is based on the following concept.

\begin{definition}
Let $\gamma$ be a $(d-1)$-box with associated multiplicity $\mu_\gamma$. Given an iterative spline mesh $(\cM')_{\m{p}}^{\mu'}$ in $\RR^d$ with constant meshrectangles, we say that $\gamma$ \textbf{slots into} $(\cM')_{\m{p}}^{\mu'}$ if
\begin{itemize}
\item the refined spline mesh $\cM_{\m{p}}^\mu$, with $\cM = \cM' + \gamma$ and $\mu = \mu' + \mu_\gamma$, has constant meshrectangles;
\item for every meshrectangle $\gamma' \in \cM'$ such that $\gamma \cap \gamma' = \beta \neq \varnothing$, we have $\dim \beta = d - 2$.
\end{itemize}
\end{definition}

\begin{figure}[t!]
\centering
\subfloat[]{
\includegraphics[width = 0.21\textwidth, page = 1]{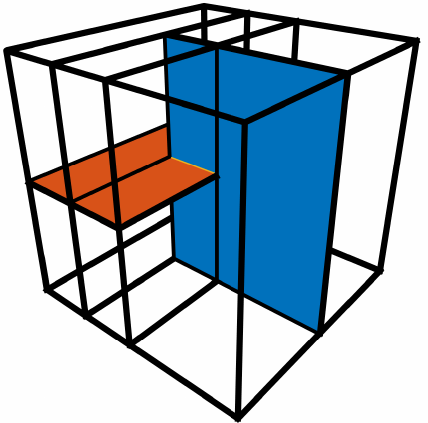}
}\quad
\subfloat[]{
\includegraphics[width = 0.21\textwidth, page = 2]{slotinto}
}\quad
\subfloat[]{
\includegraphics[width = 0.21\textwidth, page = 3]{slotinto}
}\quad
\subfloat[]{
\includegraphics[width = 0.21\textwidth, page = 4]{slotinto}
}
\caption{Examples of the insertion of a new $2$-box $\gamma$ (blue) into a mesh $\cM'$ in $\RR^3$. In (a) $\gamma$ slots into $\cM'$ because, e.g., the intersection $\beta$ with the orthogonal meshrectangle (red) is a segment (yellow), so $\dim \beta = 1$. In (b) $\gamma$ intersects the orthogonal meshrectangle (red) only in a vertex (yellow), so it does not slot into $\cM'$. In (c) it does because, e.g., the intersection with the parallel meshrectangle (red) is a segment (yellow). In (d) $\gamma$ does not slot into $\cM'$ because the intersection with the parallel meshrectangle (red) is just a vertex (yellow).}\label{fig:exslotinto}
\end{figure}

Figure~\ref{fig:exslotinto} illustrates examples of $(d-1)$-boxes that do and do not slot into a given spline mesh, for $d = 3$. Roughly speaking, the new $(d-1)$-box must intersect the existing mesh $\cM'$ in a sufficiently large way along all coordinate directions $k \in \{1, \ldots, d\}$ in order to slot into it. In the bivariate setting, if $\gamma$ preserves the constant split property, then it automatically slots into the mesh. Indeed, if $\gamma \cap \gamma' = \beta \neq \varnothing$, then $\beta$ must be a vertex, so that $\dim \beta = 0 = d - 2$.

We then arrive at the following LR rules for $d \geq 3$.
\begin{itemize}
\item \textbf{LR rule 1:} The spline mesh is an iterative spline mesh with constant meshrectangles.
\item \textbf{LR rule 2:} The spline space $\SS[\cM_{0, \m{p}}^{\mu_0}]$ over the initial tensor spline mesh $\cM_{0, \m{p}}^{\mu_0}$ is non-trivial, i.e., $\SS[\cM_{0, \m{p}}^{\mu_0}] \neq \{0\}$.
\item \textbf{LR rule 3:} The spline space in $\RR^{d-1}$ defined over the spline mesh induced on any meshrectangle of $\cM_{\m{p}}^\mu$ has vanishing homological term in its dimension formula.
\item \textbf{LR rule 4:} All $(d-1)$-boxes inserted during the construction of the iterative spline mesh must slot into the meshes into which they are inserted.
\end{itemize}

As the dimensionality increases, enforcing these conditions becomes more challenging.
The next theorem follows from \cite[Corollary~2, Lemma~17, and Theorem~6]{splinedim}.

\begin{theorem}
Let $\cM_\m{p}^\mu$ be a spline mesh in $\RR^d$ obeying LR rules 1--4. Then, the multivariate dimension formula \eqref{eq:splinedimdd} has vanishing homological term.
\end{theorem}

\section{Locally Refined B-Splines: Definition and Basic Properties}\label{sec:lr-definition}
In this section, we introduce LR meshes and LR B-splines.
We review their definition and basic properties, as well as the core algorithms for generating them.
We begin with the univariate case, and then address the bivariate case in greater detail. The extension from the bivariate to the general multivariate setting is straightforward.

One of the strengths of LR B-splines is that they reduce to standard B-splines in the univariate setting. More generally, they are nothing but standard tensor-product B-splines whenever the underlying mesh is a tensor spline mesh. This is in contrast to other adaptive spline technologies, like (truncated) hierarchical B-splines. The treatment of the univariate case thus reduces to a brief review of classical B-spline theory. More comprehensive introductions can be found in \cite{foundation} and \cite{standard}, as well as in the books \cite{deboor} and \cite{schumaker}.

\subsection{Univariate LR B-Splines}
In the univariate setting, LR B-splines coincide with standard univariate B-splines. Therefore, we simply refer to them as B-splines in this case.

Let $\pmb{\tau}_p^{\mu}$ be a global knot vector as in Definition~\ref{def:knotvector}, with $\sum_{i=1}^m \mu(\tau_i)> p+1$. To define the corresponding B-splines, it is convenient to consider
a vector $\m{x} \defeq (x_1, \ldots, x_{n+p+1})$, with non-decreasing components,
where
\begin{equation*}
n+p+1 \defeq \sum_{i=1}^m \mu(\tau_i),
\end{equation*}
and
\begin{equation*}
\underbrace{x_1 = \cdots = x_{\mu(\tau_1)}}_{=\, \tau_1}
< \underbrace{x_{\mu(\tau_1)+1} = \cdots = x_{\mu(\tau_1)+\mu(\tau_2)}}_{=\, \tau_2}
< \cdots .
\end{equation*}
The vector $\m{x}$ is equivalent to $\pmb{\tau}_p^{\mu}$ and will also be referred to as (\textbf{global}) \textbf{knot vector}.
In this notation, the elements of $\m{x}$ are called \textbf{knots}, rather than breakpoints. The multiplicity of a knot $x_r$, denoted by $\mu(x_r)$, is the number of times the value $x_r$ appears in the sequence $\m{x}$. Note that $1 \leq \mu(x_r) \leq p+1$. When the first and the last knots have multiplicity $p+1$, i.e., when $\pmb{\tau}_p^{\mu}$ is open, $\m{x}$ is called \textbf{open} as well.

\begin{definition}\label{def:bspline}
Let $\m{x} \defeq (x_1, \ldots, x_{n+p+1})$ be a knot vector. For $j = 1, \ldots, n$, the \textbf{$j$th B-spline of degree $p$ on the knot vector $\m{x}$}, denoted by $B_{j, p, \m{x}} \colon \RR \to \RR$, is defined recursively by
\begin{equation*}
B_{j, p, \m{x}}(x) \defeq
\frac{x - x_j}{x_{j+p} - x_j} \, B_{j, p-1, \m{x}}(x)
+
\frac{x_{j+p+1} - x}{x_{j + p + 1} - x_{j+1}} \, B_{j+1, p - 1, \m{x}}(x),
\end{equation*}
for $p\geq 1$, starting from
\begin{equation*}
\begin{aligned}
B_{j, 0, \m{x}}(x) &\defeq
\begin{cases}
	1, & x\in [x_j, x_{j + 1}),  \\
	0, & \text{otherwise},
\end{cases}
\quad j=1,\dots, n+p-1,
\\
B_{n+p, 0, \m{x}}(x) &\defeq \begin{cases}
		1, & x\in [x_{n+p}, x_{n+p + 1}],  \\
	0, & \text{otherwise}.
\end{cases}
\end{aligned}
\end{equation*}
Here we assume that any fraction with zero denominator is defined to be zero.
The global B-spline set is defined by $\cB_{p, \m{x}} \defeq \{B_{j, p, \m{x}}\}_{j=1}^n$.
For $ p+1\leq r\leq n$ such that $x_r < x_{r+1}$, the corresponding local B-spline set is defined by
$ \cB_{r,p,\m{x}} \defeq \{B_{j,p,\m{x}}\}_{j=r-p}^r$.
\end{definition}

\begin{figure}[t!]
\subfloat[]{
\includegraphics[width = \textwidth, page = 1]{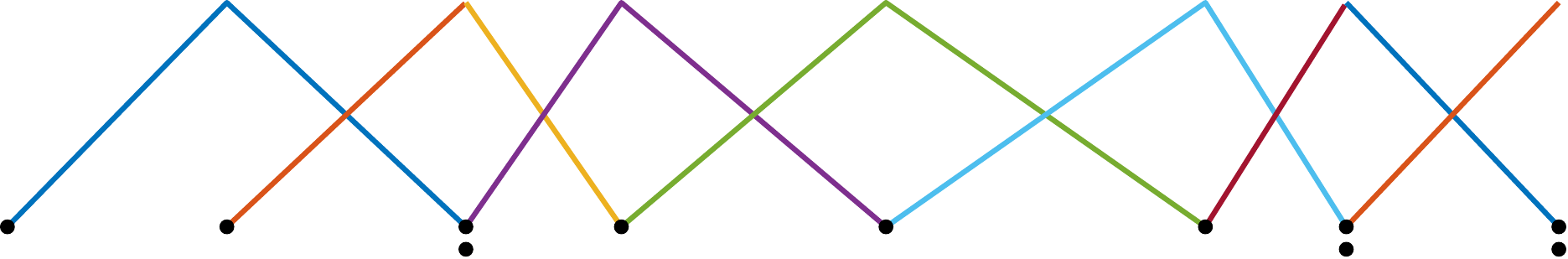}
}\\
\subfloat[]{
\includegraphics[width = \textwidth, page = 2]{unibsplines}
}
\caption{B-splines of degree $p = 1$ (a) and $p = 2$ (b) defined on a general global knot vector. The knots are depicted with black dots. In particular, those of multiplicity equal to $2$ have double dots. The distinct B-splines are drawn with different colors and only within their supports, omitting where they vanish.}
\label{fig:bsplineex}
\end{figure}

Figure~\ref{fig:bsplineex} shows two examples of B-spline sets.
The properties listed below form the foundations for working with B-splines.
They are divided into two groups. The first group (Proposition~\ref{prop:bsplinebasicproperties}) follows directly from the definition of B-splines on a knot vector $\m{x}$, while the second group (Proposition~\ref{prop:bsplineadvancedproperties}) requires a more delicate analysis. An exhaustive list of properties is beyond the scope of this work
and we refer the interested reader to \cite{foundation} and references therein.

\begin{proposition}\label{prop:bsplinebasicproperties}
Let $B_{j, p, \m{x}}$ be as in Definition~\ref{def:bspline}.
\begin{itemize}
\item \textbf{Local knot dependence:} $B_{j,p,\m{x}}$ depends only on the $p+2$ knots in the vector
$\m{x}_j \defeq (x_j, \ldots, x_{j+p+1})$, which is called the \textbf{local knot vector} of $B_{j,p,\m{x}}$.

\item \textbf{Local support:} The support of $B_{j, p, \m{x}}$ is the knot span of $\m{x}_j$, i.e.,
\begin{equation*}
\supp B_{j,p, \m{x}} = [x_j, x_{j+p+1}].
\end{equation*}
Moreover, $B_{j,p, \m{x}}(x) = 0$ if $x \in [x_r, x_{r+1})$ and $j\notin \{r-p,\dots,r\}$.

\item \textbf{Positivity:} $B_{j,p,\m{x}}$ is positive in the interior of its support, i.e.,
\begin{equation*}
B_{j,p, \m{x}}(x) > 0, \quad x \in (x_j, x_{j+p+1}).
\end{equation*}
\item \textbf{Piecewise polynomial structure:} $B_{j,p,\m{x}}$ is a piecewise polynomial of degree at most $p$, i.e.,
\begin{equation*}
B_{j,p,\m{x}}|_{[x_r, x_{r+1})} \in \PP_p, \quad r=j,\dots, j+p.
\end{equation*}

\item \textbf{Affine invariance:} For any $\alpha > 0$ and $\beta \in \RR$, it holds
\begin{equation*}
B_{j,p,\alpha \m{x} + \beta}(\alpha x + \beta) = B_{j,p,\m{x}}(x).
\end{equation*}
\end{itemize}
\end{proposition}

From the local support property, we see that the set $\cB_{r,p,\m{x}}$ contains the $p+1$ B-splines not vanishing on the interval $[x_r, x_{r+1})$.
A proof of the following properties can be found in, e.g., \cite{foundation}.

\begin{proposition}\label{prop:bsplineadvancedproperties}
Let $B_{j, p, \m{x}}$ be as in Definition~\ref{def:bspline}, and suppose $ p+1\leq r\leq n$ such that $x_r < x_{r+1}$.
\begin{itemize}
\item \textbf{Smoothness}:  For $k\in\{j,\dots, j+p+1\}$, we have that $B_{j, p, \m{x}}$ is exactly $C^{p-\mu_j(x_k)}$-continuous at $x_k$, where $\mu_j(x_k)$ is the number of times $x_k$ appears in the local knot vector $\m{x}_j$.
\item \textbf{Derivative formula:} Let $D_+$ denote the right derivative. For $p \geq 1$, it holds
\begin{equation}\label{eq:Bsplinederivative}
D_+ B_{j,p,\m{x}}(x) = p \left(
\frac{B_{j,p-1,\m{x}}(x)}{x_{j+p}-x_j}
-
\frac{B_{j+1,p-1,\m{x}}(x)}{x_{j+p+1}-x_{j+1}}
\right),
\end{equation}
with the usual convention for zero denominators.
%
%
\item \textbf{Global and local partition of unity:} It holds
\begin{equation*}
\sum_{j=1}^n B_{j,p,\m{x}}(x) = 1, \quad x\in[x_{p+1}, x_{n+1}],
\end{equation*}
and
\begin{equation*}
\sum_{j=r-p}^r B_{j,p,\m{x}}(x) = 1, \quad x \in [x_r, x_{r+1}).
\end{equation*}
\item \textbf{Global and local monomial reproduction:} For $1 \leq q \leq p$, it holds
\begin{equation*}
x^q = \sum_{j=1}^n x_{j,p}^{\overbrace{\scalebox{0.75}{$*\cdots*$}}^{q\text{ times}}} \, B_{j,p,\m{x}}(x), \quad x \in [x_{p+1}, x_{n+1}],
\end{equation*}
and
\begin{equation*}
x^q = \sum_{j=r-p}^r x_{j,p}^{\overbrace{\scalebox{0.75}{$*\cdots*$}}^{q\text{ times}}} \, B_{j,p,\m{x}}(x), \quad x \in [x_r, x_{r+1}),
\end{equation*}
where
\begin{equation*}
x_{j,p}^{*\cdots*} \defeq
 \sum (x_{j_1}x_{j_2}\cdots x_{j_q})\, \bigg/\, {p \choose q},
\end{equation*}
and the sum is taken over all combinations of $q$-tuples $(j_1,\ldots,j_q)$ in the index set $\{j+1,\ldots,j+p\}$.
For $q=1$, the coefficients $x_{j,p}^*$ are the average of the internal knots in the local knot vector $\m{x}_j$ of $B_{j, p, \m{x}}$, i.e.,
\begin{equation*}
x_{j,p}^* = \frac{x_{j+1} + \cdots + x_{j+p}}{p},
\end{equation*}
and are called the \textbf{Greville abscissae} or \textbf{Greville points}.
\item \textbf{Global and local linear independence:} If $x_{p+1}<x_{p+2}$ and $x_n<x_{n+1}$, then the functions in $\cB_{p,\m{x}}$ are linearly independent on $[x_{p+1}, x_{n+1}]$.
Moreover, the functions in $\cB_{r,p,\m{x}}$ are linearly independent on $[x_r, x_{r+1})$.
\end{itemize}
\end{proposition}

\begin{remark}\label{rmk:convexhull1d}
The positivity and (local) partition of unity properties imply that the B-splines in $\cB_{p, \m{x}}$ (and $\cB_{r, p, \m{x}}$) form a \textbf{positive partition of unity} (also called \textbf{convex partition of unity}). This has an important geometric consequence.
Consider a spline curve in $\RR^d$ ($d \geq 2$), represented in terms of $\cB_{p, \m{x}}$, i.e.,
\begin{equation*}
\mathbf{C}(x) \defeq \sum_{j=1}^n \m{P}_j \, B_{j, p, \m{x}}(x),\quad \m{P}_j \in \RR^d, \quad x\in[x_{p+1}, x_{n+1}].
\end{equation*}
It can be locally expressed as
\begin{equation*}
\mathbf{C}(x) = \sum_{j=r-p}^r \m{P}_j \, B_{j, p, \m{x}}(x),\quad x\in[x_r, x_{r+1}).
\end{equation*}
Therefore, the curve $\mathbf{C}$ lies within the \textbf{convex hull} of the $p+1$ consecutive points $\m{P}_{r-p}, \ldots, \m{P}_r$ for $x\in[x_r, x_{r+1})$.
This convex hull property allows us to predict the qualitative behavior of the curve from the displacement of the points $\{\m{P}_j\}_{j=1}^n$, which is particularly useful in geometric design.
\end{remark}

Without loss of generality, from now on we assume that the knot vector $\m{x}$ is open.
From Propositions~\ref{prop:bsplinebasicproperties} and \ref{prop:bsplineadvancedproperties} it follows that $\spn \cB_{p, \m{x}}$ has dimension $n$. Moreover, recalling Definition~\ref{def:splinespace1d}, we deduce that
$\spn \cB_{p, \m{x}} = \SS[\pmb{\tau}_p^\mu]$, which means that $ \cB_{p, \m{x}}$ is a basis for $ \SS[\pmb{\tau}_p^\mu]$. This property is known as the \emph{Curry--Schoenberg theorem} \cite{curryschoenberg}. In particular, when we focus on the restriction to $[x_r,x_{r+1})$, we have that $\spn \cB_{r,p,\m{x}} = \PP_p$.

\begin{definition}
A knot vector $\hat{\m{x}}$ is said to be a \textbf{refinement} of a knot vector $\m{x}$ if every knot in $\m{x}$ appears in $\hat{\m{x}}$ with at least the same multiplicity.
\end{definition}

For a fixed degree $p$, let $\hat{\m{x}}$ be a refinement of $\m{x}$. Then, any function $f \in \spn \cB_{p, \m{x}}$ satisfies at least the continuity constraints required to belong to $\spn \cB_{p, \hat{\m{x}}}$, so
\begin{equation*}
\spn \cB_{p, \m{x}} \subseteq \spn \cB_{p, \hat{\m{x}}}.
\end{equation*}
A B-spline given on the coarse knot vector $\m{x}$ can thus be described in terms of B-splines on the refined knot vector $\hat{\m{x}}$. This relation is called \textbf{knot insertion} (or \textbf{knot refinement}); see, e.g., \cite[Section~1.3.5]{foundation}.

\begin{theorem}\label{thm:knotinsertion}
Let $\hat{\m{x}} \defeq (\hat{x}_i)_{i=1}^{n+p+2}$ be a knot vector obtained from the knot vector $\m{x} \defeq (x_i)_{i=1}^{n+p+1}$ by inserting a single knot $\hat{x}$, with $x_r \leq \hat{x} < x_{r+1}$, $\hat x\in[x_{p+1},x_{n+1}]$.
For $j = r-p, \ldots, r$, we have
\begin{equation*}
B_{j,p,\m{x}} = \alpha_{1,j} \, B_{j,p,\hat{\m{x}}} + \alpha_{2,j} \, B_{j+1,p,\hat{\m{x}}},
\end{equation*}
where
\begin{equation*}
\alpha_{1,j} \defeq
\begin{cases}
1, & \text{if } x_{j+p} \leq \hat{x} < x_{j+p+1}, \\[4pt]
\displaystyle \frac{\hat{x} - x_j}{x_{j+p} - x_j}, & \text{if } x_j \leq \hat{x} < x_{j+p},
\end{cases}
\end{equation*}
and
\begin{equation*}
\alpha_{2,j} \defeq
\begin{cases}
1, & \text{if } x_j \leq \hat{x} \leq x_{j+1}, \\[4pt]
\displaystyle \frac{x_{j+p+1} - \hat{x}}{x_{j+p+1} - x_{j+1}}, & \text{if } x_{j+1} < \hat{x} < x_{j+p+1}.
\end{cases}
\end{equation*}
Moreover, for $j < r-p$, we have $B_{j,p,\m{x}} = B_{j,p,\hat{\m{x}}}$, while for $j > r$, we have $B_{j,p,\m{x}} = B_{j+1,p,\hat{\m{x}}}$.
\end{theorem}

Suppose $f\in\spn \cB_{p, \m{x}}$. Theorem~\ref{thm:knotinsertion} implies
\begin{equation*}
f(x) = \sum_{j=1}^n c_j B_{j,p,\m{x}}(x)
= \sum_{i=1}^{n+1} \hat{c}_i B_{i,p,\hat{\m{x}}}(x), \quad x \in [x_{p+1}, x_{n+1}],
\end{equation*}
where
\begin{equation*}
	\hat{c}_i \defeq
	\begin{cases}
		c_i, & \text{if } i \leq r-p, \\[4pt]
		\displaystyle
		\frac{\hat{x} - x_i}{x_{i+p} - x_i} \, c_i
		+
		\frac{x_{i+p} - \hat{x}}{x_{i+p} - x_i} \, c_{i-1},
		& \text{if } r-p < i \leq r, \\[4pt]
		c_{i-1}, & \text{if } i > r.
	\end{cases}
\end{equation*}
More efficient refinement algorithms are available when multiple knots are inserted simultaneously. In particular, the \emph{Oslo algorithm} \cite{oslo1,oslo2} provides an elegant way to perform multiple knot insertions in a single step.

\begin{remark}\label{rmk:pou1d}
In line with the notation of Theorem~\ref{thm:knotinsertion},
suppose $\alpha_{1,j}:=1$ and $\alpha_{2,j}:=0$ for $j=1,\ldots,r-p-1$. Moreover, suppose $\alpha_{1,j}:=0$ and $\alpha_{2,j}:=1$ for $j=r+1,\ldots,n$.
Then, each $B_{j,p,\m{x}}\in \cB_{p, \m{x}}$ is expressed in terms of $B_{j,p,\hat{\m{x}}}$ and $B_{j+1,p,\hat{\m{x}}}$, with coefficients $\alpha_{1,j}$ and $\alpha_{2,j}$, respectively.
With the convention $\alpha_{2,0} \defeq \alpha_{1,n+1} \defeq 0$, we also define
\begin{equation*}
w_j \defeq \alpha_{1,j} + \alpha_{2,j-1}, \quad j = 1, \ldots, n+1,
\end{equation*}
and consider the scaled set
\begin{equation*}
w\cB_{p,\hat{\m{x}}} \defeq \bigl\{w_j B_{j,p,\hat{\m{x}}}\bigr\}_{j=1}^{n+1}.
\end{equation*}
Starting from the partition of unity property for $\cB_{p,\m{x}}$, we deduce
\begin{equation*}
\begin{aligned}
1
&= \sum_{j=1}^n B_{j,p,\m{x}}(x) = \sum_{j=1}^n \bigl( \alpha_{1,j} B_{j,p,\hat{\m{x}}}(x) + \alpha_{2,j} B_{j+1,p,\hat{\m{x}}}(x) \bigr) \\
&= \sum_{j=1}^{n+1} \bigl(\alpha_{1,j} + \alpha_{2,j-1}\bigr) B_{j,p,\hat{\m{x}}}(x) = \sum_{j=1}^{n+1} w_j B_{j,p,\hat{\m{x}}}(x).
\end{aligned}
\end{equation*}
Thus, the scaled set $w\cB_{p,\hat{\m{x}}}$ forms a partition of unity.
Actually, the linear independence and partition of unity of the set
$\cB_{p,\hat{\m{x}}}$ imply that $w_j = 1$ for all $j$.
This property will no longer hold in the context of general spline meshes.
\end{remark}

\begin{figure}[t!]
	\centering
	\includegraphics[width = \textwidth, page = 3]{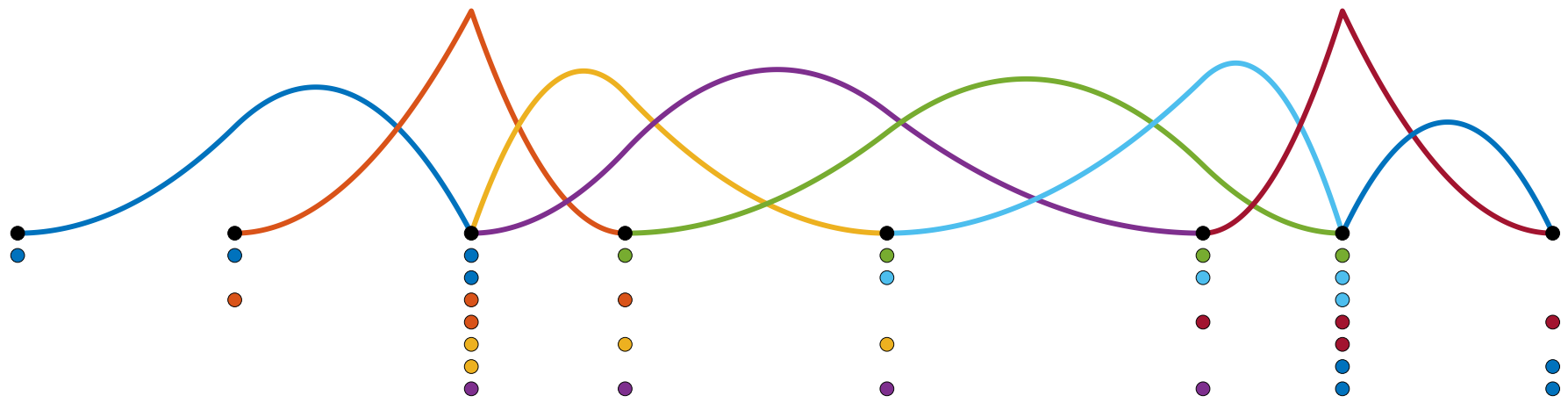}
	\caption{Quadratic B-splines defined on a knot vector as in Figure~\ref{fig:bsplineex}(b). The breakpoints are depicted with black dots. The B-splines are drawn with different colors and only within their supports, omitting where they vanish. Under each representation we report the related local knot vector with dots colored as the corresponding B-spline. Knots of multiplicity equal to $2$ in the local knot vector have double dots.
	}\label{fig:exlocalknotvector}
\end{figure}

The theory of LR B-splines relies on two fundamental ingredients: the local knot dependence of B-splines and the knot insertion relation.
The local knot dependence allows us to decouple the definition of a single B-spline from the rest of the basis functions. In particular, evaluating a function $B_{j, p, \m{x}}(x)$ does not require knowledge of the full knot vector $\m{x} \defeq (x_1, \ldots, x_{n+p+1})$, but only of its local knot vector $\m{x}_j \defeq (x_j, \ldots, x_{j+p+1})$; see
Definition \ref{def:bspline}.
This locality is exploited in efficient evaluation procedures, such as the \emph{de Boor algorithm} \cite{debooralgo}.
Similarly, when refining a B-spline, only its local knot vector is required; see Theorem~\ref{thm:knotinsertion}.

For these reasons, in the following sections we will replace the global knot vector notation with a local knot vector notation.
Figure~\ref{fig:exlocalknotvector} illustrates the local knot vectors associated with the quadratic B-splines shown in Figure~\ref{fig:bsplineex}(b).
In contrast to the representation in Figure~\ref{fig:bsplineex}, multiplicities here are defined with respect to the local knot vector, so they are the values of the functions $\mu_j$ in Proposition~\ref{prop:bsplineadvancedproperties}.

\subsection{Bivariate LR B-Splines}\label{sec:bivariatelr}
As already mentioned, the definition of LR B-splines relies on local knot vectors. Let $p \geq 0$ and let $\m{x}_B$ be a local knot vector of length $p+2$. We denote by $B[\m{x}_B]$ the unique B-spline associated with $\m{x}_B$. The degree $p$ is omitted, as it is implicitly determined by the length of $\m{x}_B$.

We consider spline meshes $\cM_{\m{p}}^\mu$ of fixed bidegree $\m{p} \defeq (p_1,p_2)$ as in Definition~\ref{def:spline-mesh}. For simplicity of notation, we sometimes drop the subscript $\m{p}$ and write $\cM^\mu$.
Let $\cM^\mu$ be a spline mesh on a rectangle $\Omega$ satisfying LR rules 1--2. Thus, $\cM^\mu$ has constant meshlines. Moreover, let $\cM_B^{\mu_B}$ be a submesh with constant meshlines, forming a rectangular grid on a rectangle $\Omega_B \subseteq \Omega$.
The multiplicity function $\mu_B \colon \cM_B \to \NN$ satisfies
\begin{equation*}
\mu_B(\gamma) \leq \mu(\gamma), \quad \forall \gamma \in \cM_B.
\end{equation*}
The submesh $\cM_B^{\mu_B}$ is a tensor spline mesh on $\Omega_B$.
Assume that this mesh consists of $p_1+2$ vertical and $p_2+2$ horizontal meshlines, where each meshline is counted with the multiplicity given by $\mu_B$.
Then, it can be represented by a pair of local knot vectors $(\m{x}_B, \m{y}_B)$, each of length $p_k+2$, $k=1,2$. This allows us to define the univariate B-splines $B[\m{x}_B]$ and $B[\m{y}_B]$, and the tensor-product B-spline
\begin{equation*}
B[\m{x}_B, \m{y}_B] \defeq B[\m{x}_B] \, B[\m{y}_B].
\end{equation*}
When no confusion arises, we simply write $B$ instead of $B[\m{x}_B, \m{y}_B]$. We say that $\cM_B^{\mu_B}$ is the \textbf{local tensor mesh} associated with $B$. Note that the support of $B$ is $\Omega_B$.

In order to define LR B-splines, we need the following support concepts.
\begin{definition}
Let $\cM^\mu$ be a spline mesh satisfying LR rules 1--2. We say that a B-spline $B[\m{x}_B, \m{y}_B]$ has \textbf{support on} $\cM^\mu$ if the associated  spline mesh $\cM_B^{\mu_B}$ satisfies
\begin{equation*}
\cM_B \subseteq \cM,
\quad
\mu_B(\gamma) \leq \mu(\gamma) \quad \forall \gamma \in \cM_B.
\end{equation*}
Furthermore, we say that $B$ has \textbf{minimal support} on $\cM^\mu$ if there does not exist a meshline $\gamma \in \cM$ traversing the interior of $\supp B$ such that
\begin{equation*}
\mu_B(\gamma) < \mu(\gamma),
\end{equation*}
where $\mu_B(\gamma) \defeq 0$ for all $\gamma \in \cM \setminus \cM_B$.
The set of \textbf{minimal support B-splines} (or \textbf{MS B-splines}) on $\cM^\mu$ is denoted by $\cB \defeq \cB[\cM^\mu]$.
\end{definition}

\begin{figure}[t!]
	\centering
	\subfloat[]{
		\begin{tikzpicture}[scale = 0.8]
			\draw[step = 1, thick] (1, 1) grid (4, 4);
			\draw[thick] (1.5, 2) -- (1.5, 4);
			\draw[thick] (2.5, 1) -- (2.5, 3.5);
			\draw[thick] (1.5, 2.5) -- (4, 2.5);
			\draw[thick] (1.5, 3.5) -- (4, 3.5);
		\end{tikzpicture}
	}\quad
	\subfloat[]{
		\begin{tikzpicture}[scale = 0.8]
			\filldraw[draw = red, fill = red!50!white, ultra thick] (1.5, 1.5) rectangle (3.5, 3.5);
			\draw[red, ultra thick] (3, 1.5) -- (3, 3.5);
			\draw[red, ultra thick] (2, 1.5) -- (2, 3.5);
			\draw[red, ultra thick] (1.5, 3) -- (3.5, 3);
			\draw[red, ultra thick] (1.5, 2.5) -- (3.5, 2.5);
			\draw[step = 1, thick] (1, 1) grid (4, 4);
			\draw[thick] (1.5, 2) -- (1.5, 4);
			\draw[thick] (2.5, 1) -- (2.5, 3.5);
			\draw[thick] (1.5, 2.5) -- (4, 2.5);
			\draw[thick] (1.5, 3.5) -- (4, 3.5);
		\end{tikzpicture}
	}\quad
	\subfloat[]{
		\begin{tikzpicture}[scale = 0.8]
			\filldraw[draw = red, fill = red!50!white, ultra thick] (1.5, 2) rectangle (4, 4);
			\draw[red, ultra thick] (3, 2) -- (3, 4);
			\draw[red, ultra thick] (2, 2) -- (2, 4);
			\draw[red, ultra thick] (1.5, 3.5) -- (4, 3.5);
			\draw[red, ultra thick] (1.5, 2.5) -- (4, 2.5);
			\draw[blue!70, ultra thick, dashed] (1.5, 3) -- (4, 3);
			\draw[step = 1, thick] (1, 1) grid (4, 4);
			\draw[thick] (1.5, 2) -- (1.5, 4);
			\draw[thick] (2.5, 1) -- (2.5, 3.5);
			\draw[thick] (1.5, 2.5) -- (4, 2.5);
			\draw[thick] (1.5, 3.5) -- (4, 3.5);
		\end{tikzpicture}
	}\quad
	\subfloat[]{
		\begin{tikzpicture}[scale = 0.8]
			\filldraw[draw = red, fill = red!50!white, ultra thick] (1, 1) rectangle (4, 4);
			\draw[red, step = 1, ultra thick] (1, 1) grid (4, 4);
			\draw[step = 1, thick] (1, 1) grid (4, 4);
			\draw[thick] (1.5, 2) -- (1.5, 4);
			\draw[thick] (2.5, 1) -- (2.5, 3.5);
			\draw[thick] (1.5, 2.5) -- (4, 2.5);
			\draw[thick] (1.5, 3.5) -- (4, 3.5);
		\end{tikzpicture}
	}
	\caption{Some biquadratic B-splines and their relation to the iterative spline mesh $\cM_{\m{p}}^\mu$, shown in (a), with constant meshlines of multiplicity $\mu = 1$ and bidegree $\m{p} = \pmb{2}$: a biquadratic B-spline with no support on $\cM_{\m{p}}^\mu$ in (b), a biquadratic B-spline with support on $\cM_{\m{p}}^\mu$ but no minimal support due to the meshline dashed in blue in (c), and a biquadratic B-spline with minimal support on $\cM_{\m{p}}^\mu$ in (d). Their supports and associated local tensor meshes are highlighted in red.}\label{fig:exminimalsupport}
\end{figure}
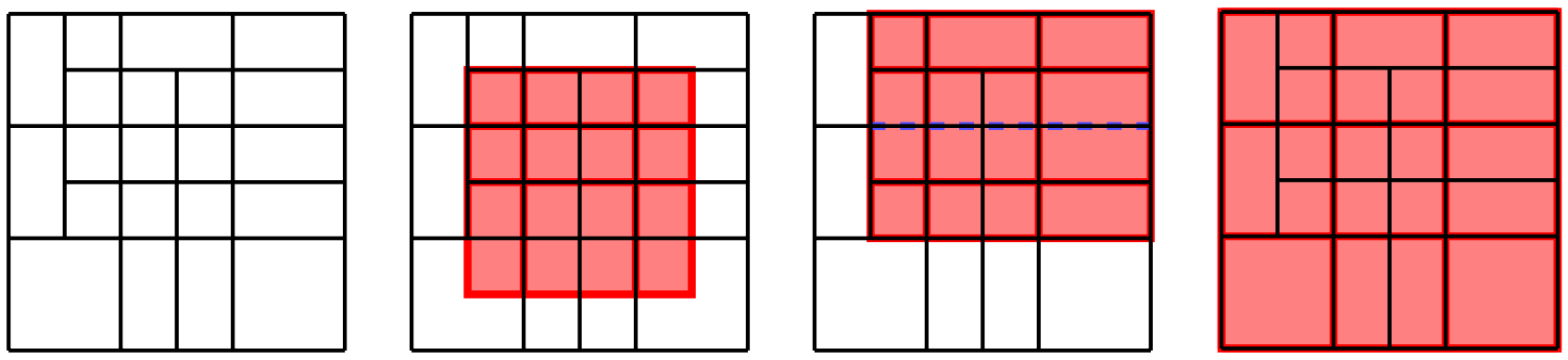

Figure~\ref{fig:exminimalsupport} illustrates B-splines with and without (minimal) support on a given mesh.
If $\cM^\mu$ is a tensor spline mesh, then the set of MS B-splines on $\cM^\mu$ coincides with the standard set of tensor-product B-splines on $\cM^\mu$.

LR meshes are particular instances of iterative spline meshes with constant meshlines, where the segment insertions are driven by the B-splines defined on the mesh. LR meshes and LR B-splines are defined in synergy \cite{lr}.

\begin{definition}\label{def:lr}
An \textbf{LR mesh} is either a tensor spline mesh or an iterative spline mesh satisfying LR rules 1--2, generated by inserting segments in such a way that each insertion breaks the minimal support property of at least one LR B-spline defined on a previous mesh.
\textbf{LR B-splines} on the LR mesh $\cM^\mu$ are defined recursively as follows.
\begin{itemize}
    \item If the underlying mesh $\cM^\mu$ is a tensor spline mesh, then the corresponding LR B-splines coincide with the tensor-product B-splines.
    \item Otherwise, LR B-splines are obtained by repeated application of the knot insertion procedure to previously defined LR B-splines on a coarser mesh, in such a way that all resulting B-splines satisfy the minimal support property.
\end{itemize}
The set of LR B-splines on $\cM^\mu$ is denoted by $\cL \defeq \cL[\cM^\mu]$.
\end{definition}

In practice, LR meshes and LR B-splines are constructed in the following manner.
Starting from an initial tensor spline mesh and the corresponding tensor-product B-splines, the mesh is refined by inserting segments that break the minimal support property of at least one B-spline supported on the mesh. Note that the ``inserted'' segments may already belong to the mesh to be refined, in which case the insertion corresponds to increasing its multiplicity. Let $B$ be such a B-spline. The insertion of a segment breaking its minimal support property results in adding one knot in either $\m{x}_B$ or $\m{y}_B$, which both consist of $p_k + 2$ consecutive knots.
By the knot insertion procedure (see Theorem~\ref{thm:knotinsertion}), $B$ can be expressed as a linear combination of two B-splines, corresponding to the first and the last $p_k + 2$ consecutive knots in the extended knot vector of length $p_k + 3$ that includes the newly inserted knot.
The set of  B-splines under consideration is then enriched by substituting $B$ by the two B-splines obtained by knot insertion.

\begin{figure}[t!]
\centering
\subfloat[]{
\begin{tikzpicture}[scale = 0.7]
\draw[thick, step = 1] (1, 1) grid (6, 4);
\draw[thick] (3, 2.5) -- (6, 2.5);
\draw[thick] (1, 3.5) -- (4, 3.5);
\draw[red, thick] (3.5, 2)node[below]{$\gamma$} -- (3.5, 3.5);
\fill (3.5, 2) circle (0.075);
\fill (3.5, 2.5) circle (0.075);
\fill (3.5, 3) circle (0.075);
\fill (3.5, 3.5) circle (0.075);
\end{tikzpicture}
}\quad
\subfloat[]{
\begin{tikzpicture}[scale = 0.7]
\fill[opacity = 0.5, orange] (1, 2) rectangle (4, 4);
\draw[ultra thick, draw = orange] (1, 2) rectangle (4, 4);
\draw[ultra thick, orange] (2, 2) -- (2, 4);
\draw[ultra thick, orange] (3, 2) -- (3, 4);
\draw[ultra thick, orange] (1, 3) -- (4, 3);
\draw[ultra thick, orange] (1, 3.5) -- (4, 3.5);
\fill[opacity = 0.5, green!50!black] (3, 2) rectangle (6, 4);
\draw[ultra thick, draw = green!50!black] (3, 2) rectangle (6, 4);
\draw[ultra thick, green!50!black] (4, 2) -- (4, 4);
\draw[ultra thick, green!50!black] (5, 2) -- (5, 4);
\draw[ultra thick, green!50!black] (3, 2.5) -- (6, 2.5);
\draw[ultra thick, green!50!black] (3, 3) -- (6, 3);
\draw[thick, step = 1] (1, 1) grid (6, 4);
\draw[thick] (3, 2.5) -- (6, 2.5);
\draw[thick] (1, 3.5) -- (4, 3.5);
\draw[red, thick] (3.5, 2)node[below]{$\gamma$} -- (3.5, 4);
\end{tikzpicture}
}\quad
\subfloat[]{
\begin{tikzpicture}[scale = 0.7]
\fill[opacity = 0.5, orange] (1, 1) rectangle (4, 3.5);
\draw[ultra thick, draw = orange] (1, 1) rectangle (4, 3.5);
\draw[ultra thick, orange] (2, 1) -- (2, 3.5);
\draw[ultra thick, orange] (3, 1) -- (3, 3.5);
\draw[ultra thick, orange] (1, 3) -- (4, 3);
\draw[ultra thick, orange] (1, 2) -- (4, 2);
\fill[opacity = 0.5, green!50!black] (3, 1) rectangle (6, 3);
\draw[ultra thick, draw = green!50!black] (3, 1) rectangle (6, 3);
\draw[ultra thick, green!50!black] (3, 2) -- (6, 2);
\draw[ultra thick, green!50!black] (3, 2.5) -- (6, 2.5);
\draw[ultra thick, green!50!black] (4, 1) -- (4, 3);
\draw[ultra thick, green!50!black] (5, 1) -- (5, 3);
\draw[thick, step = 1] (1, 1) grid (6, 4);
\draw[thick] (3, 2.5) -- (6, 2.5);
\draw[thick] (1, 3.5) -- (4, 3.5);
\draw[thick, red] (3.5, 1) -- (3.5, 3.5)node[above]{$\gamma$};
\end{tikzpicture}
}
\caption{Examples of iterative spline meshes with meshlines of multiplicity $1$ and bidegree $\m{p}=\m{2}$. The spline mesh in (a) is iterative but it is not an LR mesh because the insertion of the red $1$-meshline $\gamma$ does not refine any of the LR B-splines defined on the previous mesh. Nevertheless, LR rule 3 is satisfied. The meshes in (b)--(c) are LR meshes. The insertion of $\gamma$ (red) refines the LR B-splines whose supports and associated local tensor meshes are highlighted with colors.}\label{fig:exdifferencelrmeshanditerativemesh}
\end{figure}
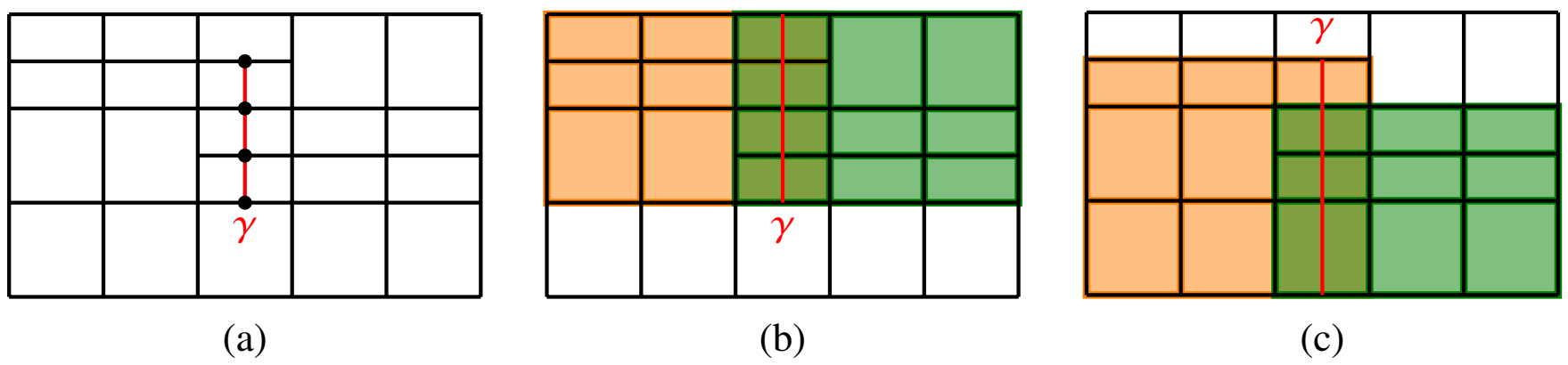

Figure~\ref{fig:exdifferencelrmeshanditerativemesh}(a)--(c) highlights the difference between a general iterative spline mesh and an LR mesh.
Note that, since the local tensor mesh $\cM_B^{\mu_B}$ of any LR B-spline $B$ consists of $p_k + 2$ $k$-meshlines for $k \in \{1,2\}$, LR rule 3 is automatically satisfied by LR meshes.

\begin{algorithm}[t!]
\caption{LR B-spline update -- $\cL = \texttt{lr}(\cL', \cM^\mu)$}\label{alg:LRgeneration}
\setlength{\columnsep}{-5cm}
\nonl\begin{multicols}{2}
\textbf{Input:}\vspace{0.1cm}
\begin{itemize}
\item[~~$\cL'\colon$] ~~~~current LR B-spline set\vspace{0.15cm}
\item[$\cM^\mu\colon$] ~~~~refinement of the LR mesh of $\cL'$
\end{itemize}\vfill\null
\columnbreak
\hspace{4cm}
\textbf{Output:}\vspace{0.1cm}

\hspace{4.25cm}\begin{minipage}{0.3\textwidth}
\begin{itemize}
\item[$\cL\colon$] ~LR B-spline set on $\cM^\mu$
\end{itemize}\end{minipage}
\vfill\null
\end{multicols}\vspace{-1cm}
\rule{0.965\textwidth}{0.025cm}

\linespread{1.35}\selectfont
Initialize $\cL\leftarrow\cL'$\;
\nonl\While{there exists $B \in \cL$ that does not have minimal support on $\cM^\mu$}{
    Select a meshline $\gamma$ of $\cM$ that breaks the minimal support property of $B$\;

    Compute the knot insertion related to $\gamma$: $B = \alpha_{1,B} B_1 + \alpha_{2,B} B_2$\;
    \nonl\For{$k = 1, 2$}{
 \nonl\If{$B_k \notin \cL$}{
    Add $B_k$ to $\cL$: $\cL \leftarrow \cL \cup \{B_k\}$\;
    }}
    Remove $B$ from $\cL$: $\cL \leftarrow \cL \setminus \{B\}$\;
}
\end{algorithm}

\begin{figure}[t!]
\centering
\subfloat[]{
\begin{tikzpicture}[scale=3.6]
\filldraw[fill=red!50!white,draw=red, line width=4pt] (.2,.4) -- (.8,.4) -- (.8,1) -- (.2,1) -- cycle;
\draw[red,line width=4pt] (.4,.4) -- (.4,1);
\draw[red,line width=4pt] (.6,.4) -- (.6,1);
\draw[red,line width=4pt] (.2,.6) -- (.8,.6);
\draw[red,line width=4pt] (.2,.8) -- (.8,.8);
\draw[step=.2] (0,.2) grid (1,1.2);
\draw (.5,.2) -- (.5,.8);
\draw[dashed] (.2,.7) -- (.8,.7);
\end{tikzpicture}
}
\subfloat[]{
\begin{tikzpicture}[scale=3.6]
\filldraw[fill=green!50!black,opacity=.5,draw=green!50!black, line width=4pt] (.2,.6) -- (.8,.6) -- (.8,1) -- (.2,1) -- cycle;
\draw[green!50!black,line width=4pt,opacity=.5] (.4,.6) -- (.4,1);
\draw[green!50!black,line width=4pt,opacity=.5] (.6,.6) -- (.6,1);
\draw[green!50!black,line width=4pt,opacity=.5] (.2,.7) -- (.8,.7);
\draw[green!50!black,line width=4pt,opacity=.5] (.2,.8) -- (.8,.8);

\filldraw[fill=orange!90!black,opacity=.5,draw=orange!90!black, line width=4pt] (.2,.4) -- (.8,.4) -- (.8,.8) -- (.2,.8) -- cycle;
\draw[orange!90!black,line width=4pt,opacity=.5] (.4,.4) -- (.4,.8);
\draw[orange!90!black,line width=4pt,opacity=.5] (.6,.4) -- (.6,.8);
\draw[orange!90!black,line width=4pt,opacity=.5] (.2,.7) -- (.8,.7);
\draw[orange!90!black,line width=4pt,opacity=.5] (.2,.6) -- (.8,.6);
\draw[step=.2] (0,.2) grid (1,1.2);
\draw (.5,.2) -- (.5,.8);
\draw (.2,.7) -- (.8,.7);
\end{tikzpicture}
}
\subfloat[]{
\begin{tikzpicture}[scale=3.6]
\filldraw[fill=green!50!black,opacity=.5,draw=green!50!black, line width=4pt] (.2,.6) -- (.8,.6) -- (.8,1) -- (.2,1) -- cycle;
\draw[green!50!black,line width=4pt,opacity=.5] (.4,.6) -- (.4,1);
\draw[green!50!black,line width=4pt,opacity=.5] (.6,.6) -- (.6,1);
\draw[green!50!black,line width=4pt,opacity=.5] (.2,.7) -- (.8,.7);
\draw[green!50!black,line width=4pt,opacity=.5] (.2,.8) -- (.8,.8);

\filldraw[fill=cyan!90!black,opacity=.5,draw=cyan!90!black, line width=4pt] (.4,.4) -- (.8,.4) -- (.8,.8) -- (.4,.8) -- cycle;
\draw[cyan!90!black,line width=4pt,opacity=.5] (.6,.4) -- (.6,.8);
\draw[cyan!90!black,line width=4pt,opacity=.5] (.5,.4) -- (.5,.8);
\draw[cyan!90!black,line width=4pt,opacity=.5] (.4,.7) -- (.8,.7);
\draw[cyan!90!black,line width=4pt,opacity=.5] (.4,.6) -- (.8,.6);

\filldraw[fill=magenta!90!black,opacity=.5,draw=magenta!90!black, line width=4pt] (.2,.4) -- (.6,.4) -- (.6,.8) -- (.2,.8) -- cycle;
\draw[magenta!90!black,line width=4pt,opacity=.5] (.4,.4) -- (.4,.8);
\draw[magenta!90!black,line width=4pt,opacity=.5] (.5,.4) -- (.5,.8);
\draw[magenta!90!black,line width=4pt,opacity=.5] (.2,.7) -- (.6,.7);
\draw[magenta!90!black,line width=4pt,opacity=.5] (.2,.6) -- (.6,.6);

\draw[step=.2] (0,.2) grid (1,1.2);
\draw (.5,.2) -- (.5,.8);
\draw (.2,.7) -- (.8,.7);
\end{tikzpicture}
}
\caption{Replacement of an LR B-spline due to the refinement of the LR mesh $\cM_{\m{p}}^\mu$ of bidegree $\m{p} = \pmb{2}$, shown in (a). Let $B$ be the LR B-spline on $\cM_{\m{p}}^\mu$ whose support and associated mesh $\cM_B^{\mu_{B}}$ are highlighted in (a). We insert a horizontal meshline $\gamma$, dashed in (a), traversing $\supp B$. Then, $B$ is replaced by the B-splines $B_1$ and $B_2$, shown in (b). There is a vertical meshline traversing $\supp B_1$ that is not part of $\cM_{B_1}^{\mu_{B_1}}$ due to the smaller support of $B_1$ (the bottom B-spline in (b)). This triggers that $B_1$ is refined as well, and replaced by two other B-splines $B_3, B_4$. The supports and local tensor meshes of the latter are highlighted in (c).
Therefore, in the end, we move from $B$, which was defined on $\cM_{\m{p}}^\mu$, to $B_2, B_3, B_4$ on the new LR mesh $(\cM + \gamma)_{\m{p}}^{\mu+\mu_\gamma}$. }\label{fig:lrgeneration}
\end{figure}
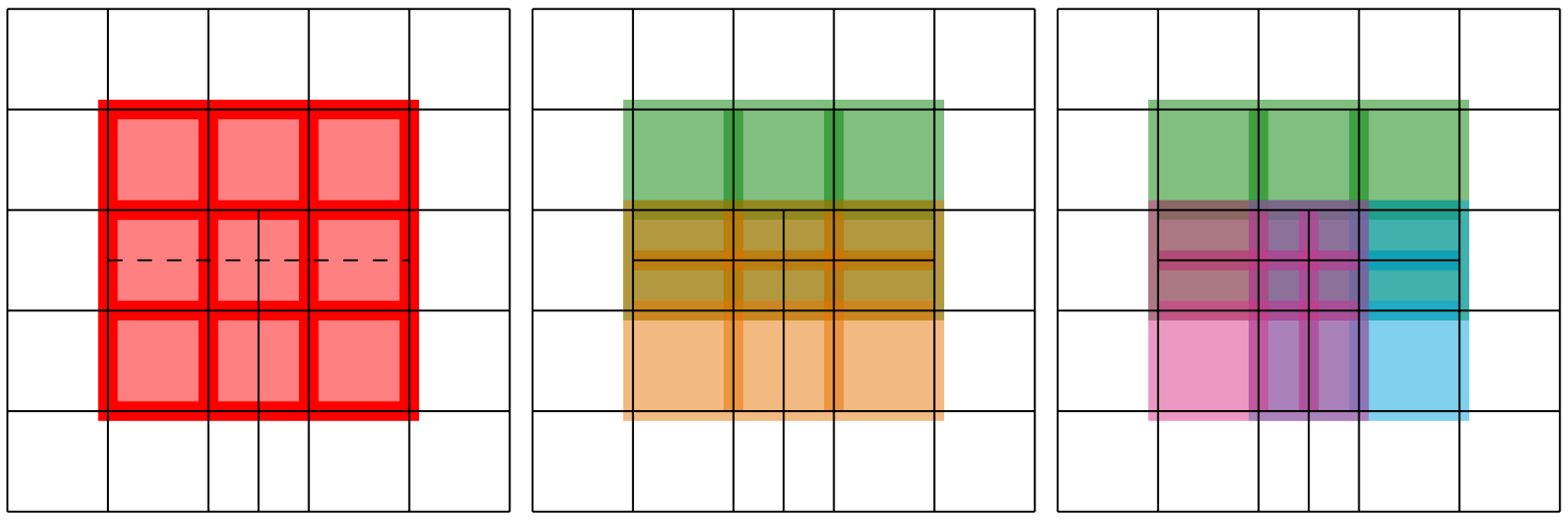

The update of the LR B-spline set induced by the insertion of a new segment in the LR mesh is schematized in Algorithm~\ref{alg:LRgeneration}.
It is important to observe that the insertion of a single segment may trigger multiple refinements of an LR B-spline (through repeated application of the knot insertion procedure), in order to recover the minimal support property.
More precisely, an LR B-spline may be successively refined, not only by the newly inserted meshline, but also by previously inserted meshlines that intersect its support. As a consequence, a single LR B-spline may be replaced by several new LR B-splines during the refinement process, as illustrated in Figure~\ref{fig:lrgeneration}.
In this situation, we say that \textbf{secondary splits} occur in the replacement of the LR B-spline.

A given LR mesh can often be obtained in various ways by changing the order of the segment insertions. The following result ensures that LR B-splines are well defined, regardless of the ordering used to construct the underlying LR mesh; see \cite[Theorem~3.4]{lr}.

\begin{theorem}
The LR B-spline set $\cL[\cM^\mu]$ does not depend on the order of insertion of the meshlines, nor on the order of the subsequent single knot insertions.
\end{theorem}

\begin{figure}[t!]
    \centering
\subfloat[]{
\begin{tikzpicture}[scale = 1.5]
\draw (0,0) -- (1.25,0) -- (1.25,1) -- (0,1) -- cycle;
\draw (0,.55) -- (.75, .55);
\draw (.5, .44) -- (1.25, .44);
\draw (.583, 0) -- (.583, .66);
\draw (.666, .33) -- (.666, 1);
\draw (0,.33) -- (1.25, .33);
\draw (0,.66) -- (1.25, .66);
\draw (.25, 0) -- (.25, 1);
\draw (.5, 0) -- (.5, 1);
\draw (.75, 0) -- (.75, 1);
\draw (1, 0) -- (1, 1);
\draw[white] (0,-.08);
\end{tikzpicture}}~
\raisebox{1.05cm}{\subfloat[]{
\begin{tabular}{c}
\begin{tikzpicture}
\filldraw[draw=red, ultra thick, fill=red!60] (0,.33) -- (.66, .33) -- (.66,1) -- (0,1) -- cycle;
\draw[red, ultra thick] (0,.55) -- (.66, .55);
\draw[red, ultra thick] (0,.66) -- (.66, .66);
\draw[red, ultra thick] (.25,.33) -- (.25, 1);
\draw[red, ultra thick] (.5,.33) -- (.5, 1);
\draw (0,0) -- (1.25,0) -- (1.25,1) -- (0,1) -- cycle;
\draw (0,.55) -- (.75, .55);
\draw (.5, .44) -- (1.25, .44);
\draw (.583, 0) -- (.583, .66);
\draw (.666, .33) -- (.666, 1);
\draw (0,.33) -- (1.25, .33);
\draw (0,.66) -- (1.25, .66);
\draw (.25, 0) -- (.25, 1);
\draw (.5, 0) -- (.5, 1);
\draw (.75, 0) -- (.75, 1);
\draw (1, 0) -- (1, 1);
\end{tikzpicture},
\begin{tikzpicture}
\filldraw[draw=red, ultra thick, fill=red!60] (.25,.33) -- (.75, .33) -- (.75,1) -- (.25,1) -- cycle;
\draw[red, ultra thick] (.25,.55) -- (.75, .55);
\draw[red, ultra thick] (.25,.66) -- (.75, .66);
\draw[red, ultra thick] (.5,.33) -- (.5, 1);
\draw[red, ultra thick] (.66,.33) -- (.66, 1);
\draw (0,0) -- (1.25,0) -- (1.25,1) -- (0,1) -- cycle;
\draw (0,.55) -- (.75, .55);
\draw (.5, .44) -- (1.25, .44);
\draw (.583, 0) -- (.583, .66);
\draw (.666, .33) -- (.666, 1);
\draw (0,.33) -- (1.25, .33);
\draw (0,.66) -- (1.25, .66);
\draw (.25, 0) -- (.25, 1);
\draw (.5, 0) -- (.5, 1);
\draw (.75, 0) -- (.75, 1);
\draw (1, 0) -- (1, 1);
\end{tikzpicture},
\begin{tikzpicture}
\filldraw[draw=red, ultra thick, fill=red!60] (.5,.33) -- (1, .33) -- (1,1) -- (.5,1) -- cycle;
\draw[red, ultra thick] (.5,.44) -- (1, .44);
\draw[red, ultra thick] (.5,.66) -- (1, .66);
\draw[red, ultra thick] (.666,.33) -- (.666, 1);
\draw[red, ultra thick] (.75,.33) -- (.75, 1);
\draw (0,0) -- (1.25,0) -- (1.25,1) -- (0,1) -- cycle;
\draw (0,.55) -- (.75, .55);
\draw (.5, .44) -- (1.25, .44);
\draw (.583, 0) -- (.583, .66);
\draw (.666, .33) -- (.666, 1);
\draw (0,.33) -- (1.25, .33);
\draw (0,.66) -- (1.25, .66);
\draw (.25, 0) -- (.25, 1);
\draw (.5, 0) -- (.5, 1);
\draw (.75, 0) -- (.75, 1);
\draw (1, 0) -- (1, 1);
\end{tikzpicture},
\begin{tikzpicture}
\filldraw[draw=red, ultra thick, fill=red!60] (.666,.33) -- (1.25, .33) -- (1.25,1) -- (.666,1) -- cycle;
\draw[red, ultra thick] (.666,.44) -- (1.25, .44);
\draw[red, ultra thick] (.666,.66) -- (1.25, .66);
\draw[red, ultra thick] (.75,.33) -- (.75, 1);
\draw[red, ultra thick] (1,.33) -- (1, 1);
\draw (0,0) -- (1.25,0) -- (1.25,1) -- (0,1) -- cycle;
\draw (0,.55) -- (.75, .55);
\draw (.5, .44) -- (1.25, .44);
\draw (.583, 0) -- (.583, .66);
\draw (.666, .33) -- (.666, 1);
\draw (0,.33) -- (1.25, .33);
\draw (0,.66) -- (1.25, .66);
\draw (.25, 0) -- (.25, 1);
\draw (.5, 0) -- (.5, 1);
\draw (.75, 0) -- (.75, 1);
\draw (1, 0) -- (1, 1);
\end{tikzpicture},
\begin{tikzpicture}
\filldraw[draw=red, ultra thick, fill=red!60] (0,0) -- (.583, 0) -- (.583,.66) -- (0,.66) -- cycle;
\draw[red, ultra thick] (0,.33) -- (.583, .33);
\draw[red, ultra thick] (0,.55) -- (.583, .55);
\draw[red, ultra thick] (.25,0) -- (.25, .66);
\draw[red, ultra thick] (.5,0) -- (.5, .66);
\draw (0,0) -- (1.25,0) -- (1.25,1) -- (0,1) -- cycle;
\draw (0,.55) -- (.75, .55);
\draw (.5, .44) -- (1.25, .44);
\draw (.583, 0) -- (.583, .66);
\draw (.666, .33) -- (.666, 1);
\draw (0,.33) -- (1.25, .33);
\draw (0,.66) -- (1.25, .66);
\draw (.25, 0) -- (.25, 1);
\draw (.5, 0) -- (.5, 1);
\draw (.75, 0) -- (.75, 1);
\draw (1, 0) -- (1, 1);
\end{tikzpicture}\\
\begin{tikzpicture}
\filldraw[draw=red, ultra thick, fill=red!60] (.25,0) -- (.75, 0) -- (.75,.66) -- (.25,.66) -- cycle;
\draw[red, ultra thick] (.25,.33) -- (.75, .33);
\draw[red, ultra thick] (.25,.55) -- (.75, .55);
\draw[red, ultra thick] (.5,0) -- (.5, .66);
\draw[red, ultra thick] (.583,0) -- (.583, .66);
\draw (0,0) -- (1.25,0) -- (1.25,1) -- (0,1) -- cycle;
\draw (0,.55) -- (.75, .55);
\draw (.5, .44) -- (1.25, .44);
\draw (.583, 0) -- (.583, .66);
\draw (.666, .33) -- (.666, 1);
\draw (0,.33) -- (1.25, .33);
\draw (0,.66) -- (1.25, .66);
\draw (.25, 0) -- (.25, 1);
\draw (.5, 0) -- (.5, 1);
\draw (.75, 0) -- (.75, 1);
\draw (1, 0) -- (1, 1);
\end{tikzpicture},
\begin{tikzpicture}
\filldraw[draw=red, ultra thick, fill=red!60] (.25,0) -- (1, 0) -- (1,1) -- (.25,1) -- cycle;
\draw[red, ultra thick] (.25,.33) -- (1, .33);
\draw[red, ultra thick] (.25,.66) -- (1, .66);
\draw[red, ultra thick] (.5,0) -- (.5, 1);
\draw[red, ultra thick] (.75,0) -- (.75, 1);
\draw (0,0) -- (1.25,0) -- (1.25,1) -- (0,1) -- cycle;
\draw (0,.55) -- (.75, .55);
\draw (.5, .44) -- (1.25, .44);
\draw (.583, 0) -- (.583, .66);
\draw (.666, .33) -- (.666, 1);
\draw (0,.33) -- (1.25, .33);
\draw (0,.66) -- (1.25, .66);
\draw (.25, 0) -- (.25, 1);
\draw (.5, 0) -- (.5, 1);
\draw (.75, 0) -- (.75, 1);
\draw (1, 0) -- (1, 1);
\end{tikzpicture},
\begin{tikzpicture}
\filldraw[draw=red, ultra thick, fill=red!60] (.5,0) -- (1, 0) -- (1,.66) -- (.5,.66) -- cycle;
\draw[red, ultra thick] (.5,.33) -- (1, .33);
\draw[red, ultra thick] (.5,.44) -- (1, .44);
\draw[red, ultra thick] (.583,0) -- (.583, .66);
\draw[red, ultra thick] (.75,0) -- (.75, .66);
\draw (0,0) -- (1.25,0) -- (1.25,1) -- (0,1) -- cycle;
\draw (0,.55) -- (.75, .55);
\draw (.5, .44) -- (1.25, .44);
\draw (.583, 0) -- (.583, .66);
\draw (.666, .33) -- (.666, 1);
\draw (0,.33) -- (1.25, .33);
\draw (0,.66) -- (1.25, .66);
\draw (.25, 0) -- (.25, 1);
\draw (.5, 0) -- (.5, 1);
\draw (.75, 0) -- (.75, 1);
\draw (1, 0) -- (1, 1);
\end{tikzpicture},
\begin{tikzpicture}
\filldraw[draw=red, ultra thick, fill=red!60] (.583,0) -- (1.25, 0) -- (1.25,.66) -- (.583,.66) -- cycle;
\draw[red, ultra thick] (.583,.33) -- (1.25, .33);
\draw[red, ultra thick] (.583,.44) -- (1.25, .44);
\draw[red, ultra thick] (.75,0) -- (.75, .66);
\draw[red, ultra thick] (1,0) -- (1, .66);
\draw (0,0) -- (1.25,0) -- (1.25,1) -- (0,1) -- cycle;
\draw (0,.55) -- (.75, .55);
\draw (.5, .44) -- (1.25, .44);
\draw (.583, 0) -- (.583, .66);
\draw (.666, .33) -- (.666, 1);
\draw (0,.33) -- (1.25, .33);
\draw (0,.66) -- (1.25, .66);
\draw (.25, 0) -- (.25, 1);
\draw (.5, 0) -- (.5, 1);
\draw (.75, 0) -- (.75, 1);
\draw (1, 0) -- (1, 1);
\end{tikzpicture}
\end{tabular}
}}~
\subfloat[]{
\begin{tikzpicture}[scale=1.5]
\filldraw[draw = blue, ultra thick, fill=blue!50] (.5,.33) -- (.75,.33) -- (.75, .66) -- (.5, .66) -- cycle;
\draw[blue, ultra thick] (.5, .44) -- (.75, .44);
\draw[blue, ultra thick] (.5, .55) -- (.75, .55);
\draw[blue, ultra thick] (.583, .33) -- (.583, .66);
\draw[blue, ultra thick] (.666, .33) -- (.666, .66);
\draw (0,0) -- (1.25,0) -- (1.25,1) -- (0,1) -- cycle;
\draw (0,.55) -- (.75, .55);
\draw (.5, .44) -- (1.25, .44);
\draw (.583, 0) -- (.583, .66);
\draw (.666, .33) -- (.666, 1);
\draw (0,.33) -- (1.25, .33);
\draw (0,.66) -- (1.25, .66);
\draw (.25, 0) -- (.25, 1);
\draw (.5, 0) -- (.5, 1);
\draw (.75, 0) -- (.75, 1);
\draw (1, 0) -- (1, 1);
\draw[white] (0,-.08);
\end{tikzpicture}
}
\caption{Consider the LR mesh $\cM_{\m{p}}^\mu$, shown in (a), with meshlines of multiplicity $1$ and bidegree $\m{p} = \pmb{2}$. The supports of all biquadratic LR B-splines on $\cM_{\m{p}}^\mu$ are depicted in (b) and the support of a biquadratic MS B-spline on $\cM_{\m{p}}^\mu$ that is not an LR B-spline in (c).}\label{fig:exmsnotlr}
\end{figure}
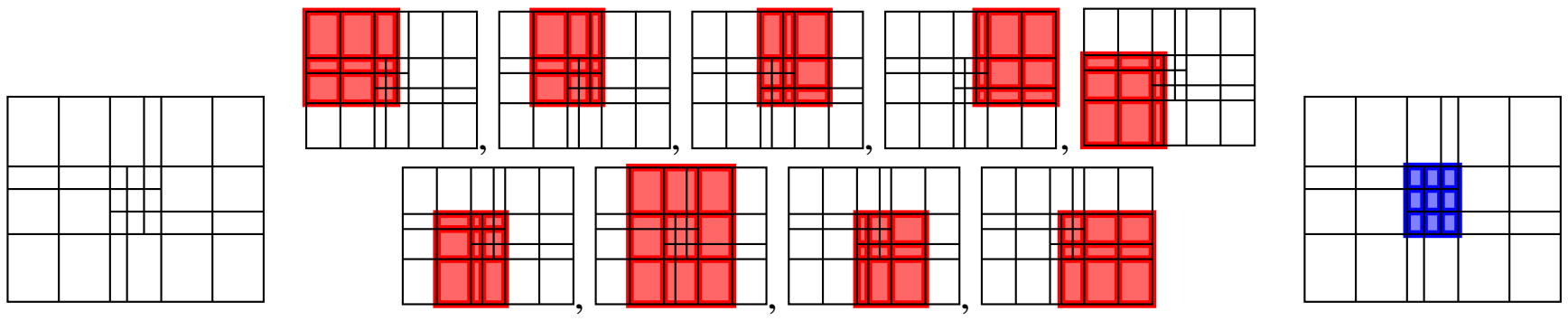

It may happen that MS B-splines cannot be obtained by refinement of previously defined LR B-splines during the refinement process of the LR mesh. As a consequence, the LR B-spline set $\cL[\cM^\mu]$ can be  a proper subset of the MS B-spline set $\cB[\cM^\mu]$ on $\cM^\mu$.
An example is shown in Figure~\ref{fig:exmsnotlr}.

The definitions of LR B-splines and MS B-splines on a given LR mesh can be unified within a graph theory framework as follows; see \cite[Section~2.3]{someproperties}.

\begin{definition}
Given an LR mesh $\cM_{\m{p}}^\mu$, the \textbf{B-spline graph} of $\cM_{\m{p}}^\mu$ is the directed graph $G[\cM_{\m{p}}^\mu] \defeq (\cV^G, \cE^G)$, where the set of vertices $\cV^G$ consists of all B-splines of bidegree $\m{p}$ with support on $\cM_{\m{p}}^\mu$, and $(B_1, B_2) \in \cE^G$ if $B_2$ is obtained from $B_1$ by a single knot insertion during the LR mesh refinement process.
Moreover, we use the following terminology:
\begin{itemize}
\item a \textbf{path} in  $G[\cM_{\m{p}}^\mu]$ is a sequence of vertices $(B_{i_1}, \ldots, B_{i_n})$ such that $(B_{i_j}, B_{i_{j+1}}) \in \cE^G$ for all $j = 1, \ldots, n-1$;
\item a vertex $B_{i_1}$ is a \textbf{predecessor} of a vertex $B_{i_n}$ if there exists a path from $B_{i_1}$ to $B_{i_n}$; in such a case, $B_{i_n}$ is called a \textbf{successor} of $B_{i_1}$;
\item a \textbf{source} is a vertex with no predecessors, while a \textbf{sink} is a vertex with no successors;
\item a vertex $B \in \cV^G$ is \textbf{reachable} from $\cW^G \subseteq \cV^G$, if there exists a path from some vertex in $\cW^G$ to $B$, or if $B \in \cW^G$.
\end{itemize}
\end{definition}

The sources of the B-spline graph are either the tensor-product B-splines defined on the initial tensor spline mesh or MS B-splines that are not obtained from refinements of previously defined LR B-splines. Note that sources can also be sinks.
For any LR mesh $\cM_{\m{p}}^\mu$, the graph $G[\cM_{\m{p}}^\mu]$ is finite, directed, and acyclic. It therefore induces a partial order relation in the set of vertices.

\begin{figure}[t!]
\centering
\begin{tikzpicture}
\matrix (m)[matrix of math nodes,column sep=0.4em, row sep=3em]{
  \includegraphics[width = 0.1\textwidth, page = 2]{bsplinegraph}
  \pgfmatrixnextcell
  \pgfmatrixnextcell
  \pgfmatrixnextcell
\includegraphics[width = 0.1\textwidth, page = 1]{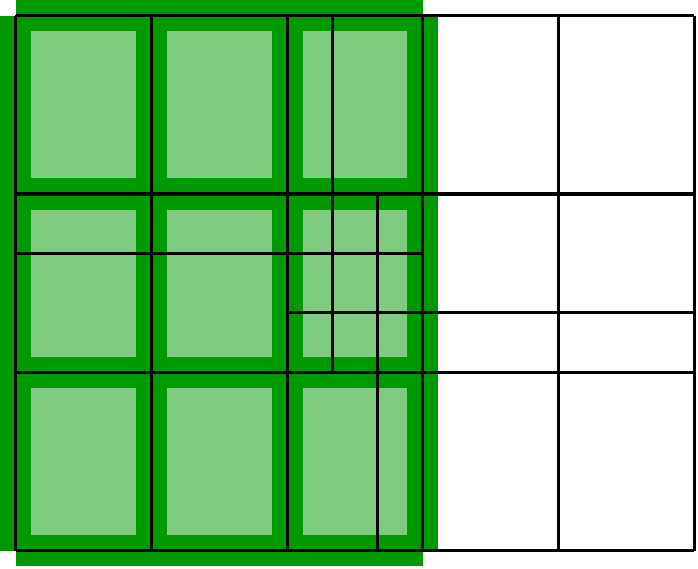}
\pgfmatrixnextcell
\includegraphics[width = 0.1\textwidth, page = 3]{bsplinegraph}
  \pgfmatrixnextcell
    \pgfmatrixnextcell
        \pgfmatrixnextcell
  \includegraphics[width = 0.1\textwidth, page = 16]{bsplinegraph}
\\
  \pgfmatrixnextcell
  \includegraphics[width = 0.1\textwidth, page = 4]{bsplinegraph}
  \pgfmatrixnextcell
  \includegraphics[width = 0.1\textwidth, page = 5]{bsplinegraph}
  \pgfmatrixnextcell
  \pgfmatrixnextcell
  \pgfmatrixnextcell
  \includegraphics[width = 0.1\textwidth, page = 14]{bsplinegraph}
  \pgfmatrixnextcell
  \includegraphics[width = 0.1\textwidth, page = 15]{bsplinegraph}
 \\
\includegraphics[width = 0.1\textwidth, page = 6]{bsplinegraph}
 \pgfmatrixnextcell
 \includegraphics[width = 0.1\textwidth, page = 7]{bsplinegraph}
 \pgfmatrixnextcell
 \includegraphics[width = 0.1\textwidth, page = 8]{bsplinegraph}
 \pgfmatrixnextcell
 \includegraphics[width = 0.1\textwidth, page = 9]{bsplinegraph}
 \pgfmatrixnextcell
 \includegraphics[width = 0.1\textwidth, page = 10]{bsplinegraph}
 \pgfmatrixnextcell
 \includegraphics[width = 0.1\textwidth, page = 12]{bsplinegraph}
 \pgfmatrixnextcell
 \includegraphics[width = 0.1\textwidth, page = 11]{bsplinegraph}
 \pgfmatrixnextcell
 \includegraphics[width = 0.1\textwidth, page = 13]{bsplinegraph}\\};

\draw[-stealth, thick] (m-1-4) -- (m-2-2);
\draw[-stealth, thick] (m-1-4) -- (m-2-3);
\draw[-stealth, thick] (m-1-5) -- (m-2-6);
\draw[-stealth, thick] (m-1-5) -- (m-2-7);
\draw[-stealth, thick] (m-2-2) -- (m-3-1);
\draw[-stealth, thick] (m-2-2) -- (m-3-2);
\draw[-stealth, thick] (m-2-3) -- (m-3-3);
\draw[-stealth, thick] (m-2-3) -- (m-3-4);

\draw[-stealth, thick] (m-2-6) -- (m-3-5);
\draw[-stealth, thick] (m-2-6) -- (m-3-6);
\draw[-stealth, thick] (m-2-7) -- (m-3-7);
\draw[-stealth, thick] (m-2-7) -- (m-3-8);
\end{tikzpicture}
\caption{The B-spline graph of the LR mesh considered in Figure~\ref{fig:exmsnotlr}.
The sources with the green support are the tensor-product B-splines defined on the initial tensor spline mesh.
The sinks with the red support are reachable from these sources, while the sink with the blue support is not.
The intermediate vertices have orange support.
The LR B-spline set $\cL$ contains the red sinks together with the isolated green source (which is a trivially reachable sink). The MS B-spline set $\cB$ additionally contains the blue sink.}\label{fig:bsplinegraph}
\end{figure}

Using this graph structure, we have the following unifying characterization of MS B-splines and LR B-splines \cite[Theorem~2.26]{someproperties}. See Figure~\ref{fig:bsplinegraph} for an example.

\begin{theorem}
Let $\cM_{\m{p}}^\mu$ be a given LR mesh and let $G[\cM_{\m{p}}^\mu]$ be the associated B-spline graph.
Moreover, let $\cW^G$ be the set of sources that are the tensor-product B-splines defined on the initial tensor spline mesh.
Then,
\begin{itemize}
\item the set of sinks agrees with the MS B-spline set $\cB[\cM_{\m{p}}^\mu]$;
\item the set of sinks reachable from $\cW^G$ agrees with the LR B-spline set $\cL[\cM_{\m{p}}^\mu]$.
\end{itemize}
\end{theorem}

Thanks to the local tensor-product structure, several properties of univariate and tensor-product B-splines are preserved by LR B-splines.

\begin{proposition}
Let $\cM_{\m{p}}^\mu$ be a given LR mesh. For any
$B \defeq B[\m{x}_B, \m{y}_B] \in \cL[\cM_{\m{p}}^\mu]$, the following properties hold.
\begin{itemize}
\item \textbf{Local support:} We have
\begin{equation*}
\supp B = [x_{B,1}, x_{B,p_1+2}] \times [y_{B,1}, y_{B,p_2+2}].
\end{equation*}
\item \textbf{Positivity:} We have
\begin{equation*}
B(x,y) > 0, \quad (x,y) \in (x_{B,1}, x_{B,p_1+2}) \times (y_{B,1}, y_{B,p_2+2}).
\end{equation*}
\item \textbf{Piecewise polynomial structure:} For $i = 1, \ldots, p_1+1$ and $j = 1, \ldots, p_2+1$, it holds
\begin{equation*}
B\big|_{[x_{B,i}, x_{B,i+1}) \times [y_{B,j}, y_{B,j+1})} \in \PP_{\m{p}}.
\end{equation*}
\item \textbf{Affine invariance:} For any $\alpha_1, \alpha_2 > 0$ and $\beta_1, \beta_2 \in \RR$, it holds
\begin{equation*}
B[\alpha_1 \m{x}_B + \beta_1,\, \alpha_2 \m{y}_B + \beta_2](\alpha_1 x + \beta_1,\, \alpha_2 y + \beta_2)
= B[\m{x}_B, \m{y}_B](x,y).
\end{equation*}

\item \textbf{Smoothness:} $B$ is $C^{p_k - \mu_B(\gamma)}$-continuous across any $k$-meshline $\gamma$ in $\cM_B^{\mu_{B}}$.

\item \textbf{Derivative formula:} Since the derivative formula for univariate B-splines (see \eqref{eq:Bsplinederivative}) is local, it extends componentwise to partial derivatives of $B$.
%
\end{itemize}
\end{proposition}

Evaluation algorithms, such as the recursive definition (Definition~\ref{def:bspline}) and the \emph{de Boor algorithm} \cite{debooralgo}, extend naturally to LR B-splines due to their locality.
On the other hand, the remaining properties listed in Proposition~\ref{prop:bsplineadvancedproperties} for the univariate case, namely partition of unity, monomial reproduction, linear independence, do not directly extend to LR B-splines.

Linear independence cannot be guaranteed in general. Both MS B-splines and LR B-splines may be linearly dependent on certain pathological LR meshes and establishing sufficient conditions for linear dependence of either of the two sets is, in general, an open problem.
On the other hand, local linear independence can be characterized in terms of the support distribution on the LR mesh.
Linear dependence and local linear independence will be analyzed in the next section.

In the remainder of this section, we focus on the other two properties.
From now on, except when explicitly stated otherwise,
we assume the initial tensor spline mesh to be open.
Clearly, $\PP_\m{p} \subseteq \spn \cL[\cM_{\m{p}}^\mu]$  because the various steps of the LR B-spline construction produce nested spaces starting from the initial tensor-product spline space, which contains $\PP_\m{p}$.
Partition of unity on general unstructured meshes can be achieved by means of proper scaling weights; see \cite[Section~7]{lr}.

\begin{theorem}\label{thm:lrpou}
Let $\cL$ be the LR B-spline set defined on an open LR mesh $\cM^\mu$. It holds
\begin{equation*}
\sum_{B \in \cL} w_B \, B = 1,
\end{equation*}
where the weights can be computed recursively during the LR B-spline construction as follows.
If $\cM^\mu$ is a tensor spline mesh, we have $w_B=1$, for all $B \in \cL$.
Otherwise, let $\cL$ be obtained from another LR B-spline set $\cL'$ (Algorithm~\ref{alg:LRgeneration}). Then,
\begin{equation*}
w_B = \sum_{B' \in \cL'} \alpha_{B,B'} \, w_{B'},
\end{equation*}
where the values $\alpha_{B,B'} \geq 0$ are those appearing in the relation
\begin{equation}\label{eq:weights}
B' = \sum_{B \in \cL} \alpha_{B,B'} \, B,
\end{equation}
derived from the repeated application of the univariate knot insertion relation (see Theorem~\ref{thm:knotinsertion}) according to the refinement process from $\cL'$ to $\cL$.
\end{theorem}

Note that, due to the possible lack of linear independence, scaling weights might not be uniquely defined. However, the procedure described in Theorem~\ref{thm:lrpou} is uniquely defined.
Starting from weights equal to $1$ on the initial tensor spline mesh, these weights are updated by accumulating the knot insertion coefficients of the LR B-splines replaced during the refinement process (see also Remark~\ref{rmk:pou1d}).
In particular, if the same B-spline appears in the refinement of several LR B-splines, its associated weight is given by the cumulative sum of the related knot insertion coefficients.

Algorithm~\ref{alg:scaledLRgeneration} extends Algorithm~\ref{alg:LRgeneration} in order to update the LR B-splines and their weights after each refinement step of the mesh. The scaled LR B-spline set is denoted by $w\cL$.

\begin{algorithm}[t!]
\caption{Scaled LR B-spline update -- \scalebox{0.9}{$[\cL, \m{w}] = \texttt{scaled\_lr}(\cL', \m{w}', \cM^\mu)$}}\label{alg:scaledLRgeneration}

\setlength{\columnsep}{-5cm}
\nonl\begin{multicols}{2}
\textbf{Input:}\vspace{0.1cm}
\begin{itemize}
\item[~~$\cL'\colon$] ~~~~~current LR B-spline set\vspace{0.1cm}
\item[\,~~$\m{w}'\colon$] ~~~~~vector of weights associated with $\cL'$\vspace{0.1cm}
\item[$\cM^\mu\colon$] ~~~~~refinement of the LR mesh of $\cL'$
\end{itemize}\vfill\null
\columnbreak

\hspace{3.5cm}
\textbf{Output:}\vspace{0.1cm}

\hspace{3.65cm}\begin{minipage}{0.4\textwidth}
\begin{itemize}
\item[$\cL\colon$] ~~LR B-spline set on $\cM^\mu$\vspace{0.1cm}
\item[\,$\m{w}\colon$] ~~vector of updated weights
\end{itemize}
\end{minipage}
\vfill\null
\end{multicols}\vspace{-1cm}

\rule{0.965\textwidth}{0.025cm}

\linespread{1.35}\selectfont

Initialize $\cL \leftarrow \cL'$\;
Initialize $\m{w} \leftarrow \m{w}'$\;

\nonl\While{there exists $B \in \cL$ that does not have minimal support on $\cM^\mu$}{
    Select a meshline $\gamma$ of $\cM$ that breaks the minimal support property of $B$\;

    Compute the knot insertion related to $\gamma$: $B = \alpha_{1,B} B_1 + \alpha_{2,B} B_2$\;
    \nonl\For{$k = 1, 2$}{
 \nonl\uIf{$B_k \in \cL$}{
    Update weight in $\m{w}$: $w_{B_k} \leftarrow w_{B_k} + w_{B} \alpha_{k,B}$\;}
    \nonl\Else{
    Define weight: $w_{B_k} \defeq w_{B} \alpha_{k,B}$\;
    Add $B_k$ to $\cL$ and $w_{B_k}$ to $\m{w}$: $\cL \leftarrow \cL \cup \{B_k\}$ and $\m{w} \leftarrow \m{w} \cup \{w_{B_k}\}$\;
    }}
    Remove $B$ from $\cL$ and $w_{B}$ from $\m{w}$: $\cL \leftarrow \cL \setminus \{B\}$ and $\m{w} \leftarrow \m{w} \setminus \{w_{B}\}$\;
}
\end{algorithm}

\begin{remark}
An LR spline surface in $\RR^d$ ($d \geq 3$), represented in terms of the scaled LR B-spline set $w\cL$, i.e.,
\begin{equation*}
\mathbf{S}(x, y) \defeq \sum_{B \in \cL} \m{P}_B \, w_B B(x, y),
\quad \m{P}_B \in \RR^d,\quad (x,y) \in \Omega,
\end{equation*}
lies in the convex hull of the control points $\{\m{P}_B\}_{B \in \cL}$, similar to the curves in the univariate case; see Remark~\ref{rmk:convexhull1d}.
\end{remark}

The same strategy can be used for the reproduction of monomials and, more generally, of LR splines (i.e., linear combinations of LR B-splines) defined on intermediate meshes arising during the iterative construction of the LR mesh.
Indeed, on the initial tensor spline mesh, the coefficients for representing monomials are known from the univariate theory; see Proposition~\ref{prop:bsplineadvancedproperties}. Then, these coefficients can be propagated to finer meshes, yielding representations of the same function in terms of LR B-splines defined on the final mesh.

\begin{theorem}\label{thm:lrrep}
Let $\cL$ be an LR B-spline set obtained from another LR B-spline set $\cL'$ (Algorithm~\ref{alg:LRgeneration}).
Suppose $f\in\spn \cL'$. Then,
\begin{equation*}
f = \sum_{B' \in \cL'} c_{B'} \, B'=\sum_{B \in \cL} c_B \, B,
\end{equation*}
where the coefficients $c_B$ are given by
\begin{equation*}
c_B = \sum_{B' \in \cL'} \alpha_{B,B'} \, c_{B'},
\end{equation*}
and the values $\alpha_{B,B'}$ are defined in \eqref{eq:weights}.
\end{theorem}

\begin{remark}
Theorems~\ref{thm:lrpou} and~\ref{thm:lrrep} show that both partition of unity and monomial reproduction are achieved through recursive procedures. In general, there is no closed-form expression, and the construction depends on the refinement history from the initial tensor spline mesh.
However, if additional constraints are imposed on the LR mesh, namely the \emph{N$_2$S property} (see Section~\ref{sec:loclinind} for its definition), then the situation simplifies significantly: local linear independence is ensured and the resulting LR B-splines retain all key features of tensor-product B-splines. In particular, all weights $w_B$ are equal to $1$ and the monomial reproduction takes the same form as in the tensor-product case, i.e., for $0 \leq i_k \leq p_k$, it holds
\begin{equation}\label{eq:explicitmonomials}
x^{i_1} y^{i_2}
= \sum_{B \in \cL}
x_{B,p_1}^{\overbrace{\scalebox{0.75}{$*\cdots*$}}^{i_1 \text{ times}}}
y_{B,p_2}^{\overbrace{\scalebox{0.75}{$*\cdots*$}}^{i_2 \text{ times}}}
\, B(x, y),
\quad (x,y) \in \Omega.
\end{equation}
Note that the coefficients in \eqref{eq:explicitmonomials} are simply obtained by multiplying the corresponding coefficients from the univariate case; see Proposition~\ref{prop:bsplineadvancedproperties}.
\end{remark}

The existence of an algorithmic procedure to define scaling weights ensuring partition of unity,
coupled with a lower probability to generate linear dependence relations in the function set (this will be analyzed in the next section), are the main reasons why LR B-splines are preferred over MS B-splines on a given LR mesh in practical applications.

The multivariate definitions of MS B-splines and LR B-splines in $\RR^d$ extend naturally from the bivariate case. The only additional requirement is the enforcement of LR rule 4 when inserting new $(d-1)$-boxes.
Despite the higher dimensionality, a $k$-meshrectangle traversing the support of a B-spline still induces a univariate knot insertion in the $k$th local knot vector. Consequently, the construction of LR meshes and LR B-splines follows the same principles as in the bivariate setting and all properties and procedures discussed above extend to the general multivariate case.

\section{Spanning Properties and Linear Independence}\label{sec:spanning-linind}
In this section, we analyze three important aspects related to the span of MS B-splines and LR B-splines.
\begin{itemize}
\item \textbf{Completeness:} MS B-splines, and thus LR B-splines, are spline functions on a given LR mesh. However, their span  does not always coincide with the full spline space on the same mesh. We provide sufficient conditions ensuring completeness.
\item \textbf{Linear dependence:} LR meshes allow for a high flexibility in segment insertion during refinement. This may lead to pathological configurations in which linear dependence relations arise among MS B-splines or LR B-splines. We derive necessary conditions for the occurrence of such dependencies and discuss strategies to detect and eliminate them.
\item \textbf{Local linear independence:} It is possible to construct non-tensor LR meshes for which LR B-splines are locally linearly independent. We provide both analytical and geometric characterizations of this property.
\end{itemize}

We mainly focus on the univariate and bivariate settings. However, some statements naturally extend to the general multivariate case.

\subsection{Completeness: The Hand-in-Hand Principle}
In the univariate setting, the B-spline set $\cB_{p,\m{x}}$ defined on the knot vector $\m{x}$ spans the full spline space $\SS[\pmb{\tau}_p^{\mu}]$ defined on the partition identified by $\pmb{\tau}_p^{\mu}$. Therefore, the set $\cB_{p,\m{x}}$ is complete.
This result extends directly to tensor-product spline spaces.

Any LR mesh $\cM^\mu$ is constructed through a sequence of refinements starting from an initial tensor spline mesh. On this initial mesh, tensor-product B-splines span the full spline space. With the aim of preserving this property throughout the refinement process, we introduce the hand-in-hand principle.

\begin{definition}
Let $(\cM')^{\mu'}$ be a given LR mesh and let $\cL'$ be the corresponding LR B-spline set (resp., let $\cB'$ be the corresponding MS B-spline set). Suppose $\spn \cL' = \SS[(\cM')^{\mu'}]$ (resp., $\spn \cB' = \SS[(\cM')^{\mu'}]$).
Let $\gamma$ be a segment to be inserted. Define the refined mesh $\cM^\mu$ by
\begin{equation*}
\cM \defeq \cM' + \gamma,
\quad
\mu \defeq \mu' + \mu_\gamma,
\end{equation*}
and
let $\cL$ be the LR B-spline set (resp., let $\cB$ be the MS B-spline set) on $\cM^\mu$.
We say that $\cM^\mu$ goes \textbf{LR-wise} (resp., \textbf{MS-wise}) \textbf{hand-in-hand} with $(\cM')^{\mu'}$ if $\spn \cL = \SS[\cM^\mu]$ (resp., $\spn \cB = \SS[\cM^\mu]$).
\end{definition}

Since $\cL \subseteq \cB$, if $\cM^\mu$ goes LR-wise hand-in-hand with $(\cM')^{\mu'}$, then it also goes MS-wise hand-in-hand. The converse does not hold in general.
The completeness throughout the construction of an LR mesh is guaranteed if all intermediate meshes satisfy the hand-in-hand property.

The following result characterizes the hand-in-hand property in terms of the dimension of the spline spaces involved when inserting  a segment in a given mesh; see \cite[Theorem~5.10]{lr}.
\begin{theorem}\label{thm:handinhand}
Let $(\cM')^{\mu'}$ be a given LR mesh and let $\cL'$ be the corresponding LR B-spline set (resp., let $\cB'$ be the corresponding MS B-spline set). Suppose $\spn \cL' = \SS[(\cM')^{\mu'}]$ (resp., $\spn \cB' = \SS[(\cM')^{\mu'}]$).
Let $\gamma$ be a segment to be inserted in the $k$th direction. Define the refined mesh $\cM^\mu$ by
\begin{equation*}
\cM \defeq \cM' + \gamma,
\quad
\mu \defeq \mu' + \mu_\gamma,
\end{equation*}
and
let $\cL$ be the LR B-spline set (resp., let $\cB$ be the MS B-spline set) on $\cM^\mu$. Let $\pmb{\tau}_{p_{3-k}}^{\tilde{\mu}_{3-k}}$ be the global knot vector induced on $\gamma$, with modified multiplicity $\tilde{\mu}_{3-k}$ as in Theorem~\ref{thm:splinedimincrease2d}.
Define
\begin{equation*}
\cL[\gamma] \defeq \{ B \in \cL : \gamma \text{ traverses } \mathring{\supp B} \}
\end{equation*}
\begin{equation*}
\bigl(\, \text{resp., }
\cB[\gamma] \defeq \{ B \in \cB : \gamma \text{ traverses } \mathring{\supp B} \}\, \bigr),
\end{equation*}
and for each $B \in \cL[\gamma]$ (resp., $B \in \cB[\gamma]$), let $B_\gamma$ denote the restriction to $\gamma$ of the univariate B-spline in the $(3-k)$th direction associated with $B$.
Then, the refined mesh $\cM^\mu$ goes LR-wise (resp., MS-wise) hand-in-hand with $(\cM')^{\mu'}$ if and only if
\begin{equation*}
\spn \{ B_\gamma : B \in \cL[\gamma] \}
= \SS[\pmb{\tau}_{p_{3-k}}^{\tilde{\mu}_{3-k}}]
\end{equation*}
\begin{equation*}
\bigl(\, \text{resp., }
\spn \{ B_\gamma : B \in \cB[\gamma] \}
= \SS[\pmb{\tau}_{p_{3-k}}^{\tilde{\mu}_{3-k}}]\, \bigr).
\end{equation*}
\end{theorem}

Since each segment insertion must refine at least one LR B-spline, we have that $\cL[\gamma]$ is non-empty. If the dimension of the spline space increases by one, then the hand-in-hand property is automatically satisfied as stated in the following corollary; see \cite[Theorem~5.2]{lr}.

\begin{corollary}\label{cor:handinhand}
Let $(\cM')^{\mu'}$ be a given LR mesh and let $\cM^\mu$ be obtained by inserting a segment $\gamma$. If
\begin{equation*}
\dim \SS[\cM^\mu] = \dim \SS[(\cM')^{\mu'}] + 1,
\end{equation*}
then $\cM^\mu$ goes LR-wise, and thus also MS-wise, hand-in-hand with $(\cM')^{\mu'}$.
\end{corollary}

\begin{figure}[t!]
\centering
\subfloat[]{
\begin{tikzpicture}[scale = 1]
\draw[thick, step = 1] (0, 1) grid (3, 4);
\draw[thick] (0, 2.5) -- (2, 2.5);
\draw[thick] (1.5, 2) -- (1.5, 4);
\draw[thick] (1.75, 1) -- (1.75, 2.5);
\draw[thick, red] (1.5, 2.25) --node[above]{$\gamma$} (3, 2.25);
\draw[dashed] (1.5, 3) -- (1.5, 4.25);
\draw[dashed] (1.75, 2.5) -- (1.75, 4.25);
\draw[dashed] (2, 3) -- (2, 4.25);
\draw[dashed] (3, 3) -- (3, 4.25);
\draw[thick, red] (1.5, 4.25) --node[above]{$\pmb{\tau}_1^{\mu_1}$} (3, 4.25);
\fill (1.5, 4.25) circle (0.045);
\fill (1.75, 4.25) circle (0.045);
\fill (2, 4.25) circle (0.045);
\fill (3, 4.25) circle (0.045);
\fill (1.5, 2.25) circle (0.045);
\fill (1.75, 2.25) circle (0.045);
\fill (2, 2.25) circle (0.045);
\fill (3, 2.25) circle (0.045);
\end{tikzpicture}
}\quad\raisebox{0.12\textwidth}{
\subfloat[]{
$
\begin{array}{c}
\includegraphics[width = 0.25\textwidth, page = 1]{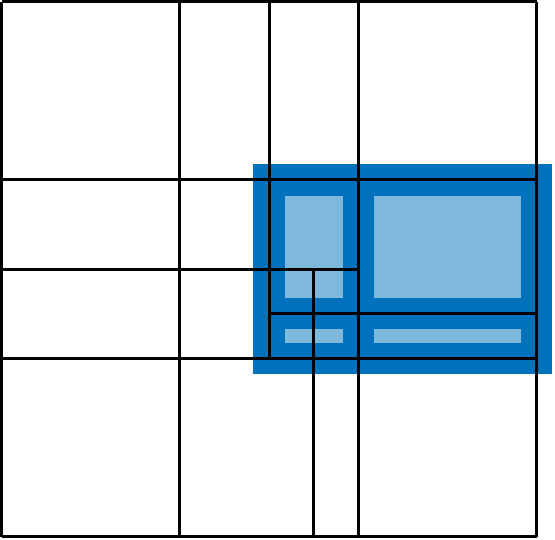}\\\\
\includegraphics[width = 0.25\textwidth, page = 2]{handinhand}
\end{array}
$
}}\quad
\subfloat[]{
\includegraphics[width = 0.333\textwidth, page = 3]{handinhand}
}
\caption{Example of an LR mesh that goes MS-wise but not LR-wise hand-in-hand.
The bidegree is $\m{p} = \pmb{1}$ and all meshlines have multiplicity $1$. The newly inserted segment $\gamma$ is highlighted in (a), together with the associated global knot vector $\pmb{\tau}_1^{\mu_1}$, which has $4$ breakpoints. From \eqref{eq:splinedim1d}, we deduce that $\dim \SS[\pmb{\tau}_1^{\mu_1}] = 2$.
The MS B-splines in $\cB[\gamma]$ are shown in (b). Only the upper one is an LR B-spline, so the cardinality of $\cL[\gamma]$ is $1$ and the mesh does not go LR-wise hand-in-hand. The corresponding univariate restrictions are depicted in (c). The two functions are clearly linearly independent, so the mesh goes MS-wise hand-in-hand after the insertion of $\gamma$.}\label{fig:handinhand}
\end{figure}

\begin{figure}[t!]
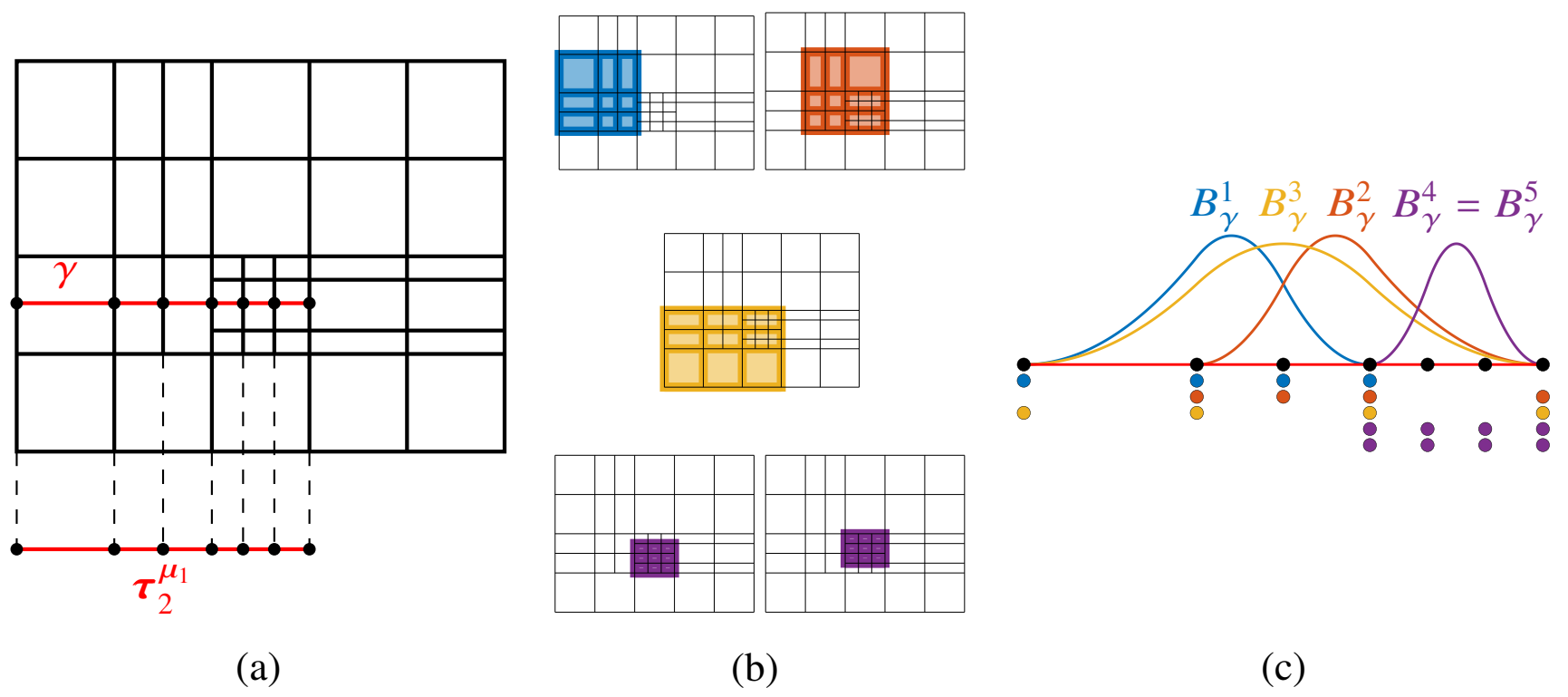

\centering
\subfloat[]{
\begin{tikzpicture}[scale = 2.85]
\draw[thick] (-.25,0) -- (1, 0) -- (1, 1) -- (-.25, 1) --cycle;
\draw[thick] (0,0) -- (0,1);
\draw[thick] (-.25,.25) -- (1, .25);
\draw[thick] (-.25, .5) -- (1, .5);
\draw[thick] (-.25, .75) -- (1, .75);
\draw[thick] (.25, 0) -- (.25, 1);
\draw[thick] (.5, 0) -- (.5, 1);
\draw[thick] (.75, 0) -- (.75, 1);
\draw[thick] (.125, .25) -- (.125, 1);
\draw[thick] (.25,.44) -- (1, .44);
\draw[thick] (.25, .31) -- (1, .31);
\draw[thick] (.33, .25) -- (.33, .5);
\draw[thick] (.41, .25) -- (.41, .5);
\draw[thick, red] (-.25,.38) -- (.5,.38);
\draw[red] (-0.125, 0.38)node[above]{$\gamma$};
\draw[dashed] (-0.25, 0) -- (-0.25, -0.25);
\draw[dashed] (0, 0) -- (0, -0.25);
\draw[dashed] (0.125, 0.25) -- (0.125, -0.25);
\draw[dashed] (0.25, 0) -- (0.25, -0.25);
\draw[dashed] (0.33, 0.25) -- (0.33, -0.25);
\draw[dashed] (0.41, 0.25) -- (0.41, -0.25);
\draw[dashed] (0.5, 0) -- (0.5, -0.25);
\draw[thick, red] (-.25, -.25) --node[below]{$\pmb{\tau}_2^{\mu_1}$} (0.5,-.25);
\fill (0,.38) circle (.016);
\fill (.125, .38) circle (.016);
\fill (.25,.38) circle (.016);
\fill (.33, .38) circle (.016);
\fill (.41,.38) circle (.016);
\fill (.5, .38) circle (.016);
\fill (-.25,.38) circle (.016);
\fill (0,-0.25) circle (.016);
\fill (.125, -0.25) circle (.016);
\fill (.25,-0.25) circle (.016);
\fill (.33, -0.25) circle (.016);
\fill (.41,-0.25) circle (.016);
\fill (.5, -0.25) circle (.016);
\fill (-.25,-0.25) circle (.016);
\end{tikzpicture}
}~\raisebox{0.185\textwidth}{
\subfloat[]{
$
\begin{array}{c}
\includegraphics[width = 0.125\textwidth, page = 4]{handinhand}~
\includegraphics[width = 0.125\textwidth, page = 5]{handinhand}\\\\
\includegraphics[width = 0.125\textwidth, page = 6]{handinhand}\\\\
\includegraphics[width = 0.125\textwidth, page = 7]{handinhand}~
\includegraphics[width = 0.125\textwidth, page = 8]{handinhand}
\end{array}
$}}~
\subfloat[]{\raisebox{0.1\textwidth}{
\begin{tikzpicture}
\node at (0, 0) {\includegraphics[width = 0.333\textwidth, page = 9]{handinhand}};
\node at (0, 1) {\textcolor[RGB]{236.895, 176.97, 31.875}{$B_\gamma^3$}};
\node at (-0.5, 1) {\textcolor[RGB]{0, 113.985, 188.955}{$B_\gamma^1$}};
\node at (0.5, 1) {\textcolor[RGB]{216.75, 82.8750, 24.99}{$B_\gamma^2$}};
\node at (1.35, 1) {\textcolor[RGB]{125.97, 46.92, 141.78}{$B_\gamma^4 = B_\gamma^5$}};
\end{tikzpicture}}
}
\caption{Example of an LR mesh that does not go hand-in-hand, neither LR-wise nor MS-wise.
The bidegree is $\m{p} = \pmb{2}$ and all meshlines have multiplicity $1$. The newly inserted segment $\gamma$ is highlighted in (a), together with the associated global knot vector $\pmb{\tau}_2^{\mu_1}$, which has $7$ breakpoints, and so the increase in dimension is $4$. The MS B-splines in $\cB[\gamma]$ are shown in (b). All of them are also LR B-splines. The corresponding $5$ univariate restrictions are depicted in (c). They are not linearly independent. Indeed, $B_\gamma^4 = B_\gamma^5$ and $B_\gamma^3$ can be expressed, using knot insertion, as a linear combination of $B_\gamma^1$ and $B_\gamma^2$. The $3$ linearly independent functions in $\{B_\gamma: B \in \cL[\gamma]\}$ are not sufficient to span a space of dimension $4$. Thus, the mesh does not go hand-in-hand, neither LR-wise nor MS-wise.}\label{fig:handinhand1}
\end{figure}

Theorem~\ref{thm:handinhand} reduces the verification of the hand-in-hand property to a univariate spanning condition. Since all functions $B_\gamma$ belong to $\SS[\pmb{\tau}_{p_{3-k}}^{\tilde{\mu}_{3-k}}]$, we always have
\begin{equation*}
\dim\spn  \{ B_\gamma : B \in  \cL[\gamma]\}\leq\dim \spn \{ B_\gamma : B \in \cB[\gamma] \}
\leq \dim \SS[\pmb{\tau}_{p_{3-k}}^{\tilde{\mu}_{3-k}}].
\end{equation*}
We distinguish two cases in which the inequality is strict:
\begin{itemize}
\item the cardinality of $\cL[\gamma]$ (resp., $\cB[\gamma]$) is smaller than $\dim \SS[\pmb{\tau}_{p_{3-k}}^{\tilde{\mu}_{3-k}}]$;
\item the cardinality of $\cL[\gamma]$ (resp., $\cB[\gamma]$) is at least $\dim \SS[\pmb{\tau}_{p_{3-k}}^{\tilde{\mu}_{3-k}}]$, but the functions $\{B_\gamma: B \in \cL[\gamma]\}$ (resp., $\{B_\gamma: B \in \cB[\gamma]\}$) are linearly dependent.
\end{itemize}
Figure~\ref{fig:handinhand} shows an example in which the mesh goes MS-wise hand-in-hand but not LR-wise hand-in-hand, due to a deficit in the cardinality of $\cL[\gamma]$.
Even if the cardinality of such set is sufficiently large, the univariate B-splines $\{B_\gamma\}$ may fail to span the entire spline space $\SS[\pmb{\tau}_{p_{3-k}}^{\tilde{\mu}_{3-k}}]$ due to linear dependence. An example is shown in Figure~\ref{fig:handinhand1}.

Heuristically, LR B-splines defined on LR meshes constructed through automatic refinement procedures typically satisfy the completeness property. It is conjectured that completeness holds whenever local linear independence is enforced.
To date, two mesh refinement strategies in the LR B-spline literature guarantee the completeness property: the \emph{hierarchical LR} \cite{hlr} and the \emph{effective grading} \cite{eg} procedures. Both approaches also ensure local linear independence of the resulting LR B-splines.

\subsection{Linear Dependence}
In this section, we present two geometric necessary conditions for linear dependence relations to occur.
These conditions impose strong structural constraints on the mesh.
We also discuss strategies for detecting and removing linear dependence relations. The most effective approach combines two techniques: the \emph{peeling algorithm} and the \emph{tensor expansion algorithm}.

Let $\cB_D \subseteq \cB$ be a subset of linearly dependent MS B-splines on a spline mesh $\cM^\mu$. We denote by $\cR_D \subset \RR^2$ the union of their supports.
We say that these functions are \textbf{actively linearly dependent} if they appear with non-zero coefficients in a linear dependence relation.

To state the first necessary condition for the existence of a linear dependence relation, we introduce the notion of nested B-splines.

\begin{definition}\label{def:nested}
Let $\cM^\mu$ be a given LR mesh, and let $B_1$ and $B_2$ be two distinct MS B-splines defined on $\cM^\mu$. We say that $B_2$ is \textbf{nested} in $B_1$, and write $B_2 \preceq B_1$, if
\begin{itemize}
\item $\supp B_2 \subseteq \supp B_1$, and
\item for every meshline $\gamma \in \cM^\mu$ such that $\gamma \subset \partial \supp B_1 \cap \partial \supp B_2$, it holds that
\begin{equation*}
\mu_{B_2}(\gamma) \leq \mu_{B_1}(\gamma).
\end{equation*}
\end{itemize}
\end{definition}

\begin{figure}[t!]
\centering
\subfloat[]{
\begin{tikzpicture}[scale=2.5]
\draw (0,0) -- (1,0) -- (1,1) -- (0,1) -- cycle;
\draw (.02,0) -- (.02,1);
\draw (.33, 0) -- (.33, 1);
\draw (.66, 0) -- (.66, 1);
\draw (0, .3) -- (1, .3);
\draw (0,.7) -- (1, .7);
\draw (.5, 0) -- (.5,.7);
\draw (0,.5) -- (.66,.5);
\end{tikzpicture}
}\quad
\subfloat[]{
\begin{tikzpicture}[scale=2.5]
\filldraw[fill=red!50, draw = red, ultra thick] (.02,0) -- (1,0) -- (1,1) -- (.02,1) -- cycle;
\draw[red, ultra thick] (.02,.33) -- (1,.33);
\draw[red, ultra thick] (.02,.66) -- (1,.66);
\draw[red, ultra thick] (.33,0) -- (.33,1);
\draw[red, ultra thick] (.66,0) -- (.66,1);
\draw (0,0) -- (1,0) -- (1,1) -- (0,1) -- cycle;
\draw (.02,0) -- (.02,1);
\draw (.33, 0) -- (.33, 1);
\draw (.66, 0) -- (.66, 1);
\draw (0, .33) -- (1, .33);
\draw (0,.66) -- (1, .66);
\draw (.5, 0) -- (.5,.66);
\draw (0,.5) -- (.66,.5);
\end{tikzpicture}
}\quad
\subfloat[]{
\begin{tikzpicture}[scale=2.5]
\filldraw[fill=green!50!black!50!white, draw = green!50!black, ultra thick] (.02,0) -- (.66,0) -- (.66,.66) -- (.02,.66) -- cycle;
\draw[green!50!black, ultra thick] (.02,.33) -- (.66,.33);
\draw[green!50!black, ultra thick] (.02,.5) -- (.66,.5);
\draw[green!50!black, ultra thick] (.5,0) -- (.5,.66);
\draw[green!50!black, ultra thick] (.33,0) -- (.33,.66);
\draw (0,0) -- (1,0) -- (1,1) -- (0,1) -- cycle;
\draw (.02,0) -- (.02,1);
\draw (.33, 0) -- (.33, 1);
\draw (.66, 0) -- (.66, 1);
\draw (0, .33) -- (1, .33);
\draw (0,.66) -- (1, .66);
\draw (.5, 0) -- (.5,.66);
\draw (0,.5) -- (.66,.5);
\end{tikzpicture}
}\quad
\subfloat[]{
\begin{tikzpicture}[scale=2.5]
\filldraw[fill=cyan!50, draw = cyan, ultra thick] (0,0) -- (.5,0) -- (.5,.66) -- (0,.66) -- cycle;
\draw[cyan, ultra thick] (0,.33) -- (.5,.33);
\draw[cyan, ultra thick] (0,.5) -- (.5,.5);
\draw[cyan, ultra thick] (.33,0) -- (.33,.66);
\draw[cyan, ultra thick] (.02,0) -- (.02,.66);
\draw (0,0) -- (1,0) -- (1,1) -- (0,1) -- cycle;
\draw (.02,0) -- (.02,1);
\draw (.33, 0) -- (.33, 1);
\draw (.66, 0) -- (.66, 1);
\draw (0, .33) -- (1, .33);
\draw (0,.66) -- (1, .66);
\draw (.5, 0) -- (.5,.66);
\draw (0,.5) -- (.66,.5);
\end{tikzpicture}
}\caption{Example~of nested LR B-splines on the LR mesh $\cM_{\m{p}}^\mu$ of bidegree $\m{p} = \pmb{2}$, shown in (a). All meshlines have multiplicity $1$, except those on the left boundary of $\cM_{\m{p}}^\mu$, highlighted with double lines, which have multiplicity $2$. Three biquadratic LR B-splines $B_1, B_2, B_3$ are depicted in (b)--(d) through their supports and associated local tensor meshes. All knots have multiplicity $1$, except for the first one in the local knot vector $\m{x}_{B_3}$, which has multiplicity $2$. Therefore, we have $B_2 \preceq B_1$, while $B_3 \npreceq B_2$ and $B_3 \npreceq B_1$. Note that $\supp B_3 \subseteq \supp B_2$ and $\supp B_3 \subseteq \supp B_1$.}
\label{fig:nestedex}
\end{figure}
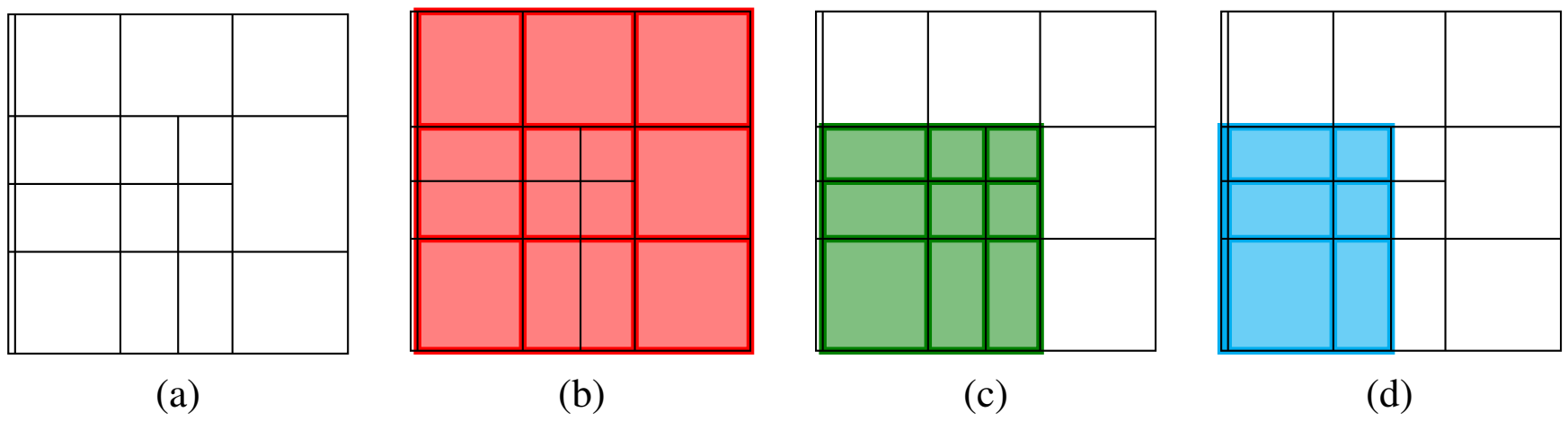

\begin{figure}[t!]
\centering
\subfloat[]{
\begin{tikzpicture}[scale=3.5]
\filldraw[fill=red!50, draw=red, ultra thick] (.25,.2) -- (1, .2) -- (1, .3) -- (1.25,.3) -- (1.25, .7) -- (1.15, .7) -- (1.15,.8) -- (.5, .8) -- (.5, .7) -- (.4, .7) -- (.4, .8) -- (.25,.8) -- cycle;
\draw (.75,.5) node{$\cR_D$};
\draw (.25,.2) -- (.55, .2) -- (.55,.4) -- (.25,.4) -- cycle;
\fill[pattern=north west lines, pattern color=black] (.25,.2) -- (.55, .2) -- (.55,.4) -- (.25,.4) -- cycle;
\draw (.25,.2) -- (.65, .2) -- (.65,.5) -- (.25,.5) -- cycle;
\fill[pattern=north east lines, pattern color=black] (.25,.2) -- (.65, .2) -- (.65,.5) -- (.25,.5) -- cycle;
\end{tikzpicture}}\qquad
\subfloat[]{
\begin{tikzpicture}[scale=3.5]
\filldraw[fill=red!50, draw=red, ultra thick] (.25,.2) -- (1, .2) -- (1, .3) -- (1.25,.3) -- (1.25, .7) -- (1.15, .7) -- (1.15,.8) -- (.5, .8) -- (.5, .7) -- (.4, .7) -- (.4, .8) -- (.25,.8) -- cycle;
\draw (.75,.5) node{$\cR_D$};
\draw (.35,.3) -- (.55,.3);
\draw (.45,.3) -- (.45,.5);

\draw (.85,.6) -- (1.05,.6);
\draw (.95,.6) -- (.95,.7);

\draw (.8,.25) -- (.8,.35);
\draw (.65,.3) -- (.8,.3);

\draw (.3,.5) -- (.3,.7);
\draw (.3,.6) -- (.5,.6);

\draw (.7,.7) -- (.9,.7);
\draw (.8,.7) -- (.8,.55);
\end{tikzpicture}}
\caption{Visualization of the necessary conditions for linear dependence in the LR B-spline set. Let $\cR_D$ be the union of the supports of the actively linearly dependent LR B-splines. Then, (a) at every convex corner of $\cR_D$ at least one B-spline is nested into another (Proposition~\ref{prop:nestednesscond}), and (b) there are at least five relevant T-vertices in the interior of $\cR_D$, with at least one for each orientation (Proposition~\ref{prop:fiveverticescond}).}
\label{fig:lindepconditions}
\end{figure}
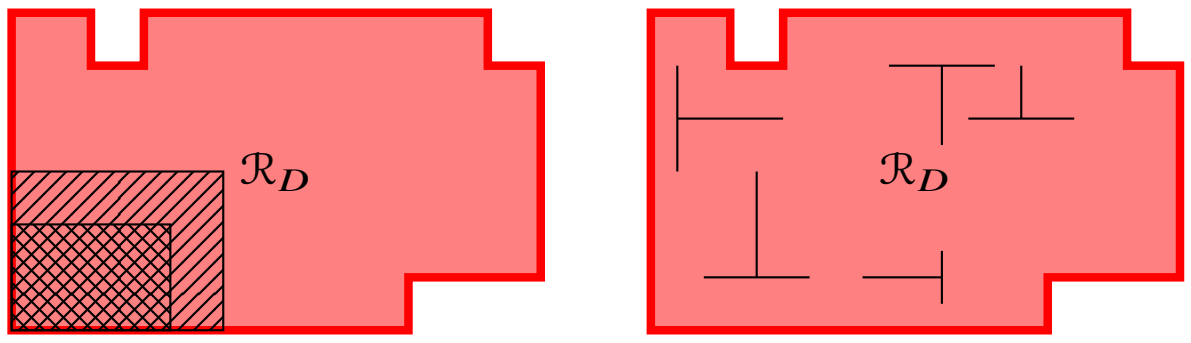

The notion of nested MS B-splines was originally introduced in a different form \cite{someproperties}. Definition~\ref{def:nested} is equivalent to that formulation (see \cite[Appendix]{n2s2}) but is simpler to describe.
Figure~\ref{fig:nestedex} illustrates an example of nested LR B-splines.
The first necessary condition for linear dependence is the following \textbf{nestedness condition}; see Figure~\ref{fig:lindepconditions}(a) and \cite[Proposition~7.3]{lindep}.

\begin{proposition}\label{prop:nestednesscond}
Let $\m{v}$ be any convex corner of $\cR_D$. Then, there exist $B_1, B_2 \in \cB_D$ such that
\begin{itemize}
\item $\m{v} \in \supp B_1 \cap \supp B_2$,
\item $\mu_{B_1}(\gamma) = \mu_{B_2}(\gamma)$ for all $\gamma \in \cM_{B_1} \cap \cM_{B_2}$,
\item $B_2 \preceq B_1$.
\end{itemize}
\end{proposition}

The nestedness condition in Proposition~\ref{prop:nestednesscond} cannot be satisfied on LR meshes of bidegree $\m{p} = \pmb{0}$. Indeed, nesting one B-spline into another during the construction of an LR mesh would require terminating a meshline in the interior of an element, which is not allowed. Therefore, both MS B-splines and LR B-splines of bidegree $\m{p} = \pmb{0}$ are linearly independent on any LR mesh.

To state the second condition for linear dependence, we introduce the notions of relevant T-vertices and their orientations.

\begin{definition}
A vertex $\m{v} \in \cR_D$ is called \textbf{relevant} for $\cB_D$ if it corresponds to a pair of knots in at least one MS B-spline in $\cB_D$.
A vertex $\m{v} \in \cR_D$ is a \textbf{relevant T-vertex} for $\cB_D$  if it is the intersection of three minimal meshlines belonging to local tensor meshes $\cM_B$, with $B\in\cB_D$, such that two of them are  collinear and the third minimal meshline, say $\gamma$, is orthogonal to the other two. We assign an orientation to $\m{v}$ according to the position of $\gamma$:
\begin{itemize}
\item \textbf{downward} if $\gamma$ lies below $\m{v}$,
\item \textbf{upward} if $\gamma$ lies above $\m{v}$,
\item \textbf{rightward} if $\gamma$ lies to the right of $\m{v}$,
\item \textbf{leftward} if $\gamma$ lies to the left of $\m{v}$.
\end{itemize}
\end{definition}

We can now state the second necessary condition for linear dependence,  the so-called \textbf{five T-vertices condition}; see Figure~\ref{fig:lindepconditions}(b) and \cite[Theorem~7.12]{lindep}.

\begin{proposition}\label{prop:fiveverticescond}
For a linear dependence relation among LR B-splines (resp., MS B-splines) of a set $\cB_D$, there exist at least five (resp., four) relevant T-vertices in $\mathring{\cR_D}$, with at least one vertex for each orientation. Moreover, each relevant T-vertex corresponds to a pair of knots shared by at least two LR B-splines (resp., MS B-splines) in $\cB_D$.
\end{proposition}

The results in Propositions~\ref{prop:nestednesscond} and~\ref{prop:fiveverticescond} are the main ingredients to prove the following lower bound theorem; see \cite[Propositions~7.15 and~8.2]{lindep}.

\begin{theorem}\label{thm:lowerbounds}
Let $\cM^\mu_{\m{p}}$ be a given LR mesh of bidegree $\m{p} \defeq (p_1, p_2)$ and let $\cB_D$ be a non-empty set of actively linearly dependent MS B-splines. It holds
\begin{itemize}
\item if $p_1 + p_2 \geq 2$, then
$| \cB_D |\geq 6$;
\item if $p_k \geq 2$ for some  $k \in \{1,2\}$, and  $\cB_D\subset\cL[\cM^\mu_{\m{p}}]$, then
$| \cB_D |\geq 8$.
\end{itemize}
Both bounds are sharp.
\end{theorem}

The configuration in Figure~\ref{fig:exmsnotlr} shows an example of linear dependence in the MS B-spline set for bidegree $\m{p} = \pmb{2}$. More examples are detailed below; see also \cite{lindep}.

\begin{example}{Example: Linear Dependence of MS B-Splines}%
Consider the LR meshes shown in Figure~\ref{fig:lindepms} for different bidegrees.
Each configuration is obtained by starting from an initial tensor spline mesh supporting $3$ tensor-product B-splines. By inserting local meshlines (first horizontal, then vertical), the spline space dimension increases to $9$, according to Theorem~\ref{thm:splinedimincrease2d}. The MS B-spline set contains $10$ functions, of which $6$ are actively linearly dependent, highlighted in red and blue in the figure.

\begin{figure}[t!]
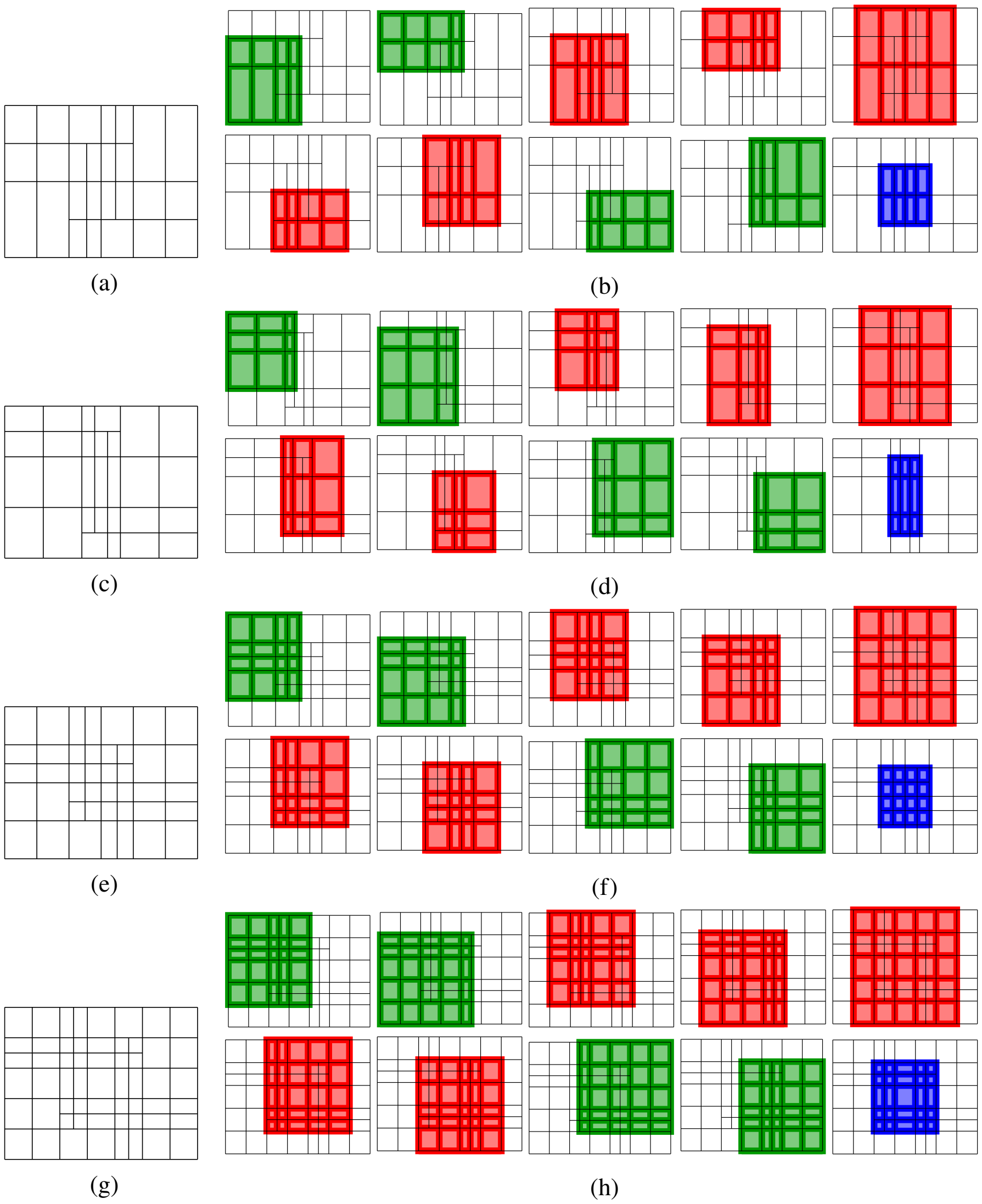

\centering
\newcounter{p}
\setcounter{p}{1}

\subfloat[]{%
\includegraphics[width=0.2\textwidth,page=\thep]{lindep}%
\stepcounter{p}
}
~\raisebox{0.1225\textwidth}{
\subfloat[]{
$\begin{array}{c}
\includegraphics[width=0.15\textwidth,page=\thep]{lindep}\stepcounter{p}~
\includegraphics[width=0.15\textwidth,page=\thep]{lindep}\stepcounter{p}~
\includegraphics[width=0.15\textwidth,page=\thep]{lindep}\stepcounter{p}~
\includegraphics[width=0.15\textwidth,page=\thep]{lindep}\stepcounter{p}~
\includegraphics[width=0.15\textwidth,page=\thep]{lindep}\stepcounter{p}~
\\
\includegraphics[width=0.15\textwidth,page=\thep]{lindep}\stepcounter{p}~
\includegraphics[width=0.15\textwidth,page=\thep]{lindep}\stepcounter{p}~
\includegraphics[width=0.15\textwidth,page=\thep]{lindep}\stepcounter{p}~
\includegraphics[width=0.15\textwidth,page=\thep]{lindep}\stepcounter{p}~
\includegraphics[width=0.15\textwidth,page=\thep]{lindep}\stepcounter{p}~
\end{array}$
}}\\
\subfloat[]{%
\includegraphics[width=0.2\textwidth,page=\thep]{lindep}%
\stepcounter{p}
}
~\raisebox{0.1225\textwidth}{
\subfloat[]{
$\begin{array}{c}
\includegraphics[width=0.15\textwidth,page=\thep]{lindep}\stepcounter{p}~
\includegraphics[width=0.15\textwidth,page=\thep]{lindep}\stepcounter{p}~
\includegraphics[width=0.15\textwidth,page=\thep]{lindep}\stepcounter{p}~
\includegraphics[width=0.15\textwidth,page=\thep]{lindep}\stepcounter{p}~
\includegraphics[width=0.15\textwidth,page=\thep]{lindep}\stepcounter{p}~
\\
\includegraphics[width=0.15\textwidth,page=\thep]{lindep}\stepcounter{p}~
\includegraphics[width=0.15\textwidth,page=\thep]{lindep}\stepcounter{p}~
\includegraphics[width=0.15\textwidth,page=\thep]{lindep}\stepcounter{p}~
\includegraphics[width=0.15\textwidth,page=\thep]{lindep}\stepcounter{p}~
\includegraphics[width=0.15\textwidth,page=\thep]{lindep}\stepcounter{p}~
\end{array}$
}}\\
\subfloat[]{%
\includegraphics[width=0.2\textwidth,page=\thep]{lindep}%
\stepcounter{p}
}
~\raisebox{0.1225\textwidth}{
\subfloat[]{
$\begin{array}{c}
\includegraphics[width=0.15\textwidth,page=\thep]{lindep}\stepcounter{p}~
\includegraphics[width=0.15\textwidth,page=\thep]{lindep}\stepcounter{p}~
\includegraphics[width=0.15\textwidth,page=\thep]{lindep}\stepcounter{p}~
\includegraphics[width=0.15\textwidth,page=\thep]{lindep}\stepcounter{p}~
\includegraphics[width=0.15\textwidth,page=\thep]{lindep}\stepcounter{p}~
\\
\includegraphics[width=0.15\textwidth,page=\thep]{lindep}\stepcounter{p}~
\includegraphics[width=0.15\textwidth,page=\thep]{lindep}\stepcounter{p}~
\includegraphics[width=0.15\textwidth,page=\thep]{lindep}\stepcounter{p}~
\includegraphics[width=0.15\textwidth,page=\thep]{lindep}\stepcounter{p}~
\includegraphics[width=0.15\textwidth,page=\thep]{lindep}\stepcounter{p}~
\end{array}$
}}\\
\subfloat[]{%
\includegraphics[width=0.2\textwidth,page=\thep]{lindep}%
\stepcounter{p}
}
~\raisebox{0.1225\textwidth}{
\subfloat[]{
$\begin{array}{c}
\includegraphics[width=0.15\textwidth,page=\thep]{lindep}\stepcounter{p}~
\includegraphics[width=0.15\textwidth,page=\thep]{lindep}\stepcounter{p}~
\includegraphics[width=0.15\textwidth,page=\thep]{lindep}\stepcounter{p}~
\includegraphics[width=0.15\textwidth,page=\thep]{lindep}\stepcounter{p}~
\includegraphics[width=0.15\textwidth,page=\thep]{lindep}\stepcounter{p}~
\\
\includegraphics[width=0.15\textwidth,page=\thep]{lindep}\stepcounter{p}~
\includegraphics[width=0.15\textwidth,page=\thep]{lindep}\stepcounter{p}~
\includegraphics[width=0.15\textwidth,page=\thep]{lindep}\stepcounter{p}~
\includegraphics[width=0.15\textwidth,page=\thep]{lindep}\stepcounter{p}~
\includegraphics[width=0.15\textwidth,page=\thep]{lindep}\stepcounter{p}~
\end{array}$
}}
\caption{LR meshes exhibiting linear dependence in the MS B-spline set.
The meshes have constant multiplicity equal to $1$.
The bidegrees are $\m{p} = (3,1)$ in (a)--(b), $\m{p} = \pmb{2}$ in  (c)--(d), $\m{p} = \pmb{3}$ in (e)--(f), and $\m{p} = \pmb{4}$ in (g)--(h).
For each configuration, the B-splines in red and blue are linearly dependent. The red functions are LR B-splines, while the blue ones are MS B-splines but not LR B-splines. The green functions are LR B-splines that are not actively involved in the linear dependence relation.}
\label{fig:lindepms}
\end{figure}

A similar configuration can be achieved for bidegree $\m{p} = \pmb{1}$. In this case, however, not only local insertions but also two additional meshline extensions are required. Moreover, in order to build a proper LR mesh during these extensions, we need a finer initial tensor spline mesh supporting $9$ tensor-product B-splines.
The resulting configuration is shown in Figure~\ref{fig:lindepms1}. After the local insertions and extensions, the spline space dimension increases to $17$, while the number of MS B-splines becomes $18$. Exactly $6$ MS B-splines are actively linearly dependent.

\begin{figure}[t!]
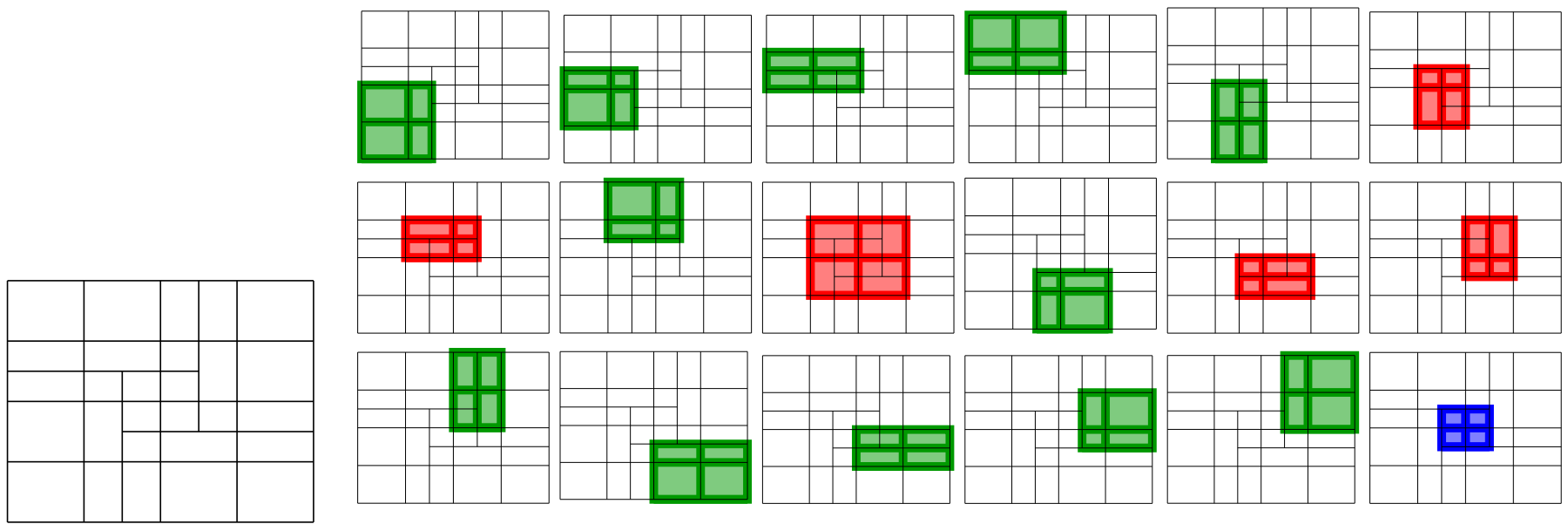

\centering
\setcounter{p}{1}

\subfloat[]{%
\includegraphics[width=0.2\textwidth,page=\thep]{lindep1}%
\stepcounter{p}
}
~\raisebox{0.1625\textwidth}{
\subfloat[]{
$\begin{array}{c}
\includegraphics[width=0.125\textwidth,page=\thep]{lindep1}\stepcounter{p}~
\includegraphics[width=0.125\textwidth,page=\thep]{lindep1}\stepcounter{p}~
\includegraphics[width=0.125\textwidth,page=\thep]{lindep1}\stepcounter{p}~
\includegraphics[width=0.125\textwidth,page=\thep]{lindep1}\stepcounter{p}~
\includegraphics[width=0.125\textwidth,page=\thep]{lindep1}\stepcounter{p}~
\includegraphics[width=0.125\textwidth,page=\thep]{lindep1}\stepcounter{p}~
\\
\includegraphics[width=0.125\textwidth,page=\thep]{lindep1}\stepcounter{p}~
\includegraphics[width=0.125\textwidth,page=\thep]{lindep1}\stepcounter{p}~
\includegraphics[width=0.125\textwidth,page=\thep]{lindep1}\stepcounter{p}~
\includegraphics[width=0.125\textwidth,page=\thep]{lindep1}\stepcounter{p}~
\includegraphics[width=0.125\textwidth,page=\thep]{lindep1}\stepcounter{p}~
\includegraphics[width=0.125\textwidth,page=\thep]{lindep1}\stepcounter{p}~
\\
\includegraphics[width=0.125\textwidth,page=\thep]{lindep1}\stepcounter{p}~
\includegraphics[width=0.125\textwidth,page=\thep]{lindep1}\stepcounter{p}~
\includegraphics[width=0.125\textwidth,page=\thep]{lindep1}\stepcounter{p}~
\includegraphics[width=0.125\textwidth,page=\thep]{lindep1}\stepcounter{p}~
\includegraphics[width=0.125\textwidth,page=\thep]{lindep1}\stepcounter{p}~
\includegraphics[width=0.125\textwidth,page=\thep]{lindep1}\stepcounter{p}~
\end{array}$
}}
\caption{An LR mesh exhibiting linear dependence in the MS B-spline set of bidegree $\m{p} = \pmb{1}$. The B-splines in red and blue are linearly dependent. The red functions are LR B-splines, while the blue one is an MS B-spline but not an LR B-spline. The green functions are LR B-splines that are not actively involved in the linear dependence relation.}
\label{fig:lindepms1}
\end{figure}
\end{example}
\vspace*{-0.3cm} 
\begin{example}{Example: Linear Dependence of LR B-Splines}%
Consider the LR meshes shown in Figure~\ref{fig:lindeplr} for different bidegrees. These meshes are obtained by refining those in Figure~\ref{fig:lindepms} through the insertion of an additional meshline. This refinement increases the spline space dimension by $2$ but the number of LR B-splines by $3$. Consequently, the spline space has dimension $11$, while the number of LR B-splines is $12$. Exactly $8$ LR B-splines are actively linearly dependent.

\begin{figure}[t!]
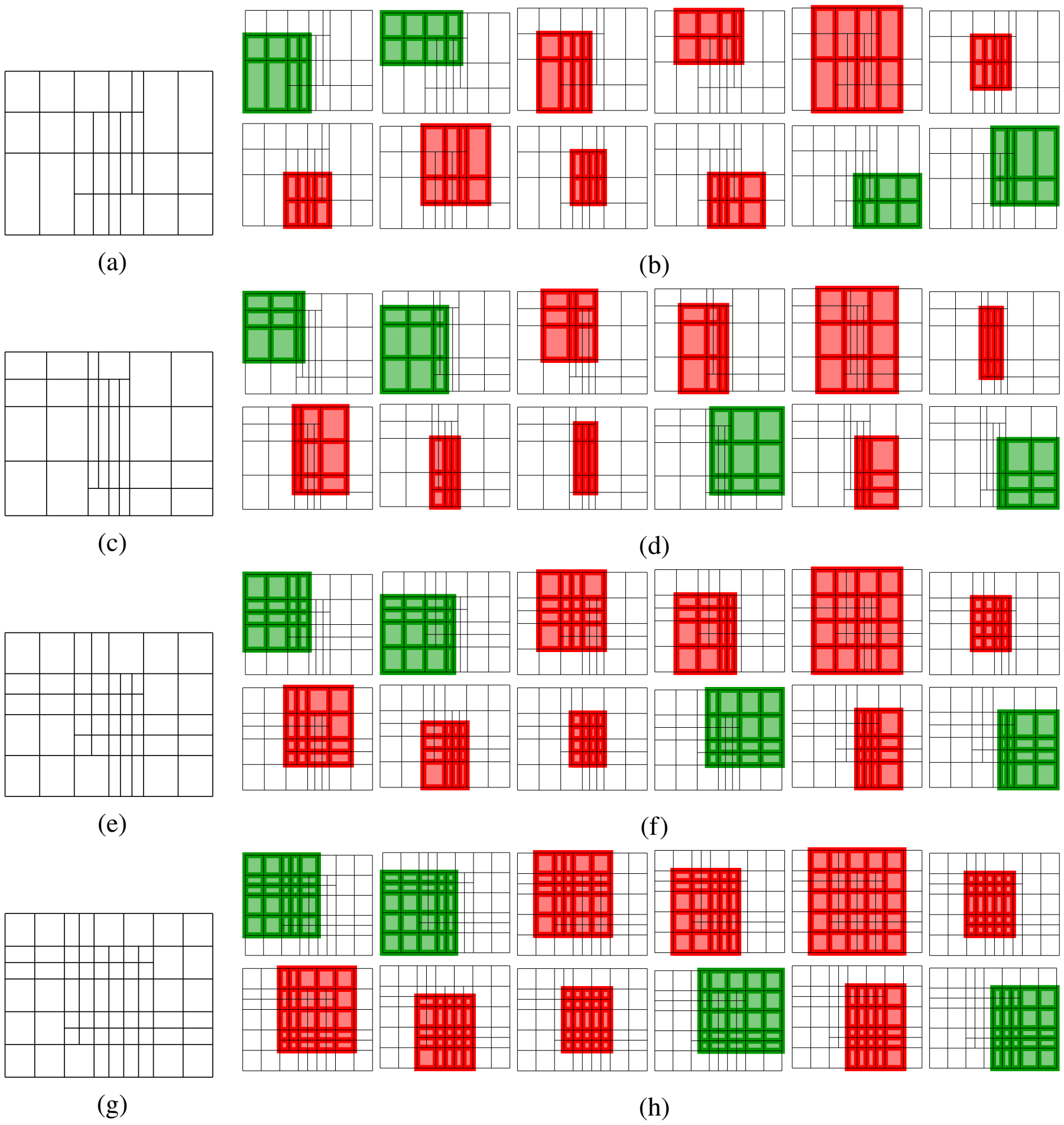

\centering
\setcounter{p}{14}

\subfloat[]{%
\includegraphics[width=0.2\textwidth,page=\thep]{lindep2}%
\stepcounter{p}
}
~\raisebox{0.1025\textwidth}{
\subfloat[]{
$\begin{array}{c}
\includegraphics[width=0.125\textwidth,page=\thep]{lindep2}\stepcounter{p}~
\includegraphics[width=0.125\textwidth,page=\thep]{lindep2}\stepcounter{p}~
\includegraphics[width=0.125\textwidth,page=\thep]{lindep2}\stepcounter{p}~
\includegraphics[width=0.125\textwidth,page=\thep]{lindep2}\stepcounter{p}~
\includegraphics[width=0.125\textwidth,page=\thep]{lindep2}\stepcounter{p}~
\includegraphics[width=0.125\textwidth,page=\thep]{lindep2}\stepcounter{p}~
\\
\includegraphics[width=0.125\textwidth,page=\thep]{lindep2}\stepcounter{p}~
\includegraphics[width=0.125\textwidth,page=\thep]{lindep2}\stepcounter{p}~
\includegraphics[width=0.125\textwidth,page=\thep]{lindep2}\stepcounter{p}~
\includegraphics[width=0.125\textwidth,page=\thep]{lindep2}\stepcounter{p}~
\includegraphics[width=0.125\textwidth,page=\thep]{lindep2}\stepcounter{p}~
\includegraphics[width=0.125\textwidth,page=\thep]{lindep2}\stepcounter{p}~
\end{array}$
}}\\
\subfloat[]{%
\includegraphics[width=0.2\textwidth,page=\thep]{lindep2}%
\stepcounter{p}
}
~\raisebox{0.1025\textwidth}{
\subfloat[]{
$\begin{array}{c}
\includegraphics[width=0.125\textwidth,page=\thep]{lindep2}\stepcounter{p}~
\includegraphics[width=0.125\textwidth,page=\thep]{lindep2}\stepcounter{p}~
\includegraphics[width=0.125\textwidth,page=\thep]{lindep2}\stepcounter{p}~
\includegraphics[width=0.125\textwidth,page=\thep]{lindep2}\stepcounter{p}~
\includegraphics[width=0.125\textwidth,page=\thep]{lindep2}\stepcounter{p}~
\includegraphics[width=0.125\textwidth,page=\thep]{lindep2}\stepcounter{p}~
\\
\includegraphics[width=0.125\textwidth,page=\thep]{lindep2}\stepcounter{p}~
\includegraphics[width=0.125\textwidth,page=\thep]{lindep2}\stepcounter{p}~
\includegraphics[width=0.125\textwidth,page=\thep]{lindep2}\stepcounter{p}~
\includegraphics[width=0.125\textwidth,page=\thep]{lindep2}\stepcounter{p}~
\includegraphics[width=0.125\textwidth,page=\thep]{lindep2}\stepcounter{p}~
\includegraphics[width=0.125\textwidth,page=\thep]{lindep2}\stepcounter{p}~
\end{array}$
}}\\
\subfloat[]{%
\includegraphics[width=0.2\textwidth,page=\thep]{lindep2}%
\stepcounter{p}
}
~\raisebox{0.1025\textwidth}{
\subfloat[]{
$\begin{array}{c}
\includegraphics[width=0.125\textwidth,page=\thep]{lindep2}\stepcounter{p}~
\includegraphics[width=0.125\textwidth,page=\thep]{lindep2}\stepcounter{p}~
\includegraphics[width=0.125\textwidth,page=\thep]{lindep2}\stepcounter{p}~
\includegraphics[width=0.125\textwidth,page=\thep]{lindep2}\stepcounter{p}~
\includegraphics[width=0.125\textwidth,page=\thep]{lindep2}\stepcounter{p}~
\includegraphics[width=0.125\textwidth,page=\thep]{lindep2}\stepcounter{p}~
\\
\includegraphics[width=0.125\textwidth,page=\thep]{lindep2}\stepcounter{p}~
\includegraphics[width=0.125\textwidth,page=\thep]{lindep2}\stepcounter{p}~
\includegraphics[width=0.125\textwidth,page=\thep]{lindep2}\stepcounter{p}~
\includegraphics[width=0.125\textwidth,page=\thep]{lindep2}\stepcounter{p}~
\includegraphics[width=0.125\textwidth,page=\thep]{lindep2}\stepcounter{p}~
\includegraphics[width=0.125\textwidth,page=\thep]{lindep2}\stepcounter{p}~
\end{array}$
}}\\
\subfloat[]{%
\includegraphics[width=0.2\textwidth,page=\thep]{lindep2}%
\stepcounter{p}
}
~\raisebox{0.1025\textwidth}{
\subfloat[]{
$\begin{array}{c}
\includegraphics[width=0.125\textwidth,page=\thep]{lindep2}\stepcounter{p}~
\includegraphics[width=0.125\textwidth,page=\thep]{lindep2}\stepcounter{p}~
\includegraphics[width=0.125\textwidth,page=\thep]{lindep2}\stepcounter{p}~
\includegraphics[width=0.125\textwidth,page=\thep]{lindep2}\stepcounter{p}~
\includegraphics[width=0.125\textwidth,page=\thep]{lindep2}\stepcounter{p}~
\includegraphics[width=0.125\textwidth,page=\thep]{lindep2}\stepcounter{p}~
\\
\includegraphics[width=0.125\textwidth,page=\thep]{lindep2}\stepcounter{p}~
\includegraphics[width=0.125\textwidth,page=\thep]{lindep2}\stepcounter{p}~
\includegraphics[width=0.125\textwidth,page=\thep]{lindep2}\stepcounter{p}~
\includegraphics[width=0.125\textwidth,page=\thep]{lindep2}\stepcounter{p}~
\includegraphics[width=0.125\textwidth,page=\thep]{lindep2}\stepcounter{p}~
\includegraphics[width=0.125\textwidth,page=\thep]{lindep2}\stepcounter{p}~
\end{array}$
}}
\caption{LR meshes exhibiting linear dependence in the LR B-spline set.
The meshes have constant multiplicity equal to $1$.
The bidegrees are 
$\m{p} = (3,1)$ in (a)--(b), $\m{p} = \pmb{2}$ in (c)--(d), $\m{p} = \pmb{3}$ (e)--(f), and $\m{p} = \pmb{4}$ in (g)--(h).
For each configuration, the LR B-splines in red are linearly dependent. The LR B-splines in green are not actively involved in the linear dependence relation.}
\label{fig:lindeplr}
\end{figure}
\end{example}

We now address the problem of detecting whether and where a linear dependence relation occurs on a given LR mesh, focusing on the LR B-spline set. To this end, we introduce the notions of overloaded elements and overloaded B-splines. We recall that, for any open LR mesh $\cM^\mu_{\m{p}}$, we have $\PP_{\m{p}}\subseteq \cL[\cM^\mu_{\m{p}}]$.

\begin{definition}
An element of the box partition $\mathcal{E}$ is an \textbf{overloaded element} if it is contained in the support of more B-splines than necessary to span the local polynomial space $\PP_{\m{p}}$, i.e., more than $(p_1+1)(p_2+1)$ supports.
A B-spline whose support consists only of overloaded elements is called an \textbf{overloaded B-spline}.
\end{definition}

\begin{figure}[t!]
\centering
\subfloat[]{
\includegraphics[width = 0.495\textwidth, page = 1]{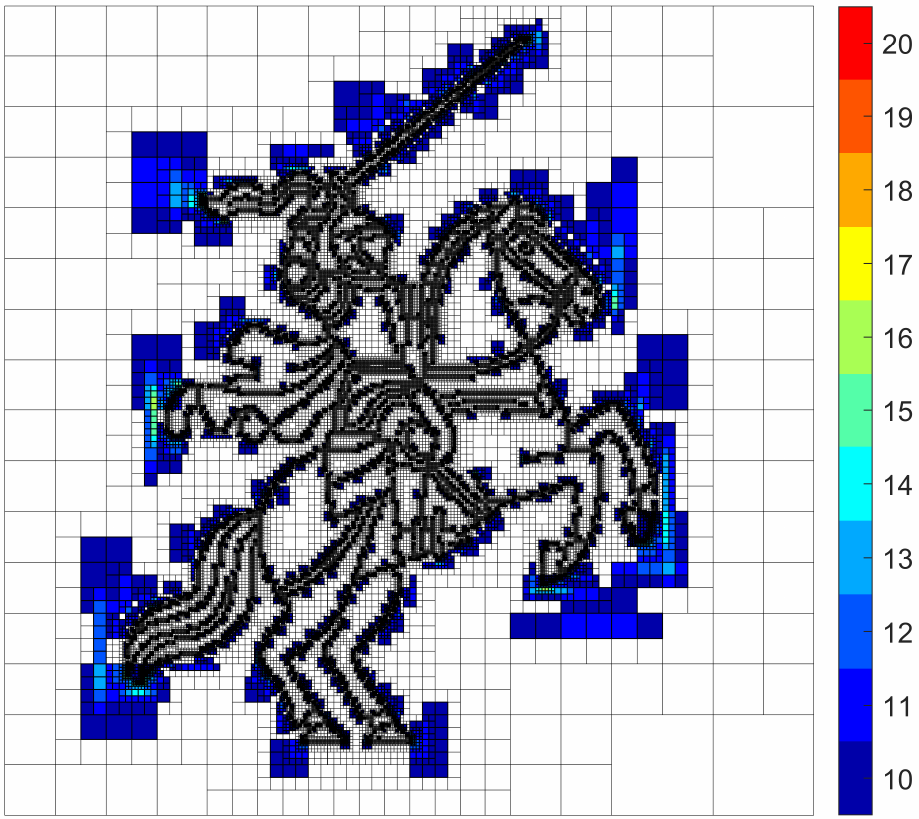}
}
\subfloat[]{
\includegraphics[width = 0.495\textwidth, page = 2]{peel}
}
\caption{Examples of overloaded elements in LR meshes of bidegree $\m{p} = \pmb{2}$. The highlighted cells are overloaded, i.e., they are contained in more than $(p_1+1)(p_2+1) = 9$ B-spline supports. The color scale indicates the number of supports containing each element.}
\label{fig:overknight}
\end{figure}

Figure~\ref{fig:overknight} shows examples of overloaded elements on LR meshes.

In a linear dependence relation, any redundant B-spline can be removed without affecting the spanning properties over the elements in its support. Consequently, only overloaded B-splines can participate in linear dependence relations \cite[Lemma~6.2]{lr}.
Therefore, a linear dependence relation must involve at least two overloaded B-splines on every element. Moreover,
if an overloaded B-spline exhibits reduced continuity across a meshline and participates in a linear dependence relation, then such reduced smoothness must also be shared by at least one other overloaded B-spline.
Therefore,  if on a given element (or along a meshline with a prescribed multiplicity) there is only one overloaded B-spline, then such a B-spline cannot be active in a linear dependence relation.

This observation forms the basis of the so-called \textbf{peeling algorithm}; see \cite[Algorithm~6.3]{lr}. Starting from the collection $\cB_O$ of overloaded B-splines, the algorithm iteratively removes (``peels away'') those B-splines that do not share any element or meshline (with the same multiplicity) with another overloaded B-spline.
As B-splines are removed, the peeling condition may then hold for other B-splines, which can then be removed in subsequent iterations. This process is repeated until no further removals are possible.
If the remaining collection contains fewer than $8$ functions, then no linear dependence relation can exist, thanks to Theorem~\ref{thm:lowerbounds}. A matrix description of the peeling algorithm can be found in \cite[Section~6]{lr}.

The peeling algorithm also localizes regions where linear dependence may occur. Distinct regions correspond to distinct subcollections of overloaded B-splines. If any of these subcollections contains fewer than $8$ elements, it cannot create a linear dependence relation.

\begin{figure}[t!]
\centering
\subfloat[]{
\includegraphics[width = 0.372\textwidth, page = 1]{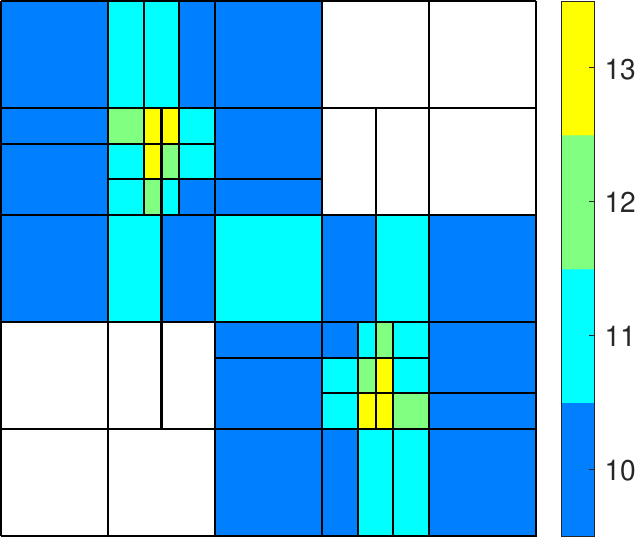}
}
\subfloat[]{
\includegraphics[width = 0.314\textwidth, page = 2]{texpansion}
}
\subfloat[]{
\includegraphics[width = 0.314\textwidth, page = 3]{texpansion}
}
\caption{Example of tensor expansion for bidegree $\m{p} = \pmb{2}$. The overloaded elements are colored in (a). The corresponding overloaded B-splines form a set of $16$ functions that cannot be eliminated by the peeling algorithm. The associated submesh $\cM_O$ is shown in (b). The resulting tensor mesh $\cM_T$ obtained by extending all meshlines in $\cM_O$ is shown in (c). The overloaded B-splines are then refined through knot insertions and expressed in the tensor-product basis on $\cM_T$. In this example, the functions are linearly dependent.}
\label{fig:tensorexpansion}
\end{figure}

The peeling algorithm may fail to conclude linear independence. In this case, it identifies a set of \emph{candidate} subcollections $\cB_O^1, \ldots, \cB_O^\ell$ of overloaded B-splines where no further removal is possible. For each of them, linear independence can be tested using the so-called \textbf{tensor expansion algorithm}; see Figure~\ref{fig:tensorexpansion}.

The latter algorithm consists of extending all meshlines of the local tensor meshes $\cM_B$, with $B\in\cB_O^i$, to obtain a local tensor mesh $\cM_T^i$. During this extension process, each B-spline in $\cB_O^i$ is possibly refined through repeated knot insertions and expressed in the tensor-product B-spline basis $\cB_T^i$ defined on $\cM_T^i$.
This yields a change-of-basis matrix whose rows contain the coefficients of the B-splines in $\cB_O^i$ expressed in the basis $\cB_T^i$.  From basic linear algebra, we obtain that the functions in $\cB_O^i$ are linearly independent if and only if such matrix  has full  rank.

One could in principle apply tensor expansion directly after identifying overloaded B-splines, without using the peeling algorithm. However, the number of overloaded B-splines may be large, while many of them can be quickly discarded by peeling.
Since tensor expansion is computationally expensive, due to the loss of locality when passing to a tensor-product representation, it should be applied only to the reduced candidate sets.

In case of linear dependence, local mesh modifications (e.g., extending meshlines to further refine the overloaded B-splines) eventually eliminate the dependence.
In the worst case, a full tensor expansion is performed, yielding a tensor-product mesh where no linear dependence occurs. In practice, however, small local refinements are often sufficient.
Refinement strategies ensuring (local) linear independence will be discussed in Section~\ref{sec:refinement}.

\subsection{Local Linear Independence}\label{sec:loclinind}
Linear dependence may (rarely) occur if no constraints are imposed on the LR mesh construction. However, it is easier to check local linear independence as there exists a very elegant characterization of the corresponding LR meshes.

Functions in $\RR^d$ are \textbf{locally linearly independent} if they are linearly independent on every open subset $A \subseteq \RR^d$ intersecting their support.
Local linear independence clearly implies (global) linear independence.

For spline functions on box partitions, it suffices to consider open sets $A$ given by single elements, because of the polynomial structure of splines and the following proposition. Note that the proposition is a direct consequence of the fact that any non-trivial univariate polynomial of degree $p$ has at most $p$ zeros on the real line.

\begin{proposition}
Let $\mathfrak p,\mathfrak q$ be univariate polynomials.
Suppose that there exists a set $A \subseteq \RR$ having an accumulation point such that
$\mathfrak p (x) = \mathfrak q(x)$ for all $x \in A$.
Then, $\mathfrak p \equiv \mathfrak q$.
\end{proposition}

The characterization of LR meshes ensuring local linear independence is closely related to the absence of nested B-splines; see Definition~\ref{def:nested}.

\begin{definition}\label{def:N2S}
An LR mesh is said to satisfy the \textbf{non-nested-support} (\textbf{N$_2$S}) property if no LR B-spline is nested into another.
\end{definition}

We are now ready to provide some equivalent statements to local linear independence for LR B-splines; see \cite[Theorem~5.1, Corollary 3.8]{someproperties} and \cite[Theorem~4]{hlr}.

\begin{theorem}\label{th:lli}
Let $\cM^\mu$ be a given LR mesh, and let $\cB$ and $\cL$ be the corresponding MS B-spline set and LR B-spline set, respectively. The following statements are equivalent:
\begin{itemize}
\item $\cM^\mu$ satisfies the N$_2$S property;
\item monomials admit the explicit representation  \eqref{eq:explicitmonomials} in terms of $\cL$;
\item the LR B-splines in $\cL$ form a partition of unity (i.e., all scaling weights are equal to $1$);
\item no element is overloaded, i.e., each element belongs to the support of exactly $(p_1+1)(p_2+1)$ LR B-splines;
\item the LR B-splines in $\cL$ are locally linearly independent.
\end{itemize}
Moreover, if $\cM^\mu$ satisfies the N$_2$S property, then $\cL = \cB$.
\end{theorem}

This result shows that LR B-splines defined on LR meshes satisfying the N$_2$S property retain all key properties of tensor-product B-splines, while still allowing local refinement.
The characterization in Theorem~\ref{th:lli} does not directly provide a constructive refinement strategy ensuring the N$_2$S property. More recent works have introduced refinement procedures that guarantee this property by construction \cite{hlr,rmb,eg,n2s2}. They will be briefly described in the next section.

\subsection{Automatic Refinement Strategies}\label{sec:refinement}
In this section, we present refinement strategies to construct LR B-splines enjoying local linear independence (and thus also linear independence). As shown in Theorem~\ref{th:lli}, it suffices to avoid or eliminate nested B-splines.

The most natural way to enforce the N$_2$S property is to extend the meshlines in the local tensor meshes of nested LR B-splines so that they traverse the support of the LR B-spline containing them. These extensions trigger the refinement of the containing LR B-spline into  LR B-splines with smaller supports, which eventually eliminate the nesting configuration.

This is the basic idea underlying the \textbf{non-nested-support structured} (\textbf{N$_2$S$_2$}) \textbf{mesh refinement} strategy \cite{n2s2}. The procedure consists of two steps.
\begin{enumerate}
\item A subset of LR B-splines marked for refinement, according to a given error indicator or design criterion, is refined by dyadically splitting the largest boxes in their local tensor meshes. This procedure is referred to as \textbf{structured mesh refinement} \cite{johannessen2014}. No guarantee of the N$_2$S property is provided, and nesting configurations are likely to appear.

\item In order to restore the N$_2$S property, the refined LR B-spline set is scanned. Whenever an LR B-spline contains nested LR B-splines, the meshlines in one of the two coordinate directions in the local tensor meshes of the nested functions are extended to traverse the support of the containing LR B-spline. This induces further local refinement of the mesh. The procedure iterates until no nesting can
 be found anymore.
\end{enumerate}

\begin{figure}[t!]
	\centering
	\subfloat[]{
		\begin{tikzpicture}[scale = 0.35]
			\fill[orange!90!black, opacity=.5] (4, 4) rectangle (7, 7);
			\draw[orange!90!black, line width=4pt] (4, 4) -- (4, 7) -- (7, 7) -- (7, 4) -- cycle;
			\draw[orange!90!black, step = 1, line width=4pt] (4, 4) grid (7, 7);
			\draw[step = 1] (1, 2) grid (10, 9);
	\end{tikzpicture}}\quad
	\subfloat[]{
		\begin{tikzpicture}[scale = 0.35]
			\fill[orange!90!black, opacity=.5] (4, 4) rectangle (7, 7);
			\draw[orange!90!black, line width=4pt] (4, 4) -- (4, 7) -- (7, 7) -- (7, 4) -- cycle;
			\draw[orange!90!black, step = 1, line width=4pt] (4, 4) grid (7, 7);
			\draw[step = 1] (1, 2) grid (10, 9);
			\draw[step = 0.5] (4, 4) grid (7, 7);
	\end{tikzpicture}}\quad
	\subfloat[]{
		\begin{tikzpicture}[scale = 0.35] 
			\fill[orange!90!black, opacity=.5] (4, 4) rectangle (7, 7);
			\draw[orange!90!black, line width=4pt] (4, 4) -- (4, 7) -- (7, 7) -- (7, 4) -- cycle;
			\draw[orange!90!black, step = 1, line width=4pt] (4, 4) grid (7, 7);
			\draw[step = 1] (1, 2) grid (10, 9);
			\draw[step = 0.5] (4, 4) grid (7, 7);
			\foreach \y in {0,...,2}{
				\draw (2, {4.5 + \y}) -- (4, {4.5 + \y});
				\draw (7, {4.5 + \y}) -- (9, {4.5 + \y});
			};
	\end{tikzpicture}}\\
	\subfloat[]{
		\begin{tabular}{c}
			\begin{tikzpicture}[scale=1.25]
				\filldraw[draw=red, ultra thick, fill=red!50] (0,.495) -- (.495,.495) -- (.495,1) -- (0,1) -- cycle;
				\draw[red, ultra thick] (.165,.495) -- (.165,1);
				\draw[red, ultra thick] (.33,.495) -- (.33,1);
				\draw[red, ultra thick] (0,.66) -- (.495,.66);
				\draw[red, ultra thick] (0,.825) -- (.495,.825);
				\draw (0,0) -- (1,0) -- (1,1) -- (0,1) -- cycle;
				\draw (.33,0) -- (.33, 1);
				\draw (.66,0) -- (.66,1);
				\draw (0,.66) -- (1,.66);
				\draw (0,.33) -- (1,.33);
				\draw (0,.495) -- (.66,.495);
				\draw (0,.825) -- (.66,.825);
				\draw (.495,.33) -- (.495,1);
				\draw (.165,.33) -- (.165,1);
			\end{tikzpicture},
			\begin{tikzpicture}[scale=1.25]
				\filldraw[draw=red, ultra thick, fill=red!50] (.165,.495) -- (.66,.495) -- (.66,1) -- (.165,1) -- cycle;
				\draw[red, ultra thick] (.495,.495) -- (.495,1);
				\draw[red, ultra thick] (.33,.495) -- (.33,1);
				\draw[red, ultra thick] (.165,.66) -- (.66,.66);
				\draw[red, ultra thick] (.165,.825) -- (.66,.825);
				\draw (0,0) -- (1,0) -- (1,1) -- (0,1) -- cycle;
				\draw (.33,0) -- (.33, 1);
				\draw (.66,0) -- (.66,1);
				\draw (0,.66) -- (1,.66);
				\draw (0,.33) -- (1,.33);
				\draw (0,.495) -- (.66,.495);
				\draw (0,.825) -- (.66,.825);
				\draw (.495,.33) -- (.495,1);
				\draw (.165,.33) -- (.165,1);
			\end{tikzpicture},
			\begin{tikzpicture}[scale=1.25]
				\filldraw[draw=red, ultra thick, fill=red!50] (0,.33) -- (.495,.33) -- (.495,.825) -- (0,.825) -- cycle;
				\draw[red, ultra thick] (.165,.33) -- (.165,.825);
				\draw[red, ultra thick] (.33,.33) -- (.33,.825);
				\draw[red, ultra thick] (0,.66) -- (.495,.66);
				\draw[red, ultra thick] (0,.495) -- (.495,.495);
				\draw (0,0) -- (1,0) -- (1,1) -- (0,1) -- cycle;
				\draw (.33,0) -- (.33, 1);
				\draw (.66,0) -- (.66,1);
				\draw (0,.66) -- (1,.66);
				\draw (0,.33) -- (1,.33);
				\draw (0,.495) -- (.66,.495);
				\draw (0,.825) -- (.66,.825);
				\draw (.495,.33) -- (.495,1);
				\draw (.165,.33) -- (.165,1);
			\end{tikzpicture},\\
			\begin{tikzpicture}[scale=1.25]
				\filldraw[draw=red, ultra thick, fill=red!50] (.165,.33) -- (.66,.33) -- (.66,.825) -- (.165,.825) -- cycle;
				\draw[red, ultra thick] (.495,.33) -- (.495,.825);
				\draw[red, ultra thick] (.33,.33) -- (.33,.825);
				\draw[red, ultra thick] (.165,.66) -- (.66,.66);
				\draw[red, ultra thick] (.165,.495) -- (.66,.495);
				\draw (0,0) -- (1,0) -- (1,1) -- (0,1) -- cycle;
				\draw (.33,0) -- (.33, 1);
				\draw (.66,0) -- (.66,1);
				\draw (0,.66) -- (1,.66);
				\draw (0,.33) -- (1,.33);
				\draw (0,.495) -- (.66,.495);
				\draw (0,.825) -- (.66,.825);
				\draw (.495,.33) -- (.495,1);
				\draw (.165,.33) -- (.165,1);
			\end{tikzpicture},
			\begin{tikzpicture}[scale=1.25]
				\filldraw[draw=red, ultra thick, fill=red!50] (0,0) -- (1,0) -- (1,1) -- (0,1) -- cycle;
				\draw[red, ultra thick] (.66,0) -- (.66,1);
				\draw[red, ultra thick] (.33,0) -- (.33,1);
				\draw[red, ultra thick] (0,.33) -- (1,.33);
				\draw[red, ultra thick] (0,.66) -- (1,.66);
				\draw (0,0) -- (1,0) -- (1,1) -- (0,1) -- cycle;
				\draw (.33,0) -- (.33, 1);
				\draw (.66,0) -- (.66,1);
				\draw (0,.66) -- (1,.66);
				\draw (0,.33) -- (1,.33);
				\draw (0,.495) -- (.66,.495);
				\draw (0,.825) -- (.66,.825);
				\draw (.495,.33) -- (.495,1);
				\draw (.165,.33) -- (.165,1);
			\end{tikzpicture}.
		\end{tabular}
	}\quad
	\subfloat[]{
		\begin{tabular}{c}
			\begin{tikzpicture}[scale=1.25]
				\filldraw[draw=cyan, ultra thick, fill=cyan!50] (0,.495) -- (.495,.495) -- (.495,1) -- (0,1) -- cycle;
				\draw[cyan, ultra thick] (.165,.495) -- (.165,1);
				\draw[cyan, ultra thick] (.33,.495) -- (.33,1);
				\draw[cyan, ultra thick] (0,.66) -- (.495,.66);
				\draw[cyan, ultra thick] (0,.825) -- (.495,.825);
				\draw (0,0) -- (1,0) -- (1,1) -- (0,1) -- cycle;
				\draw (.33,0) -- (.33, 1);
				\draw (.66,0) -- (.66,1);
				\draw (0,.66) -- (1,.66);
				\draw (0,.33) -- (1,.33);
				\draw (0,.495) -- (1,.495);
				\draw (0,.825) -- (1,.825);
				\draw (.495,.33) -- (.495,1);
				\draw (.165,.33) -- (.165,1);
			\end{tikzpicture},
			\begin{tikzpicture}[scale=1.25]
				\filldraw[draw=cyan, ultra thick, fill=cyan!50] (.165,.495) -- (.66,.495) -- (.66,1) -- (.165,1) -- cycle;
				\draw[cyan, ultra thick] (.495,.495) -- (.495,1);
				\draw[cyan, ultra thick] (.33,.495) -- (.33,1);
				\draw[cyan, ultra thick] (.165,.66) -- (.66,.66);
				\draw[cyan, ultra thick] (.165,.825) -- (.66,.825);
				\draw (0,0) -- (1,0) -- (1,1) -- (0,1) -- cycle;
				\draw (.33,0) -- (.33, 1);
				\draw (.66,0) -- (.66,1);
				\draw (0,.66) -- (1,.66);
				\draw (0,.33) -- (1,.33);
				\draw (0,.495) -- (1,.495);
				\draw (0,.825) -- (1,.825);
				\draw (.495,.33) -- (.495,1);
				\draw (.165,.33) -- (.165,1);
			\end{tikzpicture},
			\begin{tikzpicture}[scale=1.25]
				\filldraw[draw=cyan, ultra thick, fill=cyan!50] (0,.33) -- (.495,.33) -- (.495,.825) -- (0,.825) -- cycle;
				\draw[cyan, ultra thick] (.165,.33) -- (.165,.825);
				\draw[cyan, ultra thick] (.33,.33) -- (.33,.825);
				\draw[cyan, ultra thick] (0,.66) -- (.495,.66);
				\draw[cyan, ultra thick] (0,.495) -- (.495,.495);
				\draw (0,0) -- (1,0) -- (1,1) -- (0,1) -- cycle;
				\draw (.33,0) -- (.33, 1);
				\draw (.66,0) -- (.66,1);
				\draw (0,.66) -- (1,.66);
				\draw (0,.33) -- (1,.33);
				\draw (0,.495) -- (1,.495);
				\draw (0,.825) -- (1,.825);
				\draw (.495,.33) -- (.495,1);
				\draw (.165,.33) -- (.165,1);
			\end{tikzpicture},
			\begin{tikzpicture}[scale=1.25]
				\filldraw[draw=cyan, ultra thick, fill=cyan!50] (.165,.33) -- (.66,.33) -- (.66,.825) -- (.165,.825) -- cycle;
				\draw[cyan, ultra thick] (.495,.33) -- (.495,.825);
				\draw[cyan, ultra thick] (.33,.33) -- (.33,.825);
				\draw[cyan, ultra thick] (.165,.66) -- (.66,.66);
				\draw[cyan, ultra thick] (.165,.495) -- (.66,.495);
				\draw (0,0) -- (1,0) -- (1,1) -- (0,1) -- cycle;
				\draw (.33,0) -- (.33, 1);
				\draw (.66,0) -- (.66,1);
				\draw (0,.66) -- (1,.66);
				\draw (0,.33) -- (1,.33);
				\draw (0,.495) -- (1,.495);
				\draw (0,.825) -- (1,.825);
				\draw (.495,.33) -- (.495,1);
				\draw (.165,.33) -- (.165,1);
			\end{tikzpicture},\\
			\begin{tikzpicture}[scale=1.25]
				\filldraw[draw=cyan, ultra thick, fill=cyan!50] (0.333,0.495) -- (1,0.495) -- (1,1) -- (0.333,1) -- cycle;
				\draw[cyan, ultra thick] (.66,0.495) -- (.66,1);
				\draw[cyan, ultra thick] (0.495,0.495) -- (0.495,1);
				\draw[cyan, ultra thick] (0.333, 0.825) -- (1,0.825);
				\draw[cyan, ultra thick] (0.333,.66) -- (1,.66);
				\draw (0,0) -- (1,0) -- (1,1) -- (0,1) -- cycle;
				\draw (.33,0) -- (.33, 1);
				\draw (.66,0) -- (.66,1);
				\draw (0,.66) -- (1,.66);
				\draw (0,.33) -- (1,.33);
				\draw (0,.495) -- (1,.495);
				\draw (0,.825) -- (1,.825);
				\draw (.495,.33) -- (.495,1);
				\draw (.165,.33) -- (.165,1);
			\end{tikzpicture},
			\begin{tikzpicture}[scale=1.25]
				\filldraw[draw=cyan, ultra thick, fill=cyan!50] (0.333,0.333) -- (1,0.333) -- (1,0.825) -- (0.333,0.825) -- cycle;
				\draw[cyan, ultra thick] (.66,0.495) -- (.66,0.825);
				\draw[cyan, ultra thick] (0.495,0.495) -- (0.495,0.825);
				\draw[cyan, ultra thick] (0.333, 0.495) -- (1,0.495);
				\draw[cyan, ultra thick] (0.333,.66) -- (1,.66);
				\draw (0,0) -- (1,0) -- (1,1) -- (0,1) -- cycle;
				\draw (.33,0) -- (.33, 1);
				\draw (.66,0) -- (.66,1);
				\draw (0,.66) -- (1,.66);
				\draw (0,.33) -- (1,.33);
				\draw (0,.495) -- (1,.495);
				\draw (0,.825) -- (1,.825);
				\draw (.495,.33) -- (.495,1);
				\draw (.165,.33) -- (.165,1);
			\end{tikzpicture},
			\begin{tikzpicture}[scale=1.25]
				\filldraw[draw=cyan, ultra thick, fill=cyan!50] (0,0) -- (1,0) -- (1,0.666) -- (0,0.666) -- cycle;
				\draw[cyan, ultra thick] (.333,0) -- (.333,0.666);
				\draw[cyan, ultra thick] (0.666,0.666) -- (0.666,0.666);
				\draw[cyan, ultra thick] (0, 0.495) -- (1,0.495);
				\draw[cyan, ultra thick] (0,0.333) -- (1,0.333);
				\draw (0,0) -- (1,0) -- (1,1) -- (0,1) -- cycle;
				\draw (.33,0) -- (.33, 1);
				\draw (.66,0) -- (.66,1);
				\draw (0,.66) -- (1,.66);
				\draw (0,.33) -- (1,.33);
				\draw (0,.495) -- (1,.495);
				\draw (0,.825) -- (1,.825);
				\draw (.495,.33) -- (.495,1);
				\draw (.165,.33) -- (.165,1);
			\end{tikzpicture}.
		\end{tabular}
	}
	\caption{One iteration of the N$_2$S$_2$ refinement for bidegree $\m{p} = \pmb{2}$. A given LR B-spline is marked for refinement in (a). Its support and associated tensor mesh are highlighted. The first step of the N$_2$S$_2$ strategy applies a structured mesh refinement, shown in (b). The second step restores the N$_2$S property by suitably extending meshlines in one coordinate direction to eliminate nesting among the LR B-splines. The resulting mesh, obtained by extending the horizontal meshlines, is shown in (c), which is also the final mesh produced by the N$_2$S$_2$ refinement. Four LR B-splines in the lower-right corner of the orange region of the mesh in (b) are nested in the LR B-spline with larger support overlapping the corner, shown in (d). The extensions of (c) restore the N$_2$S property by replacing the B-spline with the largest support in (d) by the three shown in the bottom row of (e).}\label{fig:N2S2iteration}
\end{figure}
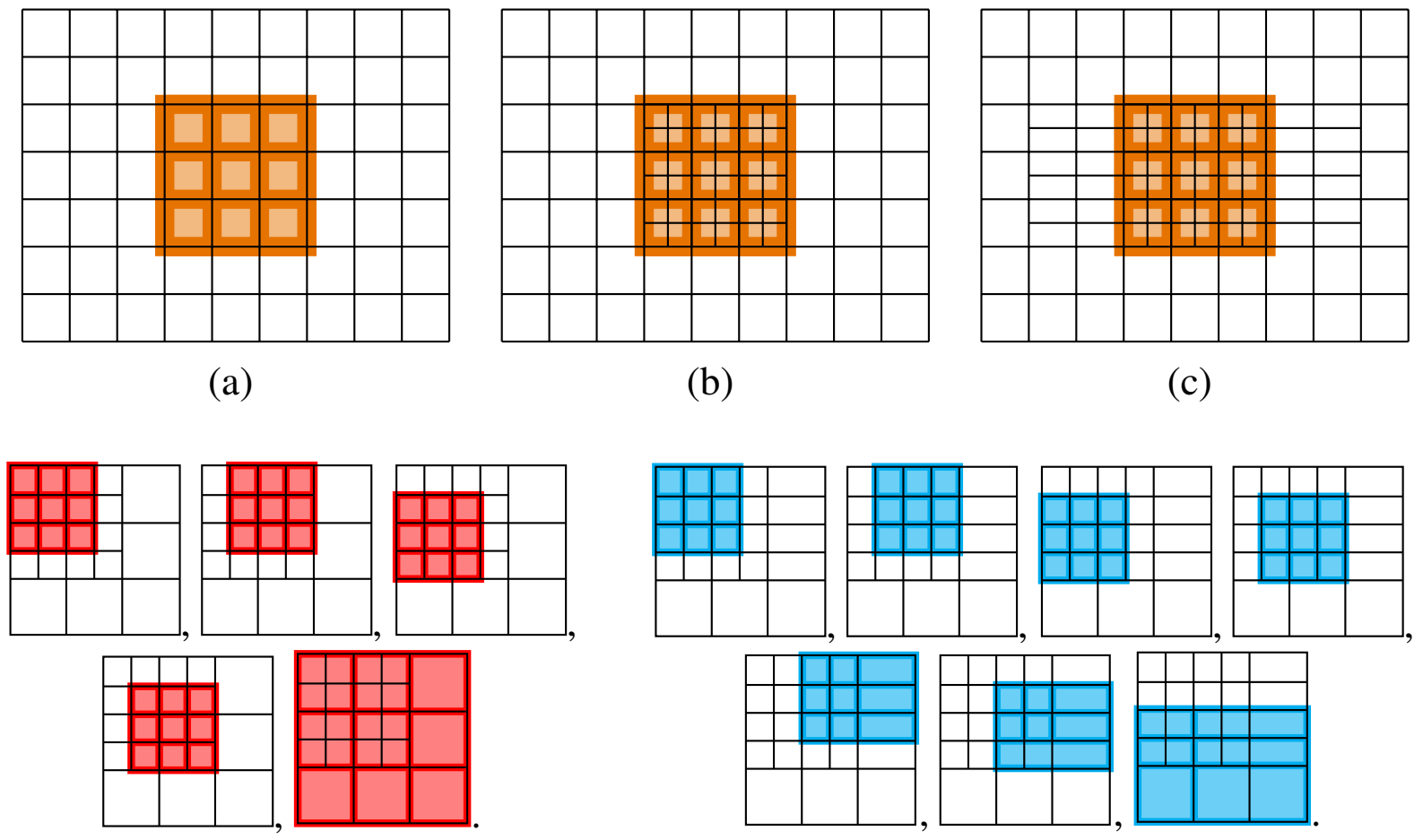

The LR B-splines newly generated in step 2 may themselves become nested in other LR B-splines, thereby propagating the N$_2$S restoration process. However, this propagation is proven to terminate \cite[Theorem~3.5]{n2s2}. From empirical observations, it typically terminates after a few additional refinements.
Figure~\ref{fig:N2S2iteration} illustrates one iteration of the N$_2$S$_2$ strategy.

There also exist other refinement procedures that guarantee the N$_2$S property. We briefly review them in the following.
\begin{itemize}
	\item \textbf{Hierarchical refinement} \cite{hlr}: It guarantees the N$_2$S property and also ensures \emph{shape regularity} and \emph{local quasi-uniformity} (or \emph{grading}) of the LR mesh. These are useful properties for achieving theoretically optimal convergence rates in the numerical discretization of PDEs \cite{igaf,axioms}. Roughly speaking, shape regularity prevents mesh boxes from becoming excessively stretched, while grading ensures that neighbouring elements have comparable sizes.
	Completeness of the resulting set of LR B-splines is an additional important feature of this refinement process.
	The hierarchical construction, however, is not inherently adaptive: once the mesh has been generated, modifying the refinement pattern requires rebuilding it from scratch to preserve the N$_2$S property, shape regularity, and grading.
	\item \textbf{Effective grading} (\textbf{EG}) \textbf{refinement} \cite{eg}: It overcomes the limitation of the hierarchical construction, while retaining its pleasant structural mesh properties. The meshes produced by EG are often very similar to those obtained through the hierarchical construction but the EG procedure can be adaptively updated without reconstructing the entire mesh. Completeness of the resulting set of LR B-splines is also guaranteed.
	\item \textbf{Reachable minimal} (\textbf{RM}) \textbf{B-splines} \cite{rmb,lrrmb}: RM B-splines are a special instance of LR B-splines with reduced global smoothness. Their associated LR mesh is adjusted to satisfy the N$_2$S property through a process similar to that employed in the N$_2$S$_2$ strategy.
\end{itemize}

An extensive numerical comparison of different refinement strategies for LR B-splines can be found in \cite{comparison}.

\section{Conclusion}\label{sec:conclusion}
In this paper, we reviewed the theoretical foundations of LR B-splines, one of the most elegant extensions of univariate B-splines to local tensor structures. They reduce to standard B-splines in the univariate setting, and their formulation remains similar to that of standard tensor-product B-splines, while allowing for local refinement.
Several properties, such as local support and non-negativity, are directly inherited from tensor-product B-splines. Partition of unity can be obtained by a proper scaling procedure.
As the dimension of spline spaces on LR meshes is always stable, we also discussed the completeness of the span of the constructed LR B-spline set. Moreover, we looked into linear dependence and local linear independence, and automatic adaptive refinement strategies ensuring the latter property.

As mentioned in the introduction, there are various other spline technologies based on local tensor structures that allow for local refinement. The most popular examples are (analysis-suitable) T-splines \cite{ast,sederberg}, (truncated) hierarchical B-splines \cite{hb,thb}, and PHT-splines \cite{pht}. Numerical comparisons between LR B-splines and some of these alternatives have been presented in \cite{johannessen2015,stangeby}. When applying suitable refinement strategies \cite{eg,n2s2},
LR B-splines do not require any restriction on degree and smoothness, are locally linearly independent, and ensure good mesh quality.
Under mild assumptions, analysis-suitable T-splines enjoy the first two properties \cite{ast}, but their structure can negatively affect the localization of the refinement process.
PHT-splines require reduced smoothness, while a characterization for local linear independence is missing for (truncated) hierarchical B-splines. Nevertheless, linear independence is always guaranteed for the latter.

Regarding flexibility, when no constraints are imposed for the N$_2$S property nor for specific mesh properties such as grading, LR B-splines allow for very general mesh configurations as only the LR rules must be satisfied. T-splines likewise accommodate a wide range of mesh configurations. In contrast, the refinement procedure for hierarchical spline spaces imposes a structure on the underlying mesh, thereby restricting the range of admissible refinements. Similarly, the constraints required to ensure the local linear independence of analysis-suitable T-splines limit the admissible mesh configurations. As a consequence, performing anisotropic refinement is more complicated for both (truncated) hierarchical B-splines and analysis-suitable T-splines. Instead, it can be readily achieved for LR B-splines, even when the local linear independence is enforced through the N$_2$S$_2$ strategy, although some propagation of refinement may be required.

In summary, LR B-splines emerge as an appealing alternative to other known spline methods based on local tensor structures. Thanks to their flexibility in supporting anisotropic local refinement, their application in data approximation and numerical simulation is receiving increasing attention; see, e.g., \cite{scattereddata,hollweck,johannessen2014,occelli,n2s2,skytt}.
Additionally, recent extensions of the LR paradigm beyond the polynomial setting \cite{tmesh,raval} paves the path to further challenging applications.

\begin{acknowledgement}
The authors are members of the research group GNCS (Gruppo Nazionale per il Calcolo Scientifico) of INdAM (Istituto Nazionale di Alta Matematica).
This work has been supported by the MUR Excellence Department Project MatMod@TOV (CUP E83C23000330006) awarded to the Department of Mathematics of the University of Rome Tor Vergata, by the Department of Mathematics of the University of Rome Tor Vergata through the Project METRO (CUP  E83C25000630005), by a Project of Relevant
National Interest (PRIN) under the National Recovery and Resilience Plan (PNRR) funded by the European Union -- Next Generation EU (CUP E53D23017910001), by the Italian Research Center on High Performance Computing, Big Data and Quantum Computing (CUP E83C22003230001), and by several GNCS projects.
\end{acknowledgement}

\bibliographystyle{plain}
\bibliography{biblio}
\end{document}